\documentclass{elsarticle}

\newcommand{\dcom}[1]{\textcolor{black}{#1}}

\biboptions{sort&compress}

\usepackage[utf8]{inputenc}
\usepackage{graphicx,amssymb,amstext,amsmath}
\usepackage{amsmath,amsxtra,amssymb,latexsym, amscd,amsthm}
\usepackage{indentfirst}
\usepackage[mathscr]{eucal}
\usepackage{graphicx}
\usepackage{enumitem}
\usepackage{ulem}
\usepackage{url}
\usepackage{subfig}
\usepackage{bm}
\usepackage{multirow}
\usepackage{algorithm2e}
\RestyleAlgo{boxed}
\RestyleAlgo{boxruled}
\usepackage{tikz}
\usetikzlibrary{matrix,shapes,arrows,positioning,chains}
\tikzset{
block/.style={
    rectangle,
    draw,
    text width=17em,
    text centered,
},
decision/.style={
    rectangle,
    draw,
    text width=17em,
    text centered,
    rounded corners
},
cloud/.style={
    draw,
    ellipse,
    minimum height=2em
},
descr/.style={
    fill=white,
    inner sep=2.5pt
},
connector/.style={
    -latex,
    font=\scriptsize
},
rectangle connector/.style={
    connector,
    to path={(\tikztostart) -- ++(#1,0pt) \tikztonodes |- (\tikztotarget) },
    pos=0.5
},
rectangle connector/.default=-2cm,
straight connector/.style={
    connector,
    to path=--(\tikztotarget) \tikztonodes
}
}

\usepackage{float}

\usepackage{mathrsfs}
\usepackage{geometry}
\usepackage{fancyhdr}

\usepackage{listings}             
\lstdefinestyle{myCustomMatlabStyle}{
  language=Matlab,
  numbers=left,
  stepnumber=1,
  numbersep=10pt,
  tabsize=4,
  showspaces=false,
  showstringspaces=false
}
\newcommand{\bx}{\ensuremath{{\bm{x}}}}

\newcommand{\bg}{\ensuremath{{\bm{g}}}}
\newcommand{\dunit}{\,\ensuremath{\text{mm}^2/\text{s}}}
\newcommand{\bunit}{\,\ensuremath{\text{s/mm}^2}}
\newcommand{\tunit}{\ensuremath{\text{ms}}}
\newcommand{\kunit}{\,\ensuremath{\text{m}/\text{s}}}

\newcommand{\lunit}{\ensuremath{\mu\text{m}}}

\newcommand{\vthree}[3]{\ensuremath{\begin{bmatrix}#1\\#2\\#3\end{bmatrix}}}

\newcommand{\bn}{\ensuremath{{\bm{n}}}}

\newcommand{\ben}{\begin{equation*}}
\newcommand{\een}{\end{equation*}}
\newcommand{\be} [1] {\begin{equation} \label{#1}}
\newcommand{\ee}{\end{equation}}

\newcommand{\ignore}[1]{}

\newcommand{\bug}{\ensuremath{{\bm{u_g}}}}
\newcommand{\e}[1]{\ensuremath{\times 10^{#1}}}

\newcommand{\soutnew}[2]{#2}

\newcommand{\marginparnew}[1]{}

\graphicspath{ {Fig/} }

\begin{document}

\begin{frontmatter}



\title{SpinDoctor: a MATLAB toolbox for diffusion MRI simulation}


\hyphenation{appro-xi-ma-tion}

\author[inria]{Jing-Rebecca Li \corref{cor1}}
\ead{jingrebecca.li@inria.fr}
\author[kth]{Van-Dang Nguyen}
\author[inria]{Try Nguyen Tran}
\author[usb]{Jan Valdman}
\author[inria]{Cong-Bang Trang}
\author[inria]{Khieu Van Nguyen}
\author[inria]{Duc Thach Son Vu}
\author[inria]{Hoang An Tran}
\author[inria]{Hoang Trong An Tran}
\author[inria]{Thi Minh Phuong Nguyen}
\address[inria]{INRIA Saclay, Equipe DEFI, CMAP, Ecole Polytechnique,
Route de Saclay, 91128 Palaiseau Cedex, France}

\address[kth]{Department of Computational Science and Technology, KTH Royal Institute of Technology, Sweden}
\address[usb]{Institute of Mathematics, Faculty of Science, University of South Bohemia, \v{C}esk\'e Bud\v{e}jovice and Institute of Information Theory and Automation of the ASCR, Prague, Czech Republic}
\cortext[cor1]{Corresponding author}

\begin{abstract}  
\soutnew{
The complex transverse water proton magnetization subject to diffusion-encoding magnetic
field gradient pulses in a heterogeneous medium can be modeled by the
multiple compartment Bloch-Torrey partial differential equation (BTPDE).
A mathematical model for the time-dependent apparent diffusion coefficient (ADC),
called the H-$ADC$ model, was obtained recently  using homogenization techniques on the BTPDE.
Under the assumption of negligible water exchange between compartments, 
the H-$ADC$ model produces the ADC of a diffusion medium 
from the solution of a diffusion equation (DE) 
subject to a time-dependent Neumann boundary condition. }{}

\soutnew{
This paper describes a publicly available MATLAB toolbox called SpinDoctor that can be used 1) to solve the BTPDE to obtain
the dMRI signal (the toolbox provides a way of robustly fitting the dMRI signal to obtain the fitted ADC); 
2) to solve the DE of the H-$ADC$ model to obtain the ADC;
3) a short-time approximation formula for the ADC is also included in the toolbox 
for comparison with the simulated ADC.
The PDEs are solved by $P1$ finite elements combined with
build-in MATLAB routines for solving ordinary differential equations.
The finite element mesh generation is performed using an external package called Tetgen that is included in the toolbox.}{}

\soutnew{SpinDoctor provides built-in options of including 1) spherical cells with a nucleus; 
2) cylindrical cells with a myelin layer; 3) an extra-cellular space (ECS) enclosed either 
a) in a box or b) in a tight wrapping around the cells; 
4) deformation of canonical cells by bending and twisting.  
5) permeable membranes for the BT-PDE (the H-$ADC$ assumes negligible permeability).
Built-in diffusion-encoding pulse sequences include the Pulsed Gradient Spin Echo and the Oscillating Gradient Spin Echo. }{}   

\soutnew{In this paper, we describe how to use the SpinDoctor and illustrate with numerical examples of computing the
ADC using the BTPDE, the HADC, and the STA models in spherical cells, cylindrical cells, the myelin layer, 
and the ECS.}{}

\soutnew{}{
The complex transverse water proton magnetization subject to diffusion-encoding magnetic
field gradient pulses in a heterogeneous medium can be modeled by the
multiple compartment Bloch-Torrey partial differential equation.
Under the assumption of negligible water exchange between compartments, 
the time-dependent apparent diffusion coefficient can be directly computed 
from the solution of a diffusion equation subject to a time-dependent Neumann boundary condition. }

\soutnew{}{
This paper describes a publicly available MATLAB toolbox called SpinDoctor that can be used 
1) to solve the Bloch-Torrey partial differential equation in order to simulate the diffusion magnetic resonance imaging signal; 
2) to solve a diffusion partial differential equation to obtain directly the apparent diffusion coefficient;
3) to compare the simulated apparent diffusion coefficient with a short-time approximation formula. }

\soutnew{}{The partial differential equations are solved by $P1$ finite elements combined with
built-in MATLAB routines for solving ordinary differential equations.
The finite element mesh generation is performed using an external package called Tetgen.}

\soutnew{}{SpinDoctor provides built-in options of including 1) spherical cells with a nucleus; 
2) cylindrical cells with a myelin layer; 3) an extra-cellular space enclosed either 
a) in a box or b) in a tight wrapping around the cells; 
4) deformation of canonical cells by bending and twisting;
5) permeable membranes;
Built-in diffusion-encoding pulse sequences include the Pulsed Gradient Spin Echo and the Oscillating Gradient Spin Echo. }   

\soutnew{}{We describe in detail how to use the SpinDoctor toolbox.  We validate SpinDoctor simulations 
using reference signals computed by the Matrix Formalism method.  We compare the accuracy and computational
time of SpinDoctor simulations with Monte-Carlo simulations and show significant speed-up of SpinDoctor 
over Monte-Carlo simulations in complex geometries.  We also illustrate several extensions 
of SpinDoctor functionalities, including the incorporation of $T_2$ relaxation, the simulation of
non-standard diffusion-encoding sequences,
as well as the use of externally generated geometrical meshes.}

\end{abstract}


\end{frontmatter}


\section{Introduction}
\label{introduction}

Diffusion magnetic resonance imaging 
\soutnew{(dMRI) }{}is an imaging modality that can be used 
to probe the tissue micro-structure by encoding the incohorent motion of water molecules 
with magnetic field gradient pulses.  This motion during 
the diffusion-encoding time causes a signal attenuation 
from which the apparent diffusion coefficient \soutnew{(ADC)}{}, 
(and possibly higher order diffusion terms,
can be calculated \cite{Hahn1950,Stejskal1965,Bihan1986}.

For unrestricted diffusion, the root of the mean squared displacement of molecules 
is given by $\bar{x} = \sqrt{2\,dim \,\sigma_0 t}$,
where $dim$ is the spatial dimension, 
$\sigma_0$ is the intrinsic diffusion coefficient, and $t$ is the diffusion time. 
In biological tissue, the diffusion is usually hindered or restricted (for example, by cell membranes) 
and the mean square displacement is smaller than in the case of 
unrestricted diffusion.  This deviation from unrestricted diffusion can be used to infer information
about the tissue micro-structure.  The 
\soutnew{dMRI}{} 
experimental parameters that can be varied include 
\begin{enumerate}
\item the diffusion time (one can choose the parameters of the diffusion-encoding sequence, such as 
Pulsed Gradient Spin Echo 
\soutnew{(PGSE)}{}\cite{Stejskal1965} and Oscillating Gradient
\soutnew{(OGSE)}{}\cite{Does2003}).
\item the magnitude of the diffusion-encoding gradient (when the 
\soutnew{MRI}{magnetic resonance imaging} 
signal is acquired at low gradient magnitudes, the signal contains only information about 
the 
\soutnew{ADC}{apparent diffusion coefficient}, 
at higher values, Kurtosis imaging \cite{Jensen2005} becomes possible);
\item the direction of the diffusion-encoding gradient (many directions may be probed, as in 
\soutnew{HARDI}{high angular resolution diffusion imaging}
\cite{Tuch2002}).
\end{enumerate}

Using 
\soutnew{dMRI}{diffusion magnetic resonance imaging}
to get tissue structural information in the mamalian brain has been the focus of much experimental and modeling work in recent years 
\cite{Assaf2008,Alexander2010,Zhang2011,Zhang2012,Burcaw2015,Palombo2017a,Palombo2016,Ning2017}.  
The predominant approach up to now has been adding the 
\soutnew{dMRI}{diffusion magnetic resonance imaging} signal from simple geometrical components and extracting model parameters of interest.   Numerous biophysical models subdivide the tissue into compartments described by spheres, ellipsoids, cylinders, and the extra-cellular space 
\soutnew{(ECS)}{} 
\cite{McHugh2015,Reynaud2017,Assaf2008,Alexander2010,Zhang2011,Burcaw2015,Fieremans2011,Panagiotaki2012,Jespersen2007,Palombo2017a}.
Some model parameters of interest include axon diameter and orientation, neurite density,  dendrite structure, 
the volume fraction and size distribution of cylinder and sphere components and the effective diffusion coefficient or tensor of the 
\soutnew{(ECS)}{extra-cellular space}.

Numerical simulations can help deepen the understanding of the relationship between 
the cellular structure and the 
\soutnew{dMRI}{diffusion magnetic resonance imaging} 
signal and lead to the formulation of appropriate models.  
They can be also used to investigate the effect of different pulse sequences and tissue features on the measured signal which can be used for the development, testing, and optimization of novel 
\soutnew{MRI}{diffusion magnetic resonance imaging}
pulse sequences \cite{Ianus2016, Drobnjak2011,Mercredi2018,Rensonnet2018}.

Two main groups of approaches to the numerical simulation of 
\soutnew{dMRI}{diffusion magnetic resonance imaging} 
are 1) using random walkers to mimic the diffusion process in a geometrical configuration; 2)  solving the Bloch-Torrey PDE 
\soutnew{(BTPDE)}{}, which describes the evolution of the complex transverse water proton magnetization under the influence of diffusion-encoding magnetic field gradients pulses.  

The first group is referred to as Monte-Carlo simulations in the literature and previous works include \cite{Hughes1995, Yeh2013, Hall2009,Palombo2016,Balls2009}.  A GPU-based acceleration of Monte-Carlo simulation was proposed in \cite{Nguyen2018a,Waudby2011}. 
Some software packages using this approach include 
\begin{enumerate}
\item Camino Diffusion MRI Toolkit developed at UCL (http://cmic.cs.ucl.ac.uk/camino/);
\item DIFSIM developed at UC San Diego (http://csci.ucsd.edu/projects/simulation.html);
\item Diffusion Microscopist Simulator \cite{Yeh2013} developed at Neurospin, CEA.
\end{enumerate}

The second group relies on solving the 
\soutnew{(BTPDE)}{Bloch-Torrey PDE} 
in a geometrical configuration.
In \cite{Hagslatt2003, Loren2005, Moroney2013} a simplifying assumption 
called the narrow pulse approximation was used, where
the pulse duration was assumed to be much smaller 
than the delay between pulses.  This assumption allows the solution of 
the diffusion equation 
\soutnew{(DE)} instead of the more complicated 
\soutnew{(BTPDE)}{Bloch-Torrey PDE}.
More generally, numerical methods to solve the 
\soutnew{(BTPDE)}{Bloch-Torrey PDE}.
with arbitrary temporal profiles have been proposed in \cite{Xu2007, Li2014, Nguyen2014, Beltrachini2015}. 
The computational domain is discretized either by a Cartesian grid \cite{Xu2007, Russell2012, Li2014} or finite elements 
\cite{Hagslatt2003, Loren2005, Moroney2013, Nguyen2014, Beltrachini2015}.
The unstructured mesh of a finite element discretization appeared to be better than a Cartesian grid in both geometry description and signal approximation \cite{Nguyen2014}.
For time discretization, both explicit and implicit methods have been used. 
In \cite{Moroney2013} a second order implicit time-stepping method called the generalized $\alpha-$method was used to allow for high frequency energy dissipation. 
An adaptive explicit Runge-Kutta Chebyshev 
\soutnew{(RKC)} method of second order was used in \cite{Li2014, Nguyen2014}. 
It has been theoretically proven that the 
\soutnew{RKC}{Runge-Kutta Chebyshev} method allows for a much larger time-step compared to the standard explicit Euler method \cite{Verwer1990}.  
There is an example showing that the 
\soutnew{RKC}{Runge-Kutta Chebyshev}
method is faster than the implicit Euler method in \cite{Nguyen2014}.  
The Crank-Nicolson method 
\soutnew{(CN)} was used in \cite{Beltrachini2015} to also allow for second order convergence in time. The efficiency of 
\soutnew{dMRI}{diffusion magnetic resonance imaging} 
simulations is also improved by either a high-performance FEM computing framework \cite{Nguyen2016a, Nguyen2018} for large-scale 
\soutnew{dMRI} simulations on supercomputers or a discretization on manifolds for thin-layer and thin-tube media \cite{Nguyen2019}.

In this paper, we present a MATLAB Toolbox called SpinDoctor that is a simulation pipeline going from 
the definition of a geometrical configuration through the numerical solution of the 
\soutnew{(BTPDE)}{Bloch-Torrey PDE} 
to the fitting of the 
\soutnew{ADC}{apparent diffusion coefficient}
from the simulated signal.  It also includes two other modules for calculating the 
\soutnew{ADC}{apparent diffusion coefficient}.
The first module is a 
homogenized apparent diffusion coefficient mathematical model, which was obtained recently using homogenization techniques on the 
\soutnew{(BTPDE)}{Bloch-Torrey PDE}.
In the 
\soutnew{H-$ADC$}{homogenized} model, the 
\soutnew{ADC}{apparent diffusion coefficient} 
of a geometrical configuration can be 
computed after solving a diffusion equation 
\soutnew{(DE) }{}
subject to a time-dependent Neumann boundary condition, under the 
assumption of negligible water exchange between compartments. 
The second module computes the short time approximation 
\soutnew{(STA)} formula for the 
\soutnew{ADC}{apparent diffusion coefficient}.
The 
\soutnew{(STA)}{short time approximation}
implemented in SpinDoctor includes a recent generalization 
of this formula to account for finite pulse duration in the 
\soutnew{PGSE}{pulsed gradient spin echo}. 
Both of these two 
\soutnew{ADC}{apparent diffusion coefficient}
calculations are sensitive to the diffusion-encoding gradient direction, unlike
many previous works where the anisotropy 
\soutnew{of the ADC}{} is neglected in analytical model development. 
\soutnew{The toolbox is publicly available at:}{}
\begin{center}
\soutnew{    \url{https://github.com/jingrebeccali/SpinDoctor} }{}
\end{center}

In summary, SpinDoctor 
\soutnew{is a MATLAB Toolbox that}{}
\begin{enumerate}
\item solves the 
\soutnew{(BTPDE)}{Bloch-Torrey PDE} in three dimensions to obtain
the 
\soutnew{dMRI}{diffusion magnetic resonance imaging} 
signal;
\item robustly fits the
\soutnew{dMRI}{diffusion magnetic resonance imaging} 
signal to obtain the 
\soutnew{ADC}{apparent diffusion coefficient}; 
\item solves the 
\soutnew{HADC}{homogenized apparent diffusion coefficient} 
model in three dimensions to obtain the 
\soutnew{ADC}{apparent diffusion coefficient};
\item computes the short-time approximation 
\soutnew{(STA)} of the 
\soutnew{ADC}{apparent diffusion coefficient};
\item computes useful geometrical quantities such as the compartment volumes and surface areas; 
\item allows permeable membranes for the 
\soutnew{(BTPDE)}{Bloch-Torrey PDE}
(the 
\soutnew{HADC}{homogenized apparent diffusion coefficient} 
assumes negligible permeabilty).
\item displays the gradient-direction dependent 
\soutnew{ADC}{apparent diffusion coefficient}; 
in three dimensions using spherical harmonics interpolation;
\end{enumerate}

SpinDoctor provides the following built-in functionalities:
\begin{enumerate}
\item  placement of non-overlapping spherical cells (with an optional nucleus) of different radii close to each other;
\item  placement of non-overlapping cylindrical cells (with an optional myelin layer) of different radii close to each other in a canonical configuration where they are parallel to the $z$-axis;
\item inclusion of an extra-cellular space 
\soutnew{(ECS)}{}  that is enclosed either 
\begin{enumerate}
\item in a tight wrapping around the cells; or
\item in a rectangular box; 
\end{enumerate}
\item  deformation of the canonical configuration by bending and twisting;  

\end{enumerate}

Built-in diffusion-encoding pulse sequences include 
\begin{enumerate}
\item the Pulsed Gradient Spin Echo 
\soutnew{(PGSE)}{};
\item the Oscillating Gradient Spin Echo \soutnew{(cos-OGSE, sin-OGSE)}{(cos- and sin- type gradients)}.  
\end{enumerate}

SpinDoctor uses the following methods:
\begin{enumerate}
\item it generates a good quality surface triangulation of the user specified geometrical configuration by calling built-in MATLAB computational geometry functions;
\item it creates a good quality tetrehedra finite elements mesh from the above surface triangulation by calling Tetgen \cite{Si2015}, an external package (executable files are included in the Toolbox package); 
\item it constructs finite element matrices for linear finite elements on tetrahedra (P1) using routines from \cite{RahmanValdman2013};
\item it adds additional degrees of freedom on the compartment interfaces to allow permeability conditions for the
\soutnew{(BTPDE)}{Bloch-Torrey PDE} 
using the formalism in \cite{Nguyen2014d};
\item it solves the semi-discretized FEM equations by calling built-in MATLAB routines for solving ordinary differential equations \soutnew{(ODEs)}{}.
\end{enumerate}

\soutnew{}{The SpinDoctor toolbox has been developed in the MATLAB R2017b and requires no additional MATLAB toolboxes. }The toolbox is publicly available at:
\begin{center}
    \url{https://github.com/jingrebeccali/SpinDoctor} 
\end{center}

\section*{Abbreviations frequently used in the text}
\marginparnew{Added in the revised version}

 MRI -  magnetic resonance imaging 
  
dMRI -   diffusion magnetic resonance imaging 
  
ADC - apparent diffusion coefficient 

HADC - homogenized ADC
  
PGSE -  pulsed gradient spin echo 
  
OGSE - oscillating gradient 
   
ECS - extra-cellular space  
   
BTPDE - Bloch-Torrey partial differential equation 

PDE - partial differential equation 

ODE - ordinary differential equation 
  
HARDI - high angular resolution diffusion imaging 

STA - short time approximation 

FE - finite elements



\section{Theory}

Suppose the user would like to simulate a geometrical configuration of cells  
with an optional myelin layer or a nucleus.  If spins will be leaving the cells or if the user wants to simulate
the extra-cellular space (ECS), then the ECS will enclose the geometrical shapes.
Let $\Omega^e$ be the ECS, $\Omega^{in}_i$ the nucleus (or the axon)
and $\Omega^{out}_i$ the cytoplasm (or the myelin layer) of the $i$th cell.  
We denote the interface between $\Omega^{in}_i$ and $\Omega^{out}_i$ by $\Gamma_{i}$
and the interface between $\Omega^{out}_i$ and $\Omega^{e}$ by $\Sigma_{i}$, finally
the outside boundary of the ECS by $\Psi$.

\subsection{Bloch-Torrey PDE}

\label{PDEsofdiffusionMRI}

In diffusion MRI, a time-varying magnetic field gradient is
applied to the tissue to encode water diffusion.  
Denoting the effective 
time profile of the diffusion-encoding magnetic field gradient by $f(t)$, and letting the
vector $\bg$ contain the amplitude and direction information of the
magnetic field gradient, the complex transverse water proton magnetization 
in the rotating frame satisfies the Bloch-Torrey PDE:
\begin{alignat}{4}
\label{eq:btpde}
&\frac{\partial}{\partial t}{M^{in}_i(\bx,t)} &&= -I\gamma f(t) \bg \cdot \bx \,M^{in}_i(\bx,t) 
+ \nabla \cdot (\sigma^{in} \nabla M^{in}_i(\bx,t)), &&\bx \in \Omega^{in}_i,&\\
&\frac{\partial}{\partial t}{M^{out}_i(\bx,t)} &&= -I\gamma f(t) \bg \cdot \bx \,M^{out}_i(\bx,t) 
+ \nabla \cdot (\sigma^{out} \nabla M^{out}_i(\bx,t)),  \;&&\bx \in \Omega^{out}_i,&\\
&\frac{\partial}{\partial t}{M^e(\bx,t)} &&= -I\gamma f(t) \bg \cdot \bx \,M^e(\bx,t) 
+ \nabla \cdot (\sigma^{e} \nabla M^e(\bx,t)), & &\bx \in \Omega^{e},&
\end{alignat}
where $\gamma=2.67513\times 10^8\,\rm rad\,s^{-1}T^{-1}$ is the
gyromagnetic ratio of the water proton, $I$ is the imaginary unit,
$\sigma^l$ is the intrinsic diffusion coefficient in the compartment $\Omega^l_i$. 
The magnetization is a function of position $\bx$ and time $t$, 
and depends on the diffusion gradient vector $\bg$ and the time profile $f(t)$.
We denote the restriction of the magnetization 
in $\Omega^{in}_i$ by $M^{in}_i$, and similarly for $M^{out}_i$ and $M^{e}$.

Some commonly used time profiles (diffusion-encoding sequences) are:
\begin{enumerate}
\item 
The pulsed-gradient spin echo (PGSE) \cite{Stejskal1965} sequence,
with two rectangular pulses of duration $\delta$, separated by a time
interval $\Delta - \delta$, for which the profile $f(t)$ is
\be{eq:pgse}
f(t) =
\begin{cases}
1, \quad &t_1 \leq t \leq t_1+\delta, \\
-1,
\quad & t_1+\Delta < t \leq t_1+\Delta+\delta,\\
0, \quad & \text{otherwise,}
\end{cases}
\ee
where $t_1$ is the starting time of the first gradient pulse with $t_1
+ \Delta >T_E/2$, $T_E$ is the echo time at which the signal is measured.
\item 
The oscillating gradient spin echo (OGSE) sequence \cite{Callaghan1995,Does2003} 
was introduced to reach short diffusion times.  An OGSE
sequence usually consists of two oscillating pulses of duration
$T$, each containing $n$ periods, hence the 
frequency is $\omega = n\frac{2\pi}{T}$, separated by a time interval 
$\tau-T$.
For a cosine OGSE, the profile $f(t)$ is
\be{eq:ogse}
f(t) =
\begin{cases}
\cos{(n \frac{2\pi}{T} t)}, \quad &t_1 < t \leq t_1 + T, \\
-\cos{(n\frac{2\pi}{T} (t-\tau))},
\quad &\tau+t_1 < t \leq t_1+\tau + T,\\
0, \quad & \text{otherwise},
\end{cases}
\ee
where $\tau = T_E/2$.
\end{enumerate}

The BTPDE needs to be supplemented by interface conditions.
We recall the interface between $\Omega^{in}_i$ and $\Omega^{out}_i$ is $\Gamma_{i}$,
the interface between $\Omega^{out}_i$ and $\Omega^{e}$ is $\Sigma_{i}$, and
the outside boundary of the ECS is $\Psi$.
The two interface conditions on $\Gamma_i$ are the flux continuity
and a condition that incorporates a permeability coefficient
$\kappa^{in,out}$ across $\Gamma_i$:
:
\begin{alignat*}{3}
\sigma^{in} \nabla M^{in}_i(\bx,t) \cdot \bn^{in}_i 
&= -\sigma^{out}\nabla M^{out}_i(\bx,t) \cdot \bn^{out}_i,\quad &&\bx \in \Gamma_i,&\\
 \sigma^{in}\nabla M^{in}_i(\bx,t) \cdot \bn^{in}_i
 &= \kappa^{in,out} \Bigl(M^{out}_i(\bx,t)- M^{in}_i(\bx,t)\Bigl), \quad &&\bx \in \Gamma_i,&
\end{alignat*}
where $\bn$ is the unit outward pointing normal vector.
Similarly, between $\Omega^{out}_i$ and $\Omega^e$ we have
\begin{alignat*}{4}
&\sigma^{out} \nabla M^{out}_i(\bx,t) \cdot \bn^{out}_i 
&&= -\sigma^{e}\nabla M^{e}(\bx,t) \cdot \bn^{e},\quad &&\bx \in \Sigma_i,&\\
& \sigma^{out}\nabla M^{out}_i(\bx,t) \cdot \bn^{out}_i
 &&= \kappa^{out,e} \Bigl(M^{e}(\bx,t)- M^{out}_i(\bx,t)\Bigl), \quad&& \bx \in \Sigma_i.&
\end{alignat*}
Finally, on the outer boundary of the ECS we have 
\begin{align*}
0 = \sigma^{e}\nabla M^{e}(\bx,t) \cdot \bn^{e},\quad \bx \in \Psi.
\end{align*}
The BTPDE also needs initial conditions:
\begin{align*}
M^{in}_i(\bx,0) = \rho^{in},\quad 
M^{out}_i(\bx,0) = \rho^{out}, \quad 
M^{e}(\bx,0) = \rho^{e}. 
\label{eq:btpde_initc}
\end{align*}
where $\rho$ is the initial spin density.

The dMRI signal is measured at echo time $t=T_E > \Delta+\delta$ for
PGSE and $T_E > 2\sigma$ for OGSE\@.  This signal is the integral of
$M(\bx, T_E)$:
\be{eq:signal}
S := \int_{\bx\in \bigcup \{\Omega^{in}_i,\;\Omega^{out}_i,\;\Omega^e\}} M(\bx, T_E)\;d\bx.
\ee

In a dMRI experiment, the pulse sequence (time profile $f(t)$) is
usually fixed, while $\bg$ is varied in amplitude (and possibly also
in direction).  \soutnew{When $\bg$ varies only in amplitude (while staying in
the same direction), }{}  $S$ is \soutnew{}{usually} plotted against a quantity called
the $b$-value.  The $b$-value depends on $\bg$ and $f(t)$ and is
defined as
\ben
b(\bg) = \gamma^2 \|\bg\|^2 \int_0^{T_E} du\left(\int_0^u f(s) ds\right)^2.
\een
For PGSE, the b-value is \cite{Stejskal1965}:
\be{bvalue_pgse}
b(\bg,\delta,\Delta) = \gamma^2 \|\bg\|^2 \delta^2\left(\Delta-\delta/3\right).
\ee
For the cosine OGSE with {\it integer}\/ number of periods $n$ in each
of the two durations $\sigma$, the corresponding $b$-value is
\cite{Xu2007}:
\be{bvalue_ogse}
b(\bg,\sigma) =
\gamma^2 \|\bg\|^2 \frac{\sigma^3}{4 n^2\pi^2} = \gamma^2 \|\bg\|^2 \frac{\sigma}{\omega^2}.
\ee
The reason for these definitions is that in a homogeneous medium, the
signal attenuation is $e^{-\sigma b}$, where $\sigma$ is the intrinsic diffusion
coefficient.

\subsection{Fitting the ADC from the dMRI signal}

An important quantity that can be derived from the dMRI signal is 
the ``Apparent Diffusion Coefficient'' (ADC), which
gives an indication of the root mean squared distance travelled by water
molecules in the gradient direction $\bg/\|\bg\|$, averaged over all starting positions:

\be{eq:ADCdef}
  ADC  := \left. -\frac{\partial}{\partial b} \log{\frac{S(b)}{S(0)}}\right\vert_{b=0}. 
\ee
We numerically compute $ADC$ by a polynomial fit of 
$$
\log{S(b)} = c_0+c_1b+\cdots+c_n b^n,
$$
increasing $n$ from 1 onwards until we get the value of $c_1$ to be stable within a numerical
tolerance.

\subsection{
\soutnew{H-$ADC$}{HADC}
 model}
In a previous work \cite{schiavi2016}, a PDE model for the time-dependent 
ADC was obtained starting from the Bloch-Torrey equation, using homogenization techniques.
In the case of negligible water exchange between compartments (low permeability),
there is no coupling between the compartments, at least to the quadratic order in $\bg$, which is the ADC term.  
The ADC in compartment $\Omega$ is given by 
\begin{equation}
\label{eq:Deff}
HADC  =   \sigma -  \frac{1}{\int_0^{TE} F(t)^2 dt } \int_0^{TE} F(t)  \; h(t) \; dt, 
\end{equation}
where 
\soutnew{}{$F(t) = \int_0^{t} f(s)  \; ds,$ and}
\begin{equation}
\label{eq:ht_def}
h(t) = \frac{1}{|\Omega|}\int_{\partial \Omega}\omega(\bx,t)\left(\bug\cdot\bn\right)\,ds
\end{equation}
is a quantity related to the directional gradient of a function $\omega$ that is the solution of the homogeneous diffusion equation 
\soutnew{(DE)}{} with Neumann boundary condition and zero initial condition:
\begin{equation}
\label{eq:wProblem}
\begin{split}
\frac{\partial}{\partial t} \omega(\bx,t) 
-\nabla\left(\sigma \nabla \omega(\bx,t)\right)&=0,  \quad \quad \quad \quad \quad \quad \quad \quad \bx \in \Omega,\\
\sigma \nabla \omega(\bx,t) \cdot \bn &= \sigma F(t) \,\bug\cdot\bn, \quad \quad \quad \bx \in \partial \Omega, \\
\omega(\bx,0) &= 0,  \quad \quad \quad \quad \quad \quad \quad \quad \bx \in \Omega,
\end{split}
\end{equation}
$\bn$ being the outward normal and $t\in[0,TE]$, \soutnew{}{$\bug$ is the unit gradient direction}.  The above set of equations, (\ref{eq:Deff})-(\ref{eq:wProblem}), comprise the
homogenized model that we call the 
\soutnew{H-$ADC$}{HADC} model.

\subsection{Short diffusion time approximation of the ADC}\label{secShortADC}

A well-known formula for the ADC in the short diffusion time regime is the following
short time approximation (STA) \cite{Mitra1992,Mitra1993}:
\ben
STA = \sigma \left( 1 -  \frac{4\sqrt{\sigma}}{3 \; \sqrt{\pi}}\sqrt{\Delta} \frac{{A}}{dim \;V} \right),
\een
where $\dfrac{A}{V}$ is the surface to volume ratio
and $\sigma$ is the intrinsic diffusivity coefficient.  In the above formula the 
pulse duration $\delta$ is assumed to be very small compared to $\Delta$. 
A recent correction to the above formula \cite{schiavi2016}, taking into account the 
finite pulse duration $\delta$ and the gradient direction $\bug$, is the following:
\be{eq:STA_FP}
STA = \sigma  \left[ 1- \frac{4 \sqrt{\sigma}}{3 \; \sqrt{\pi} } 
					C_{\delta,\Delta}
					 \frac{A_\bug}{V } \right],
\ee
where 
\ben
A_\bug = \int_{\partial \Omega} \left(\bug \cdot \bn\right)^2\,ds,
\een
and
\begin{align*} 
C_{\delta,\Delta} &= \dfrac{4}{35} 
					 \dfrac{\left(\Delta+\delta \right)^{7/2}+ \left(\Delta - \delta \right)^{7/2} - 2 \left( \delta^{7/2}+\Delta^{7/2} \right)  }{\delta^2 \left(\Delta-\delta/3\right)}  \label{eqn_def_C_delDEL}
				 = \sqrt{\Delta} \left( 1 + \dfrac{1}{3} \dfrac{\delta}{\Delta} - \dfrac{8}{35} \left( \dfrac{\delta}{\Delta}\right)^{3/2} + \cdots  \right).
\end{align*}
When $\delta \ll \Delta$, the value $C_{\delta, \Delta}$ is approximately $\sqrt{\Delta}$.

\section{Method}\label{Method}
Below is a chart describing the work flow of SpinDoctor.
\begin{figure}[!htb]
\begin{tikzpicture}
\matrix (m)[matrix of nodes, column  sep=2cm,row  sep=8mm, align=center, nodes={rectangle,draw, anchor=center} ]{
   |[block]| {Read cells parameters}              &   \\
   |[block]| {Create cells \\(canonical configuration)}              &  |[decision]| {Plot cells } \\
   |[block]| {Read simulation domain parameters}              &  \\ 
    |[block]| {Create surface triangulation\\ (canonical configuration)}          &            |[decision]| {Plot surface triangulation}                                  \\
    |[block]| {Create FE mesh on canonical configuration; \\bend and twist the FE mesh nodes by analytical transformation.}         &           |[decision]| {Plot FE mesh}                                   \\
   |[block]| {Read experiment parameters}              &  |[decision]| {Compute STA}   \\ 
         |[block]| {Solve BTPDE}    &    |[block]| {Solve HADC}                                           \\
            |[decision]| {Plot magnetization and ADC }        &       |[decision]| {Plot HADC}                                      \\
};
\path [>=latex,->] (m-1-1) edge (m-2-1);
\path [>=latex,->] (m-2-1) edge (m-3-1);
\path [>=latex,->] (m-3-1) edge (m-4-1);
\path [>=latex,->] (m-4-1) edge (m-5-1);
\path [>=latex,->] (m-5-1) edge (m-6-1);
\path [>=latex,->] (m-6-1) edge (m-7-1);
\path [>=latex,->] (m-7-1) edge (m-8-1);
\path [>=latex,->] (m-2-1) edge (m-2-2);
\path [>=latex,->] (m-4-1) edge (m-4-2);
\path [>=latex,->] (m-5-1) edge (m-5-2);
\path [>=latex,->] (m-6-1) edge (m-7-2);
\path [>=latex,->] (m-6-1) edge (m-6-2);
\path [>=latex,->] (m-7-2) edge (m-8-2);
\end{tikzpicture}
\caption{Flow chart describing the work flow of SpinDoctor}
\end{figure}
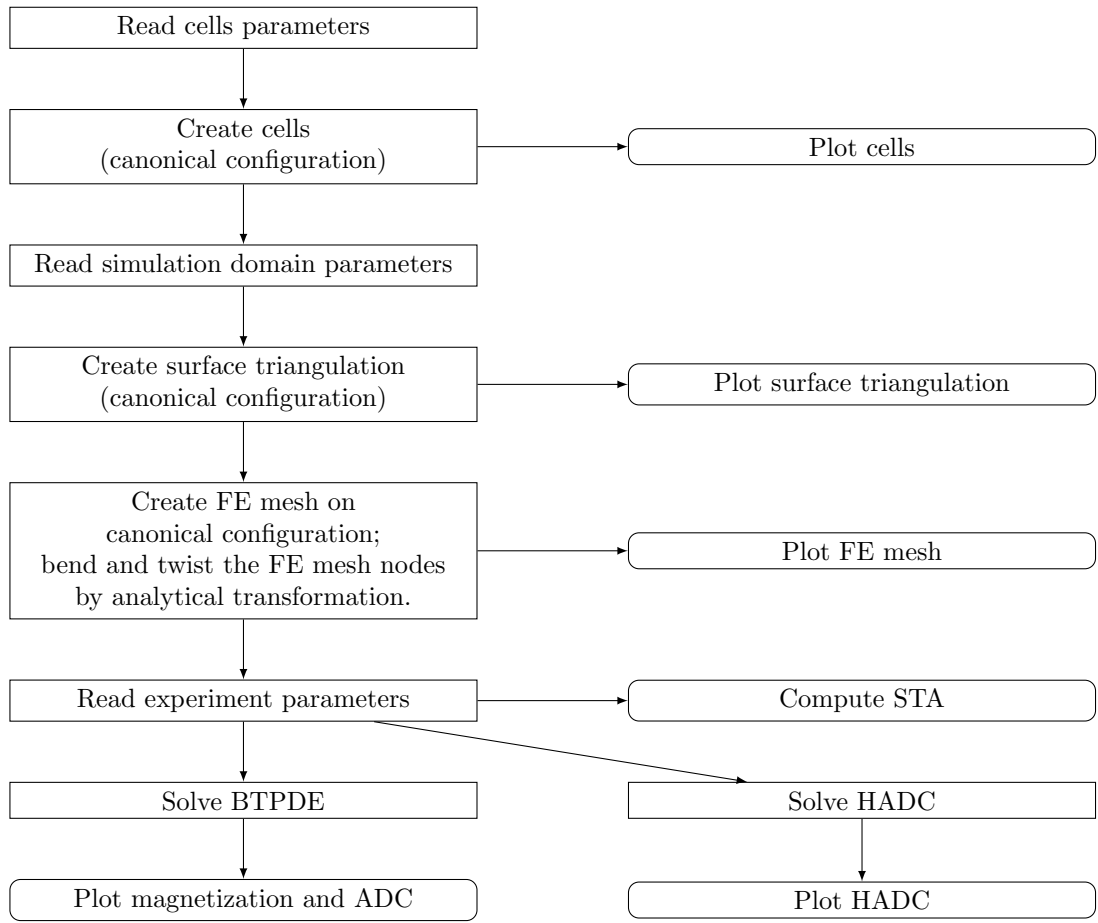

The physical units of the quantities in the input files for SpinDoctor are shown in Table \ref{table:units},
in particular, the length is in $\lunit$ and the time is in $\mu s$.
\begin{table}[!htb]
\begin{center}
\begin{tabular}{|p{4cm}|p{3cm}|}
\hline
Parameter &Unit    \\ \hline
length & \lunit \\\hline
time & $\mu s$ \\\hline
diffusion coefficient  &$\mu m^2/\mu s=\dunit$  \\\hline
permeability coefficient  &$\mu m/ \mu s = \kunit$  \\\hline
b-value & $\mu s/\mu m^2=\bunit$\\\hline
q-value & $(\mu s\mu m)^{-1}$\\\hline
\end{tabular}
\caption{Physical units of the quantities in the input files for SpinDoctor. \label{table:units}}
\end{center}
\end{table}
Below we discuss the various components of SpinDoctor in more detail.

\subsection{Read cells parameters}
The user provides an input file for the cell parameters, in the format described in Table \ref{table:params_cells}.
\begin{table}[!htb]
\begin{center}

\begin{tabular}{|p{0.5cm}|p{3cm}|p{2.5cm}|p{8cm}|}
\hline
Line&Variable name & Example & Explanation   \\
\hline
1 &cell\_shape& 1			& 1 = spheres; \newline 2 = cylinders;\\\hline
2&fname\_params\_cells& 'current\_cells'	& file name to store cells description\\\hline
3&ncell& 10										& number of cells\\\hline
4&Rmin& 1.5										& min Radius \\\hline
5&Rmax& 2.5										& max Radius \\\hline
6&dmin& 1.5										& min (\%) distance between cells: $dmin \times \frac{(Rmin+Rmax)}{2}$\\\hline
7&dmax& 2.5										& max (\%) distance between cells $dmax \times \frac{(Rmin+Rmax)}{2}$\\\hline
8&para\_deform& 0.05 0.05						& [$\alpha$\quad $\beta$]; \newline $\alpha$ defines the amount of bend; \newline $\beta$ defines the amount of twist \\\hline
9&Hcyl& 20										& height of cylinders \\\hline
\end{tabular}
\caption{Input file containing cells parameters.\label{table:params_cells}}
\end{center}
\end{table}

\subsection{Create cells (canonical configuration)}

SpinDoctor supports the placement of a group of non-overlapping cells in close vicinity to each other.  
There are two proposed configurations, one composed of spheres, the other composed of cylinders. 
The algorithm is described in Algorithm \ref{algo:create_cells}.

\begin{algorithm}[!htb]
Generate a large number of possible cell centers. \\
Compute the \soutnew{}{minimum} distance, $dist$, between the current center and previously accepted 
cells.  \\
Find the intersection of [$dist-dmax\times R_{mean}$, $dist-dmin\times R_{mean}$] and $[Rmin,Rmax]$, 
where $R_{mean} =\frac{Rmin+Rmax}{2} $.  If the intersection is not empty,
then take the middle of the intersection as the new radius and accept the new center.  Otherwise, reject the center. \\
Loop through the possible centers until get $n_{cell}$ accepted cells.
\caption{Placing $n_{cell}$ non-overlapping cells.\label{algo:create_cells}}
\end{algorithm}

\subsection{Plot cells}

\begin{figure}[!htb]
  \centering
\includegraphics[width=0.49\textwidth]{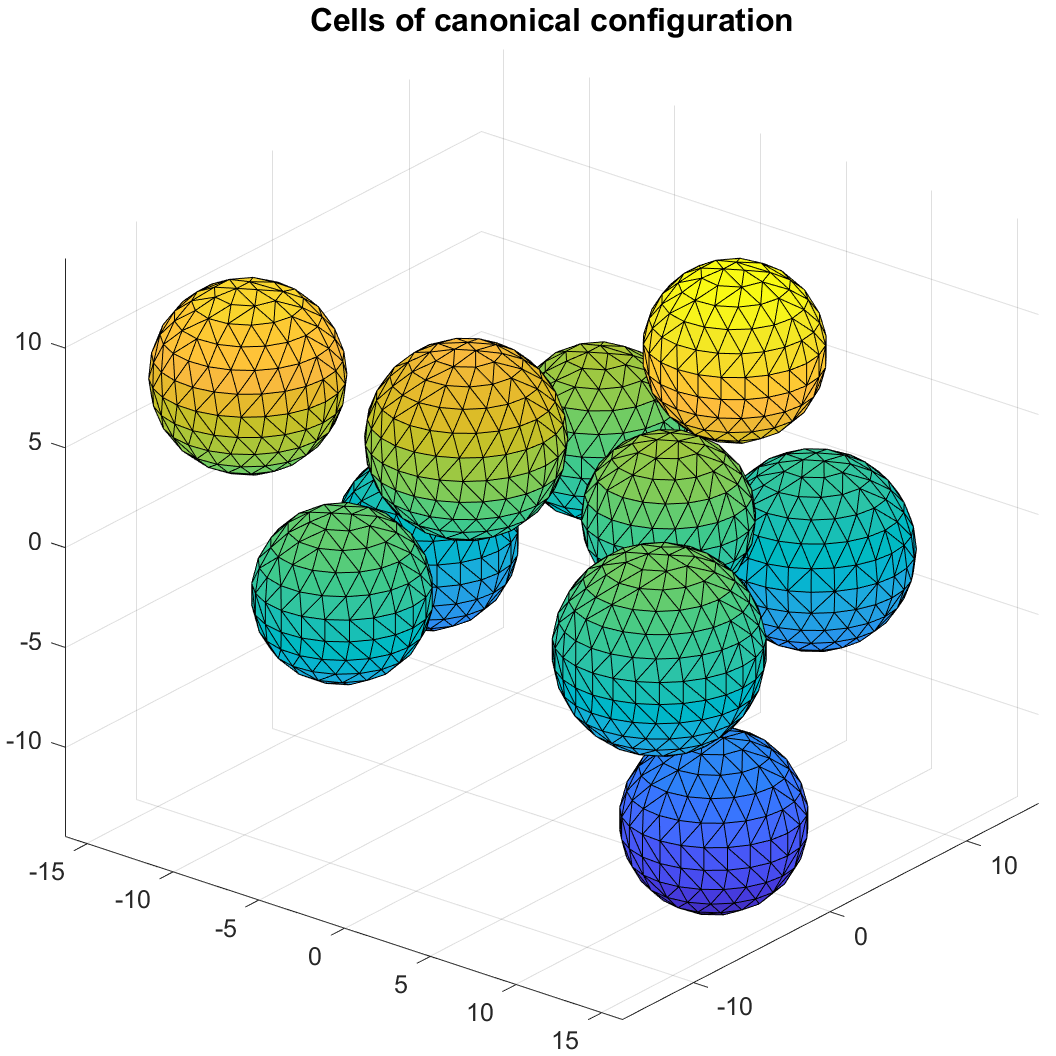} 
\includegraphics[width=0.49\textwidth]{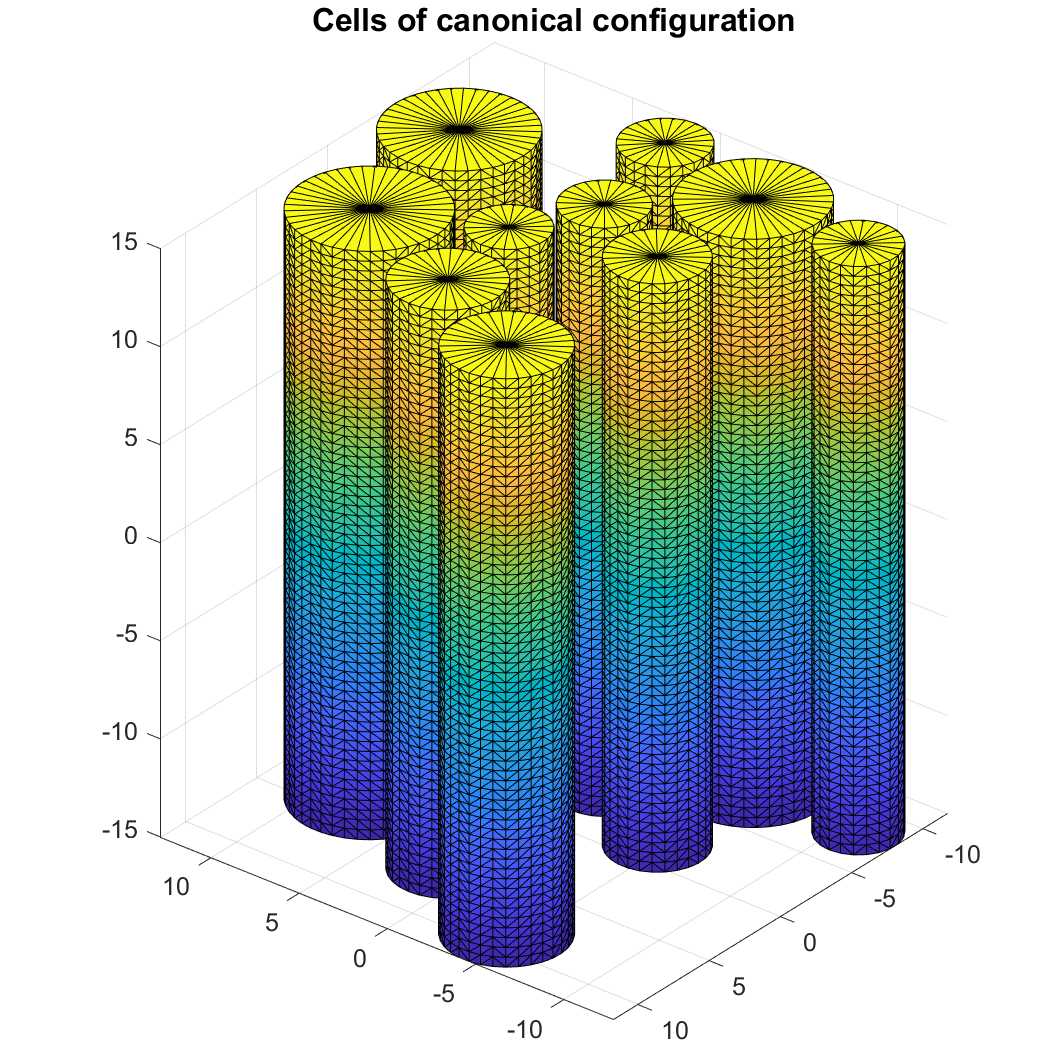}
  \caption{SpinDoctor plots cells in the canonical configuration.\label{fig:plot_cells}}
\end{figure}

SpinDoctor provides a routine to plot the cells to see if the configuration is acceptable
(see Fig. \ref{fig:plot_cells}).
\subsection{Read simulation domain parameters}
The user provides an input file for the simulation domain parameters, in the format described in Table \ref{table:params_simulation_domain}.

\begin{table}[!htb]
\begin{center}
\begin{tabular}{|p{0.5cm}|p{2.5cm}|p{3cm}|p{7cm}|}
\hline
Line&Variable name & Example & Explanation   \\
\hline
1 &Rratio& 0.0			& if Rratio is outside [0,1], it is set to 0;\newline else $Rratio = \frac{R_{in}}{R_{out}}$;\\\hline
2&include\_ECS& 2 										& 0 = no ECS; \newline 1 = box ECS; \newline 2 = tight wrap ECS; \\\hline
3&ECS\_gap& 0.3										& ECS thickness: \newline a. if box: as percentage of domain length;\newline b. if tight wrap: as percentage of mean radius \\\hline
4&dcoeff\_IN& 0.002									& diffusion coefficient in IN cmpt:\newline a. nucleus;\newline b. axon (if there is myelin); \\\hline
5&dcoeff\_OUT& 0.002									& diffusion coefficient in OUT cmpt:\newline a. cytoplasm; \newline b. axon (if there is no myelin); \\\hline
6&dcoeff\_ECS& 0.002									& diffusion coefficient in ECS cmpt;\\\hline
7&ic\_IN& 1			  							& initial spin density in In cmpt:\newline a. nucleus;\newline b. axon (if there is myelin) \\\hline
8&ic\_OUT& 1			  							& initial spin density in OUT cmpt:\newline a. cytoplasm; \newline b. axon (if there is no myelin); \\\hline
9&ic\_ECS& 1			  							& initial spin density in ECS cmpt:\\\hline
10&kappa\_IN\_OUT& 1e-5		   							& permeability between IN and OUT cmpts:\newline a. between nucleus and cytoplasm; \newline b. between axon and myelin;\\\hline
11&kappa\_OUT\_ECS& 1e-5									& permeability between OUT and ECS cmpts:\newline
a. if no nucleus: between cytoplasm and ECS;\newline b. if no myelin: between axon and ECS;\\\hline
12&Htetgen& -1										& Requested tetgen mesh size;\newline -1 = Use  tetgen default; \\\hline
13 & tetgen\_cmd&'SRC\slash TETGEN\slash tetGen\slash win64\slash tetgen'		& path to  tetgen\_cmd\\\hline
\end{tabular}
\caption{Input file of simulation domain parameters.}
\label{table:params_simulation_domain}
\end{center}
\end{table}

\subsection{Create surface triangulation}

Finite element mesh generation software requires a good surface triangulation.  This means 
the surface triangulation needs to be water-tight and does not self-intersect.  How closely these requirements
are met in floating point arithmetic has a direct impact on the quality of the finite element mesh generated.

It is often difficult to produce a good surface triangulation for arbitrary geometries.  
Thus, we restrict the allowed shapes 
to cylinders and spheres.    Below in Algorithms \ref{algo:surface_triangulation_spheres} 
and \ref{algo:surface_triangulation_cylinders} we describe how to obtain a surface triangulation 
for spherical cells with nucleus, 
cylindrical cells with myelin layer, and the ECS (box or tightly wrapped).
We describe a canonical configuration where the cylinders are placed parallel to the 
$z$-axis.  More general shapes are obtained from the canonical configuration 
by coordinate transformation 
in a later step.

\begin{algorithm}[!htb]
Suppose we have $n_{cell}$ spherical cells with nucleus.    
Denote a sphere with center $c$ and radius $R$ by $S(c,R)$, 
we use the built-in functions (convex hull, delaunnay triangulation) 
in MATLAB to get its surface triangulation, $T(c,R)$.  
Call the radii of the nucleus $r_1,\cdots,r_{ncell}$ and the radii of the cells
$R_1,\cdots,R_{ncell}$.   Then the boundaries between the cytoplasm and the nucleus are 
$$\{\Gamma_i= T(c_i,r_i)\} , i = 1,\cdots,ncell;$$ 
and between the cytoplasm and the ECS
$$\{\Sigma_i = T(c_i,R_i)\}, i = 1,\cdots,ncell;$$\\
For the box ECS, we find the coordinate limits of the set
$$
\bigcup_i S(c_i,R_i) \in [x_0,x_f]\times [y_0,y_f] \times [z_0,z_f]
$$ 
and add a gap $k = \text{ECS\_gap}\times \max\{x_f-x_0,y_f-y_0,z_f-z_0\}$ to make a box
$$
B =  [x_0-k,x_f+k]\times [y_0-k,y_f+k] \times [z_0-k,z_f+k].
$$ 
We put 2 triangles on each face of $B$ to make a surface triangulation $\Psi$ with 12 triangles.\\
For the tight-wrap ECS, we increase the cell radius by a gap size and take the union
$$
W = \bigcup_i S(c_i,R_i+ \text{ECS\_gap}\times R_{mean}),
$$ 
where $R_{mean} =\frac{Rmin+Rmax}{2} $.  We use the alphaShape function in MATLAB to find a surface triangulation $\Psi$
that contains $W$.  
\caption{Surface triangulation of spherical cells and ECS. \label{algo:surface_triangulation_spheres} }
\end{algorithm}
\newpage
\begin{algorithm}[!htb]
Suppose we have $n_{cell}$ cylindrical cells with a myelin layer, all with height $H$.  
Denote a disk with center $c$ and radius $R$ by $D(c,R)$, and the circle with 
the same center and radius by $C(c,R)$.
Let the radii of the axons be $r_1,\cdots,r_{ncell}$ and the radii of the cells be 
$R_1,\cdots,R_{ncell}$, meaning the thickness of the myelin layer is $R_i - r_i$.   \\
The boundary between the axon and the myelin layer is:
$$ C(c_i,r_i) \times [-H/2,H/2]$$
We discretize $C(c_i,r_i)$ as a polygon $P(c_i,r_i)$ and place one at $z = -H/2$ and one at $z = H/2$. 
Then we connect the corresponding vertices of $P(c_i,r_i)\times\{-H/2\}$ and $P(c_i,r_i)\times\{H/2\}$
and add a diagonal on each panel to get a surface triangulation $\Gamma_i$.\\
Between the myelin layer and the ECS we discretize $C(c_i,R_i)$ as a polygon and place one at $z=-H/2$ 
and one at $z = H/2$ to get a surface triangulation $\Sigma_i$.\\
For the box ECS, we find the coordinate limits of the union of $D(c_i,R_i)$  
and add a gap to make a rectangle in two dimensions.
Then we place the rectangle at $z = -H/2$ and at $z = H/2$ to get a box.  Finally, the box is given a surface triangulation
with 12 triangles.
\\
For tight-wrap ECS, we increase the cell radius by a gap size and take the union
$$
W = \bigcup_i D(c_i,R_i+kR_{mean}). 
$$ 
We use the alphaShape function in MATLAB to find a two dimensional polygon $Q$ that contains $W$.  
We place $Q$ at $z = -H/2$ and at $z = H/2$ and connect correponding vertices, adding a diagonal on each panel.
Suppose $Q$ is a polygon with $n$ vertices, then the surface triangulation of the side of the ECS will have $2n$ triangles.
\\The above procedure produces a surface triangulation for the boundaries that are parallel to $z$-axis.  We now 
must close the top and bottom.  The top and bottom boundaries is just the interior of $Q$.  However, the surface 
triangulation cannot be done on $Q$ directly.  We must cut out \soutnew{$D(c_i,R_i)$}{$D(c_i,r_i)$}, the disk which touches the 
axon, and $A_i = D(c_i,R_i)-D(c_i,r_i)$, the annulus which touches the myelin.  Then we triangulate
$Q - \bigcup_i D_i - \bigcup_i A_i$ using the MATLAB built-in function that triangulates a polygon
with holes to get the boundary that touches the ECS.
The surface triangulation for $A_i$ and \soutnew{$D(c_i,R_i)$ }{$D(c_i,r_i)$} are straightforward.
\caption{Surface triangulation of cylindrical cells and ECS. \label{algo:surface_triangulation_cylinders} }
\end{algorithm}

\subsection{Plot surface triangulation}
\begin{figure}[!htb]
  \centering
\includegraphics[width=0.49\textwidth]{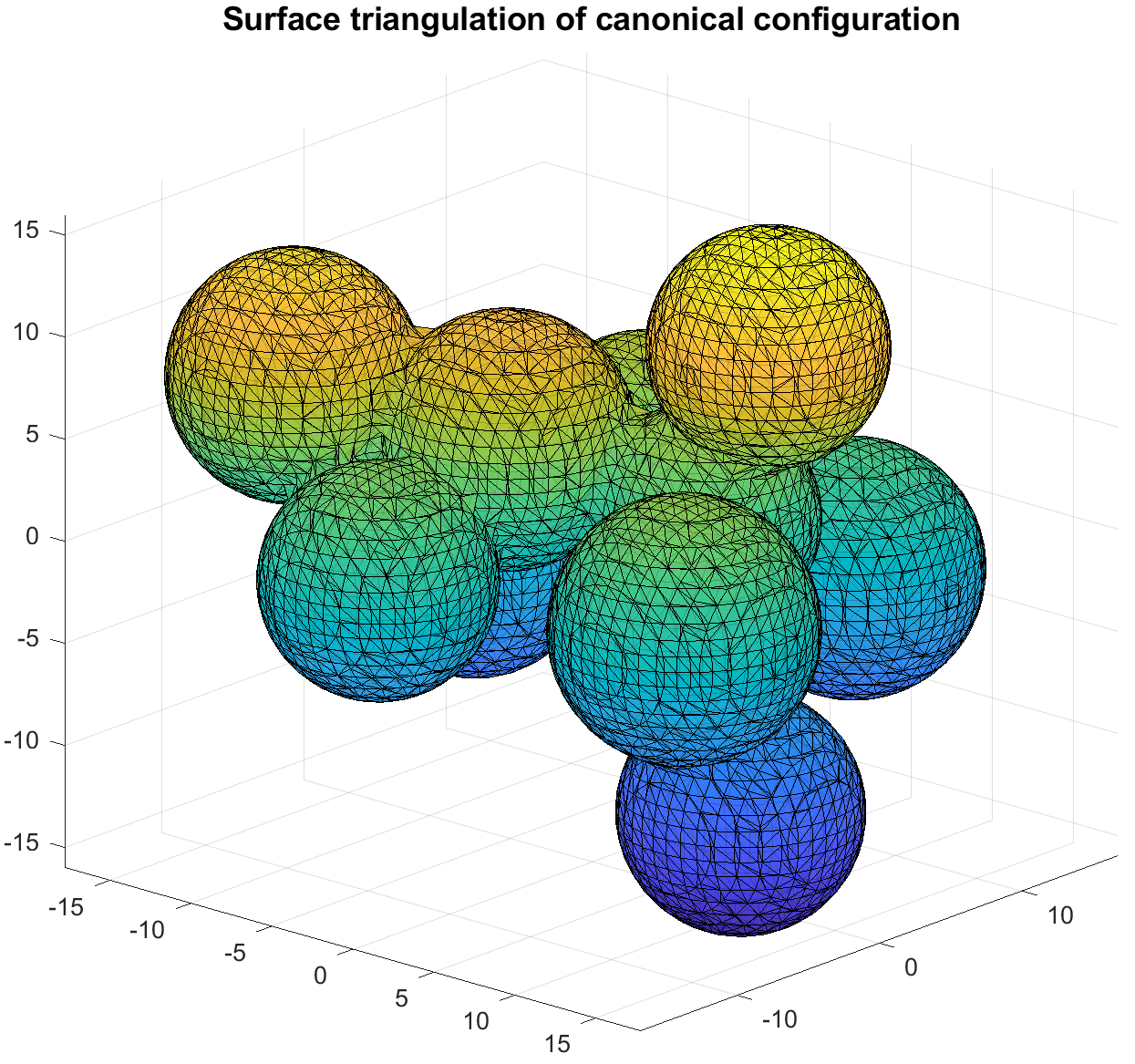} 
\includegraphics[width=0.49\textwidth]{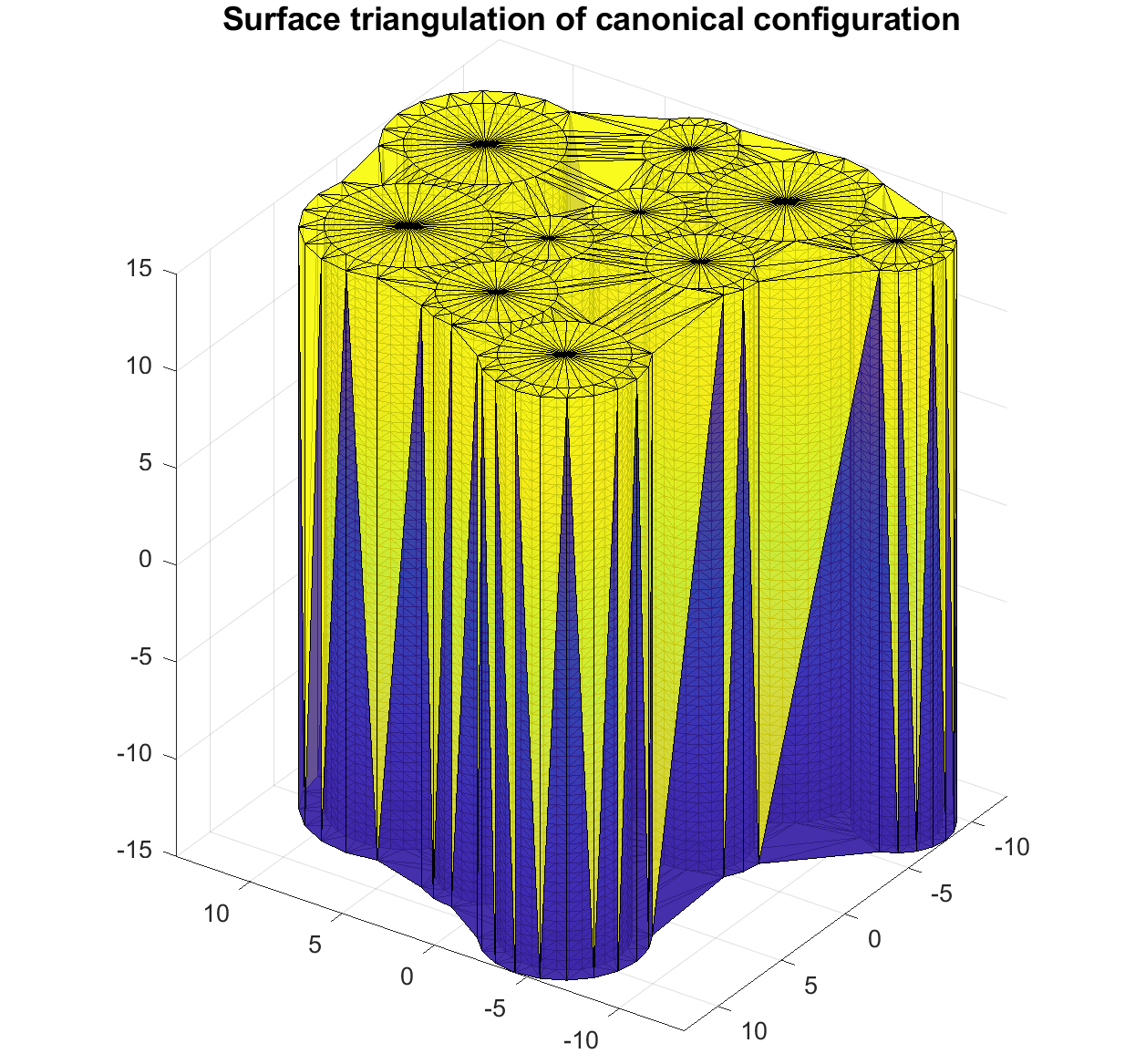}
  \caption{SpinDoctor plots the surface triangulation of the canonical configuration. 
Left: spherical cells with ECS; Right: cylindrical cells with ECS.
\label{fig:plot_surface_triangulation}}
\end{figure}

SpinDoctor provides a routine to plot the surface triangulation (see Fig. \ref{fig:plot_surface_triangulation}).

\subsection{Finite element mesh generation}

SpinDoctor calls Tetgen \cite{Si2015}, an external package (executable files are included in the toolbox package), to create 
a tetrehedra finite elements mesh from the surface triangulation generated by Algorithms
\ref{algo:surface_triangulation_spheres} 
and \ref{algo:surface_triangulation_cylinders}.  
The FE mesh is generated on the canonical configuration.  
The numbering of the compartments and boundaries used by SpinDoctor are given in
Tables \ref{table:cmpts_labels} and \ref{table:bdys_labels}. 
The labels are related to the values of the intrinsic diffusion coefficient, the initial spin
density, and the permeability requested by the user.  
Then the FE mesh nodes are deformed analytically by a coordinate transformation, described 
in Algorithm \ref{algo:bend_twist}.
\begin{table}[!htb]
\begin{center}
\begin{tabular}{|p{2cm}|p{3cm}|p{3cm}|p{3cm}|}
\multicolumn{4}{c}{Spherical cells without nucleus }\\\hline
 Cmpt & Cytoplasm & Nucleus  & ECS \\\hline
Label & OUT &  & ECS \\\hline
Number&  $[1:n_{cell}]$ &  &$n_{cell}+1$\\\hline
\multicolumn{4}{c}{Spherical cells with nucleus}\\\hline
 Cmpt & Cytoplasm & Nucleus  & ECS \\\hline
Label & OUT & IN & ECS \\\hline
Number &  $[1:n_{cell}]$& $[n_{cell}+1:2n_{cell}]$ & $2n_{cell}+1$\\\hline
\multicolumn{4}{c}{Cylindrical cells without myelin  }\\\hline
Cmpt & Axon & Myelin   & ECS \\\hline
Label & OUT &  & ECS \\\hline
Number &  $[1:n_{cell}]$ &  &$n_{cell}+1$\\\hline
\multicolumn{4}{c}{Cylindrical cells with myelin }\\\hline
Cmpt & Axon & Myelin   & ECS \\\hline
Label & IN & OUT & ECS \\\hline
Number &  $[1:n_{cell}]$& $[n_{cell}+1:2n_{cell}]$ & $2n_{cell}+1$\\\hline

\end{tabular}
\end{center}
\caption{The labels and numbers of compartments.\label{table:cmpts_labels}}
\end{table}

\begin{table}[!htb]
\begin{center}
\begin{tabular}{|p{2cm}|p{3cm}|p{3.5cm}|p{4.5cm}|}
\multicolumn{4}{c}{Spherical cells without nucleus }\\\hline
Boundary & Sphere &   & Outer ECS boundary \\\hline
Label & OUT\_ECS &  & $\kappa = 0$\\\hline 
Number &  $1:n_{cell}$ &  &$n_{cell}+1$\\\hline
\multicolumn{4}{c}{Spherical cells with nucleus}\\\hline
 Boundary & Outer sphere &   Inner sphere & Outer ECS boundary \\\hline
Label & OUT\_ECS & IN\_OUT & $\kappa = 0$\\\hline 
Number &  $1:n_{cell}$& $n_{cell}+1:2n_{cell}$ & $2n_{cell}+1$\\\hline
\multicolumn{4}{c}{Cylindrical cells without myelin}\\\hline
Boundary & Cylinder\newline side wall & Cylinder\newline top and bottom   & Outer ECS boundary\newline minus cylinder top/bottom \\\hline
Label & OUT\_ECS & $\kappa = 0$ & $\kappa = 0$\\\hline 
 Number &  $2[1:n_{cell}]-1$ &  $2[1:n_{cell}]$ &$2n_{cell}+1$\\\hline
\multicolumn{4}{c}{Cylindrical cells with myelin}\\\hline
Boundary & Inner cylinder \newline side wall &   Inner cylinder \newline top and bottom & \\\hline
Label & IN\_OUT & $\kappa = 0$ & \\\hline 
 Number &  $4[1:n_{cell}]-3$& $4[1:n_{cell}]-2$ & \\\hline
 & Outer cylinder\newline side wall &   Outer cylinder\newline top and bottom & Outer ECS boundary \newline minus cylinder top/bottom \\\hline
Label & OUT\_ECS & $\kappa = 0$ & $\kappa = 0$\\\hline 
Number &  $4[1:n_{cell}]-1$& $4[1:n_{cell}]$ & $4n_{cell}+1$\\\hline
\end{tabular}
\end{center}
\caption{The labels and numbers of boundaries.  \label{table:bdys_labels}}
\end{table}

\begin{algorithm}[!htb]
The external package Tetgen \cite{Si2015} generates the finite element mesh that keeps track of the different compartments and the interfaces between them.  The mesh is saved in several text files.  \\
The connectivity matrices of the finite elements and facets are not modified by the coordinates 
transformation described below.  The nodes are transformed by 
bending and twisting as described next.\\
The set of FE mesh nodes $\{x_i,y_i,z_i\}$ are transformed in the following ways: 

Twisting around the $z$-axis with a user-chosen twisting parameter $\alpha_{twist}$ is defined by
\ben
\begin{split}
\vthree{x}{y}{z} & \rightarrow  \begin{bmatrix} 
\cos(\alpha_{twist} z) & -\sin(\alpha_{twist} z) & 0 \\
\sin(\alpha_{twist} z) & \cos(\alpha_{twist} z) & 0 \\
0 & 0 & 1
    \end{bmatrix} \vthree{x}{y}{z}.
\end{split}
\een

Bending on the $x-z$ plane with a user-chosen bending parameter $\alpha_{bend}$ is defined by
\ben
\vthree{x}{y}{z} \rightarrow 
\vthree{x+
\alpha_{bend}z^2}{y}{z}.
\een

Given $[\alpha_{bend},\alpha_{twist}]$, bending is performed after twisting.
\caption{Bending and twisting of the FE mesh of the canonical configuration.\label{algo:bend_twist}}
\end{algorithm}

\subsection{Plot FE mesh}
\begin{figure}[!htb]
  \centering
\includegraphics[width=0.49\textwidth]{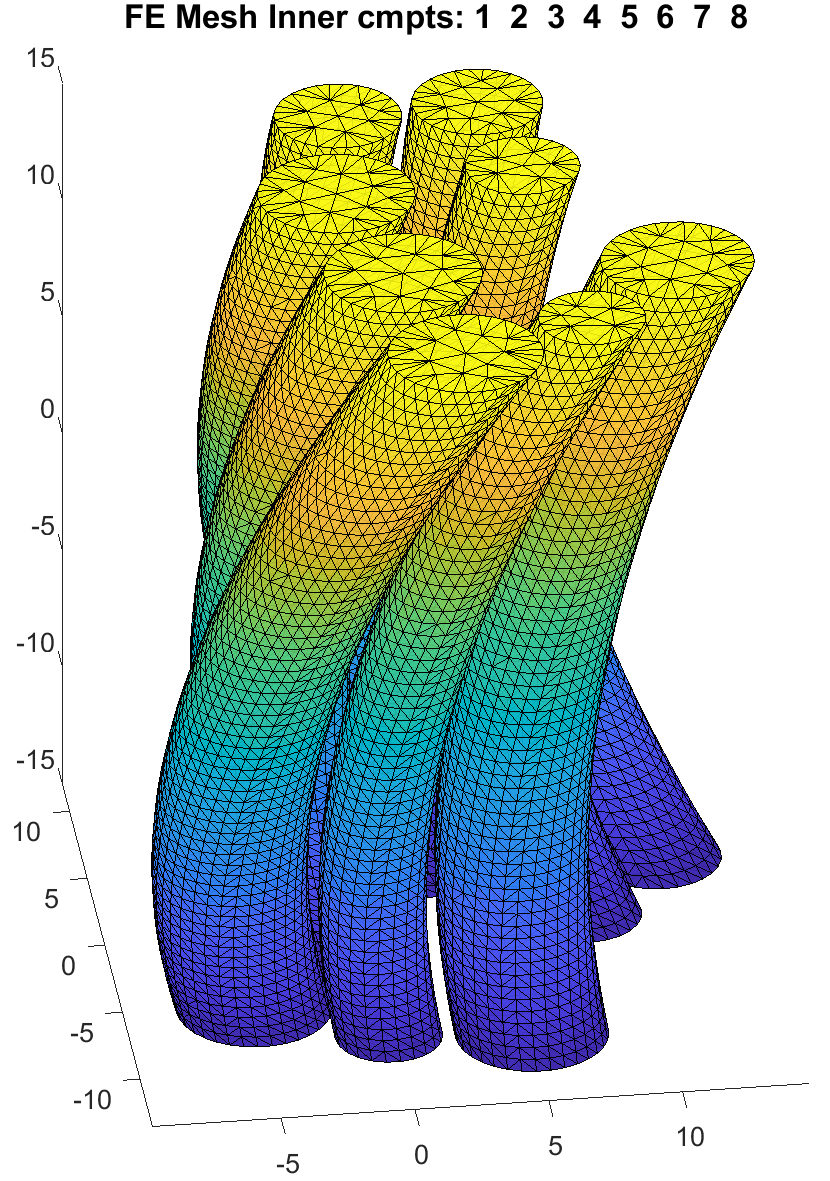}
\includegraphics[width=0.49\textwidth]{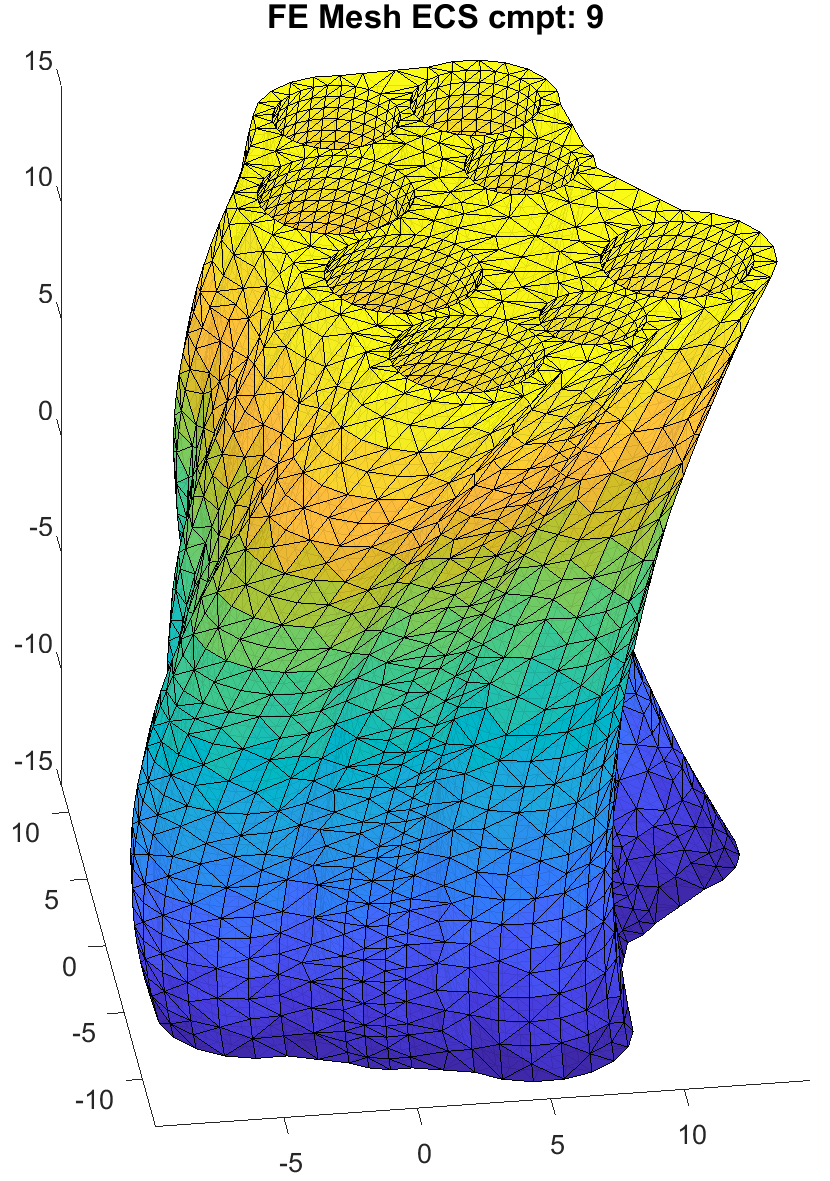}
  \label{two_profiles}
  \caption{FE mesh of cylinders and ECS after bending and twisting.  
Compartment number is 1 to 8 for the cylinders and 9 for the ECS.  \label{fig:plot_fe_mesh}}
\end{figure}

SpinDoctor provides a routine to plot the FE mesh (see Fig. \ref{fig:plot_fe_mesh} for cylinders and ECS that have been bent and twisted).

\subsection{Read experimental parameters}
The user provides an input file for the simulation experimental parameters, in the format described in Table \ref{table:params_simulation_experiment}.
\begin{table}[!htb]
\begin{center}
\begin{tabular}{|p{0.5cm}|p{2.2cm}|p{3cm}|p{7.5cm}|}
\hline
Line&Variable name & Example & Explanation   \\
\hline
1 & ngdir&20		     						& number of gradient direction;\newline if $ngdir > 1$, the gradient directions are distributed uniformly on a sphere; \newline if $ngdir = 1$, take the gradient direction from the line below;\\\hline
2 & gdir&1.0 0.0 0.0 		     				& gradient direction;  No need to normalize;\\\hline
3& nexperi &3 												& number of experiments; \\\hline
4& sdeltavec&2500 10000 10000 				& small delta;	\\\hline
5& bdeltavec&2500 10000 10000 					& big delta;	\\\hline
6& seqvec&1 2 3       								& diffusion sequence of experiment; \newline 
1 = PGSE; 2 = OGSEsin; 3 = OGSEcos; \\\hline
7& npervec&0 10 10     									& number of period of OGSE;\\\hline
8 & solve\_hadc & 1 & 0 = do not solve HADC;\newline Otherwise solve HADC; \\\hline
9& rtol\_deff, atol\_deff&1e-4 1e-4 									& $[r_{tol}\quad a_{tol}]$;
relative and absolute tolerance for HADC ODE solver;   \\\hline
10 & solve\_btpde & 1 & 0 = do not solve BTPDE;\newline Otherwise solve BTPDE; \\\hline
11& rtol\_bt, atol\_bt&1e-5 1e-5 									& $[r_{tol}\quad a_{tol}]$;
relative and absolute tolerance for BTPDE ODE solver;  \\\hline
12& nb &2 												& number of b-values;\\\hline
13& blimit &0												& 0 = specify bvec;\newline 1 = specify [bmin,bmax]; \newline 2  = specify [gmin,gmax];\\\hline
14& const\_q & 0												& 0: use input bvalues for all experiments; \newline 
1: take input bvalues for the first experiment and use the same q for the remaining experiments  \\\hline
15&bvalues&0 50 100 200  								& bvalues or [bmin, bmax] or [gmin, gmax];\newline depending on line 13;\\\hline
\end{tabular}
\caption{Input file for simulation experiment parameters.\label{table:params_simulation_experiment}}
\label{tab:experiment_parameters}
\end{center}
\end{table}

\subsection{BTPDE}
The spatial discretization of the BTPDE is based on a finite element method where
{\it  interface (ghost) elements} \cite{Nguyen2014} are used  to impose the permeable interface conditions. 
The time stepping is done using the MATLAB built-in ODE routine ode23t.
See Algorithm \ref{algo:btpde}.

\begin{algorithm}[!htb]
FE matrices are generated for each compartment by the finite element method with continuous piecewise linear basis functions (known as $P_1$).
The basis functions are denoted as $\varphi_k$ for $k=1,\dots, N_v$, where $N_v$ denotes the number of mesh nodes (vertices).  All matrices are sparse matrices. $\bm M$ and $\bm S$ are known in the FEM literature as mass and stiffness matrices  which are defined as follows:
\begin{equation*}
    \bm M_{ij} =\int_\Omega \varphi_i \varphi_j \, d\bx, \qquad \qquad  \bm S_{ij} =\int_\Omega \sigma_{i}\,\nabla \varphi_i \cdot \nabla \varphi_j \, d\bx.
\end{equation*}
$\bm J$ has a similar form as the mass matrix but it is scaled with the coefficient $\bg\cdot \bx$, we therefore call it the scaled-mass matrix
\begin{equation*}
    \bm J_{ij} =\int_\Omega  \bg\cdot\bx\, \varphi_i \varphi_j \, d\bx.
\end{equation*}
We construct the matrix based on the flux matrix $\bm Q$
\begin{equation*}
    \bm Q_{ij} =\int_{\partial \Omega}  w \,\varphi_i \varphi_j \, ds
\end{equation*}
where a scalar function $w$ is used as an interface marker.
The matrices are assembled from local element matrices and the assembly process is based on vectorized routines of \cite{RahmanValdman2013}, which replace expensive loops over elements by operations with 3-dimensional arrays. All local elements matrices in the assembly of $\bm S, \bm M, \bm J$ are evaluated at once and stored in a full matrix of size $4 \times  4 \times N_e$, where $N_e$ denotes the number of tetrahedral elements. The assembly of $\bm Q$ is even simpler; all local matrices are stored in a full matrix of size $3 \times  3 \times n_{be}$, where $n_{be}$ denotes the number of boundary triangles.
\\
Double nodes are placed at the interfaces between compartments connected by permeable membrane.  $\overline{\bm Q}$ is used to impose the interface conditions and it is associated with the {\it  interface (ghost) elements}.   Specifically, assume that the double nodes are defined in a pair of indices $\{i,\bar{i}\}$, $\overline{\bm Q}$ is defined as the following
\begin{equation*}
    \overline{\bm Q}_{ij} = 
    \begin{cases}
        \bm Q_{ij}, & \mbox{if vertex $i$ and $j$ belong to one interface}\\
        -\bm Q_{\bar{i}\bar{j}} & \mbox{if vertex $i$ and $j$ belong to two different interfaces}
    \end{cases}
\end{equation*}
\\
The fully coupled linear system has the following form
\begin{equation}\label{eq:BTPDE_Matrixform}
\bm M \frac{\partial \xi}{\partial t}=-\Bigl(I\gamma f(t)\,{\bm J}+\bm S+ \overline{\bm Q}\Bigl)\, \xi
\end{equation}
where $\xi$ is the approximation of the magnetization $M$. 
SpinDoctor calls MATLAB built-in ODE routine ode23t to solve the 
\soutnew{semi-discretized equation}{semi-discretized system of equations}.
\caption{BTPDE.\label{algo:btpde}}
\end{algorithm}

\subsection{
\soutnew{H-$ADC$}{HADC} 
model}
Similarly, the DE of the 
\soutnew{H-$ADC$}{HADC} model is discretized by finite elements.
See Algorithm \ref{algo:hadc}.

\begin{algorithm}[!htb]
Eq. (\ref{eq:wProblem}) can be discretized similarly as described for the BTPDE 
and has the matrix form
\begin{equation}\label{eq:H_ADC_model}
\bm M \frac{\partial \zeta}{\partial t}=-\bm S\, \zeta+{\bm Q} \, \bar{\zeta}
\end{equation}
where $\zeta$ is the approximation of $w$ and $\bar{\zeta}_i=\sigma_i\, F(t)\,{\bm u}_g\cdot\bn(\bx_i)$. We note that the matrices here are assembled and solved separately for each compartment.
SpinDoctor calls MATLAB built-in ODE routine ode23t to solve the semi-discretized equation.
\caption{HADC model.\label{algo:hadc}}
\end{algorithm}

\subsection{Some important output quantities}
In Table \ref{table:outputs} we list some useful quantities that are the outputs of SpinDoctor.
The braces in the "Size" column denote MATLAB cell data structure and the brackets denote MATLAB matrix data structure.
\begin{table}[!htb]
\begin{center}
\begin{tabular}{|p{3cm}|p{5.5cm}|p{5cm}|}
\hline
Variable name & Size & Explanation   \\
\hline
TOUT &  \{{nexperi}$\times${nb}$\times${Ncmpt}\}[1 $\times$nt]& ODE time discretization \\\hline
YOUT &  \{{nexperi}$\times${nb}$\times${Ncmpt}\}[Nnodes $\times$nt] & Magnetization \\\hline
MF\_cmpts&[Ncmpt $\times$ nexperi $\times$ nb]&integral of magnetization at $TE$ in each compartment. \\\hline
MF\_allcmpts&[nexperi $\times$ nb]&integral of magnetization at $TE$ summed over all compartments. \\\hline
ADC\_cmpts&[Ncmpt $\times$ nexperi ]& ADC in each compartment. \\\hline
ADC\_allcmpts&[nexperi $\times$ 1]&ADC accounting for all compartments. \\\hline
ADC\_cmpts\_dir&[ngdir $\times$ Ncmpt $\times$ nexperi ]& ADC in each compartment in each direction. \\\hline
ADC\_allcmpts\_dir&[ngdir $\times$ nexperi $\times$ 1]&ADC accounting for all compartments in each direction. \\\hline
\end{tabular}
\end{center}
\caption{Some important SpinDoctor output quantities.\label{table:outputs}}
\end{table}

\section{\soutnew{Numerical results}{SpinDoctor examples}}

\label{Numericalresults}

\soutnew{}{In this section we show some prototypical examples using the available functionalities of SpinDoctor.}
\subsection*{\soutnew{Comparison with reference solution}{}}
\soutnew{
First we compare the results from SpinDoctor with the reference solution using the Matrix Formalism method \cite{Grebenkov2010a}
for a variety of $b$-values and permeability coefficients.  Plotted in Fig. \ref{fig:reference_solution}
are the dMRI signals using the BTPDE solve in SpinDoctor
and the Matrix Formalism method.  The maximum relative errors between SpinDoctor and the 
reference signal are 2.26\%, 1.06\%, and 0.68\% for the three permeability
coefficients simulated.}{}

\subsection{Comparison of BTPDE and HADC with Short Time Approximation}

\begin{figure}[!htb]
  \centering
\includegraphics[width=0.475\textwidth]{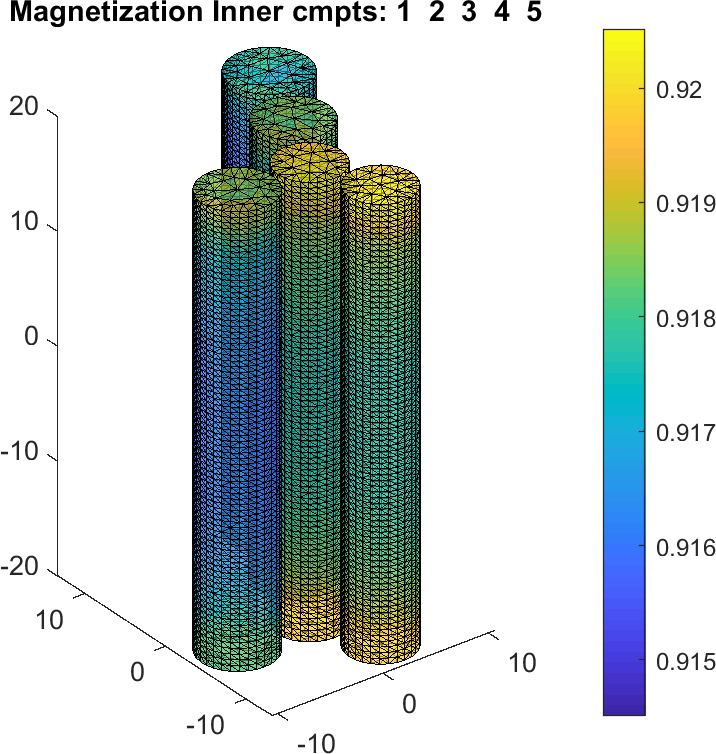} \quad \
\includegraphics[width=0.475\textwidth]{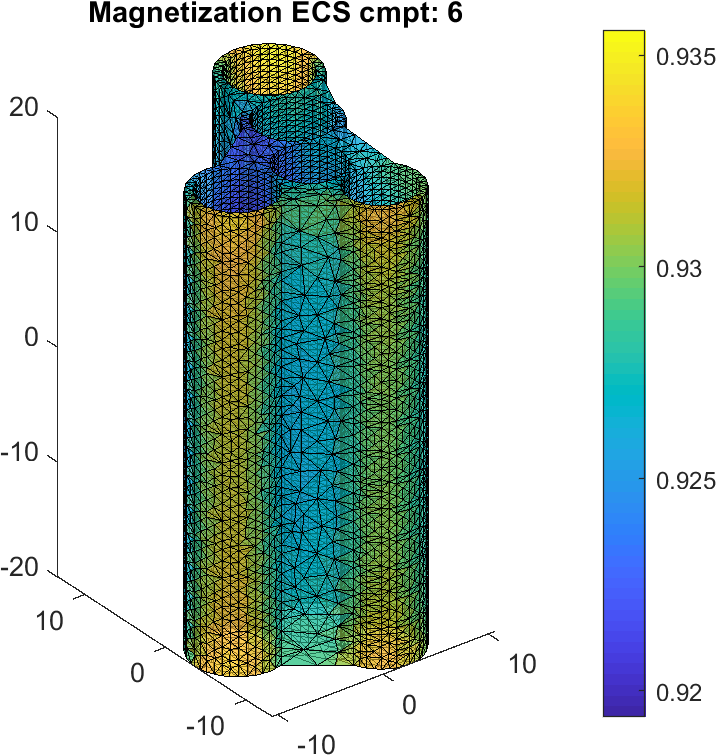}
\\
\vspace{0.8cm}
\includegraphics[width=0.32\textwidth]{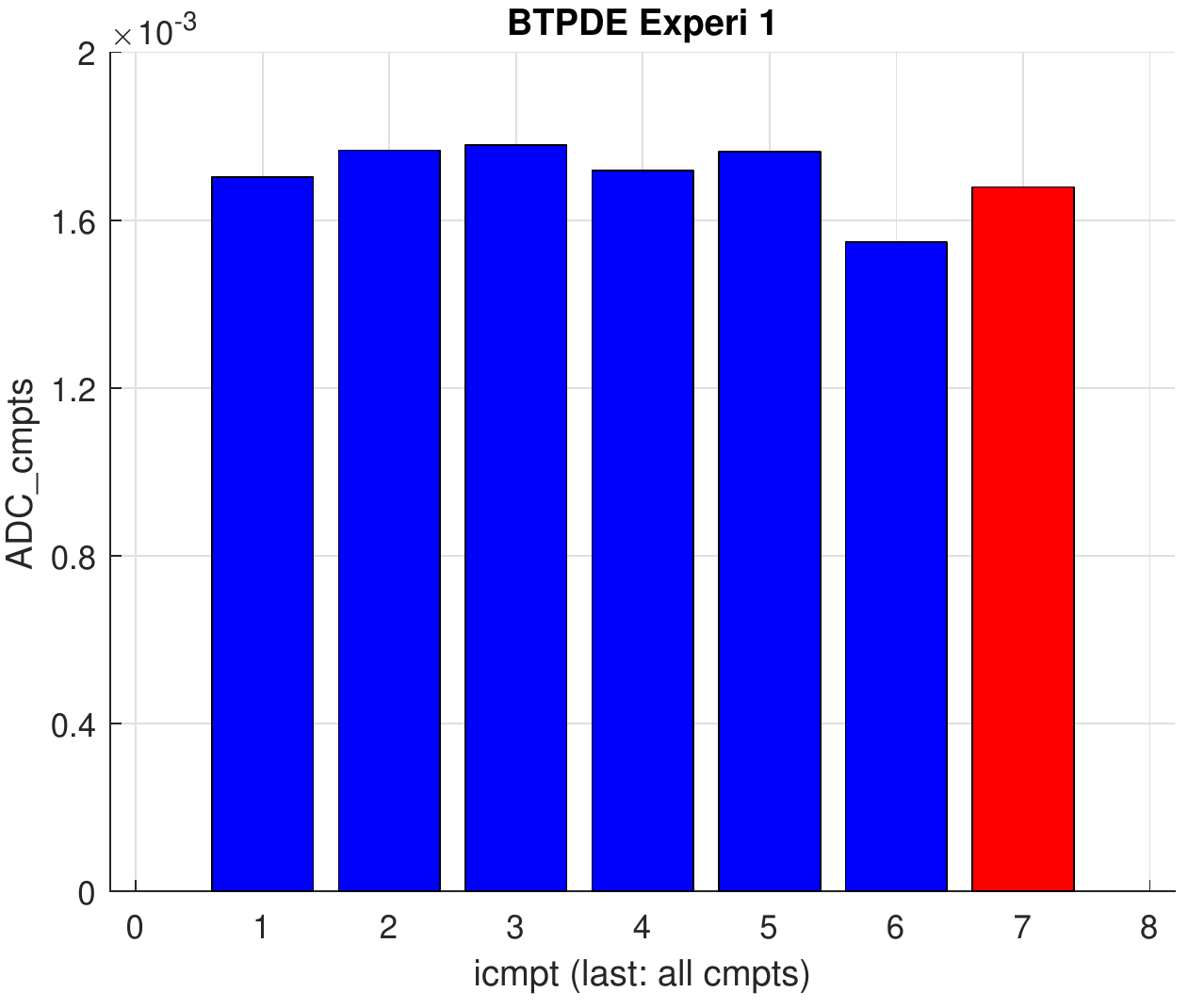} \
\includegraphics[width=0.32\textwidth]{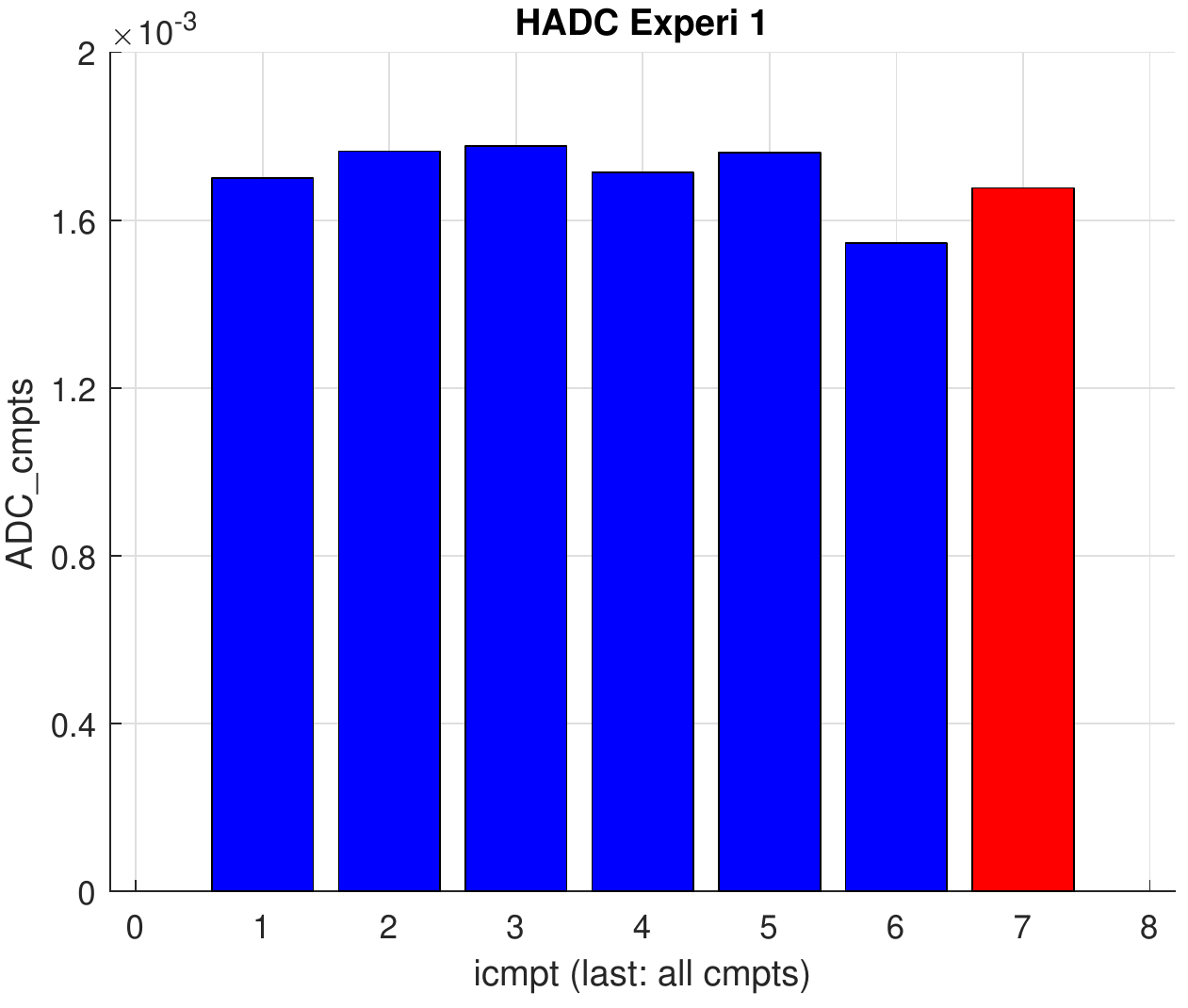} \
\includegraphics[width=0.32\textwidth]{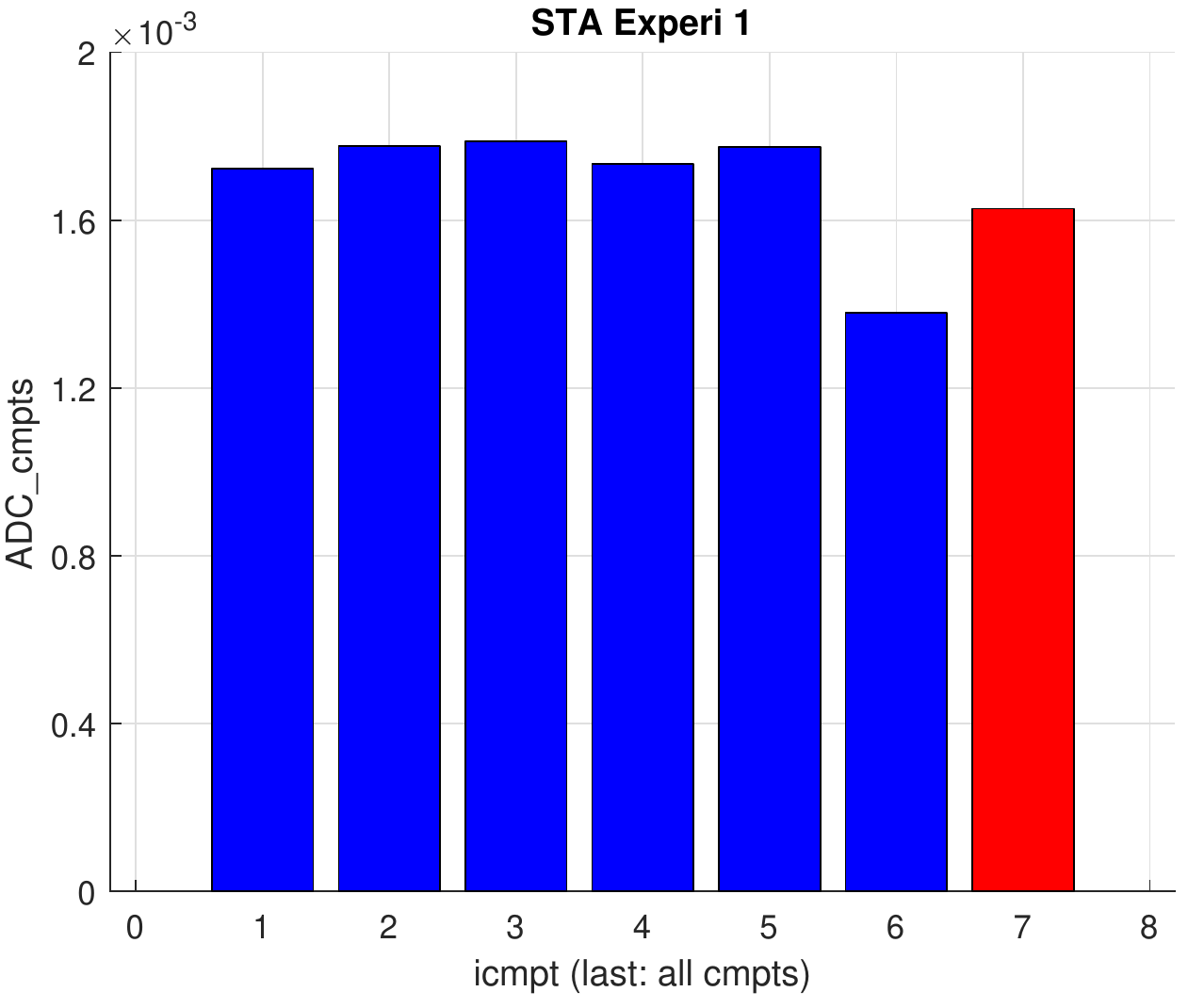}
  \label{2layer_sphere}
  \caption{Geometry: 5 cylinders, tight wrap ECS, ECS gap = 0.2, $\bug = [1,1,1]$,  $\sigma^{out}=\sigma^{ecs}=2\times10^{-3}\dunit$, $\kappa = 0\kunit$, OGSE cosine ($\delta = 14\tunit, \Delta = 14\tunit$, number of periods = 6).
The vertical bars indicate the ADC in each compartment.  
The ADC in the rightmost position is the ADC that takes into account the diffusion in all the compartments.
\label{fig:STA}}
\end{figure}

In Fig. \ref{fig:STA} we show that both BTPDE and HADC solutions match the STA values at short diffusion times 
for cylindrical cells (compartments 1 to 5).  
We also show that for the ECS (compartment 6), the STA is too low, because it does not account 
for the fact that spins in the ECS can diffuse around several cylinders.
This also shows that when the interfaces are impermeable, the BTPDE ADC and that from the 
\soutnew{H-$ADC$}{HADC} model
are identical.  The diffusion-encoding sequence here is cosine OGSE with 6 periods.

\subsection{Permeable membranes}

\begin{figure}[!htb]
  \centering
\includegraphics[width=0.45\textwidth]{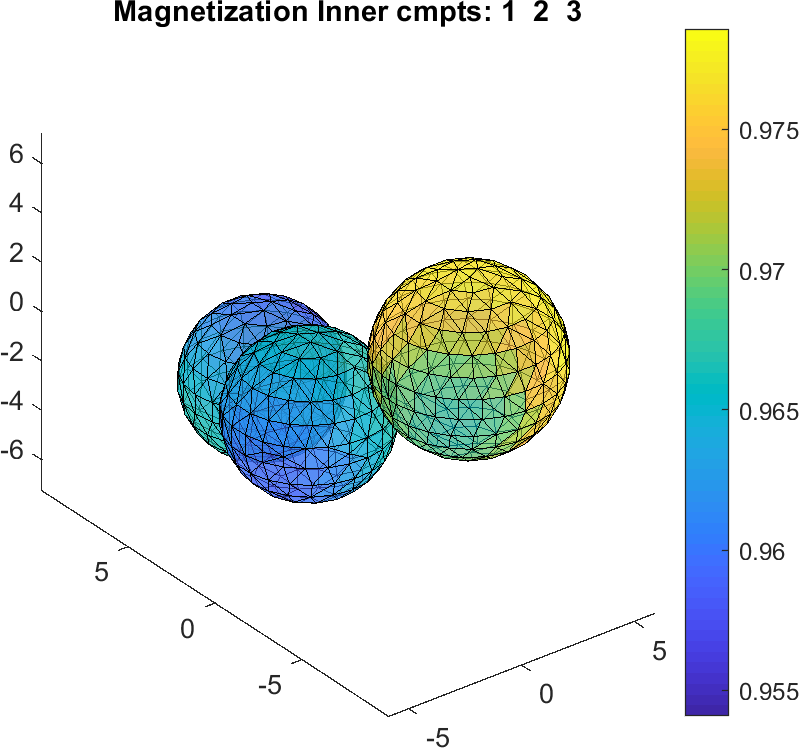} \quad \
\includegraphics[width=0.45\textwidth]{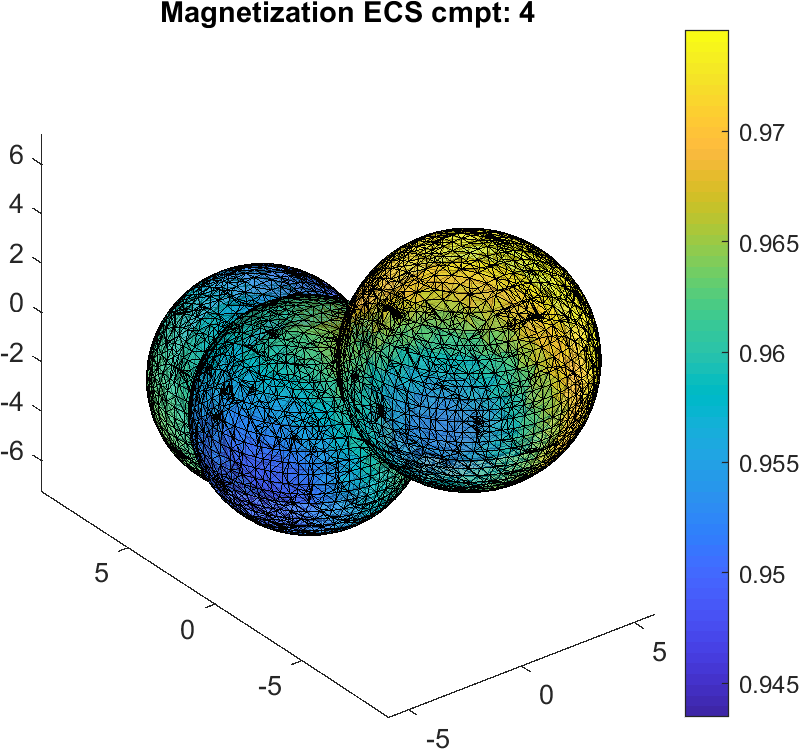}\\ 
\vspace{0.8cm}
\includegraphics[width=0.45\textwidth]{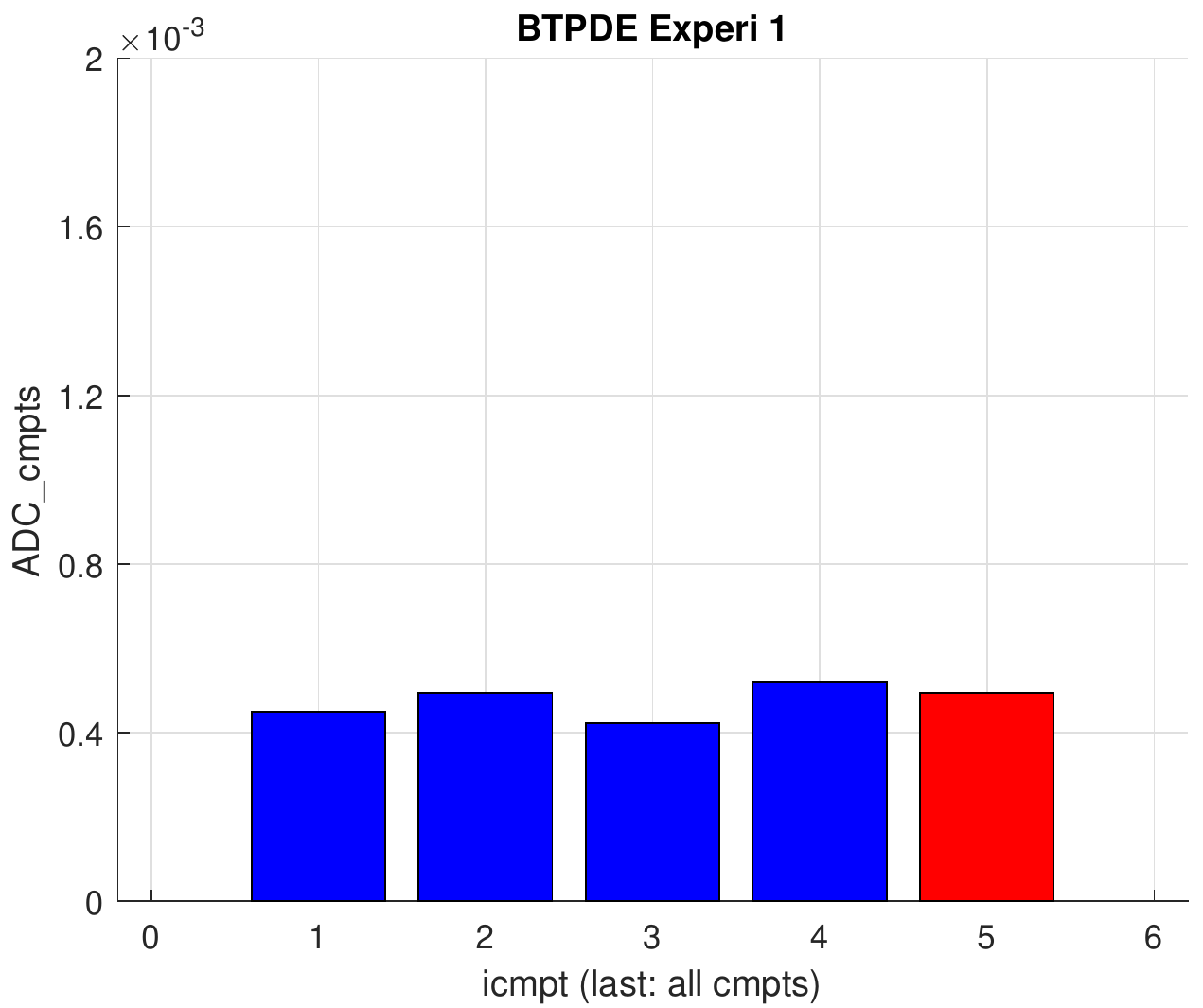} \quad 
\includegraphics[width=0.45\textwidth]{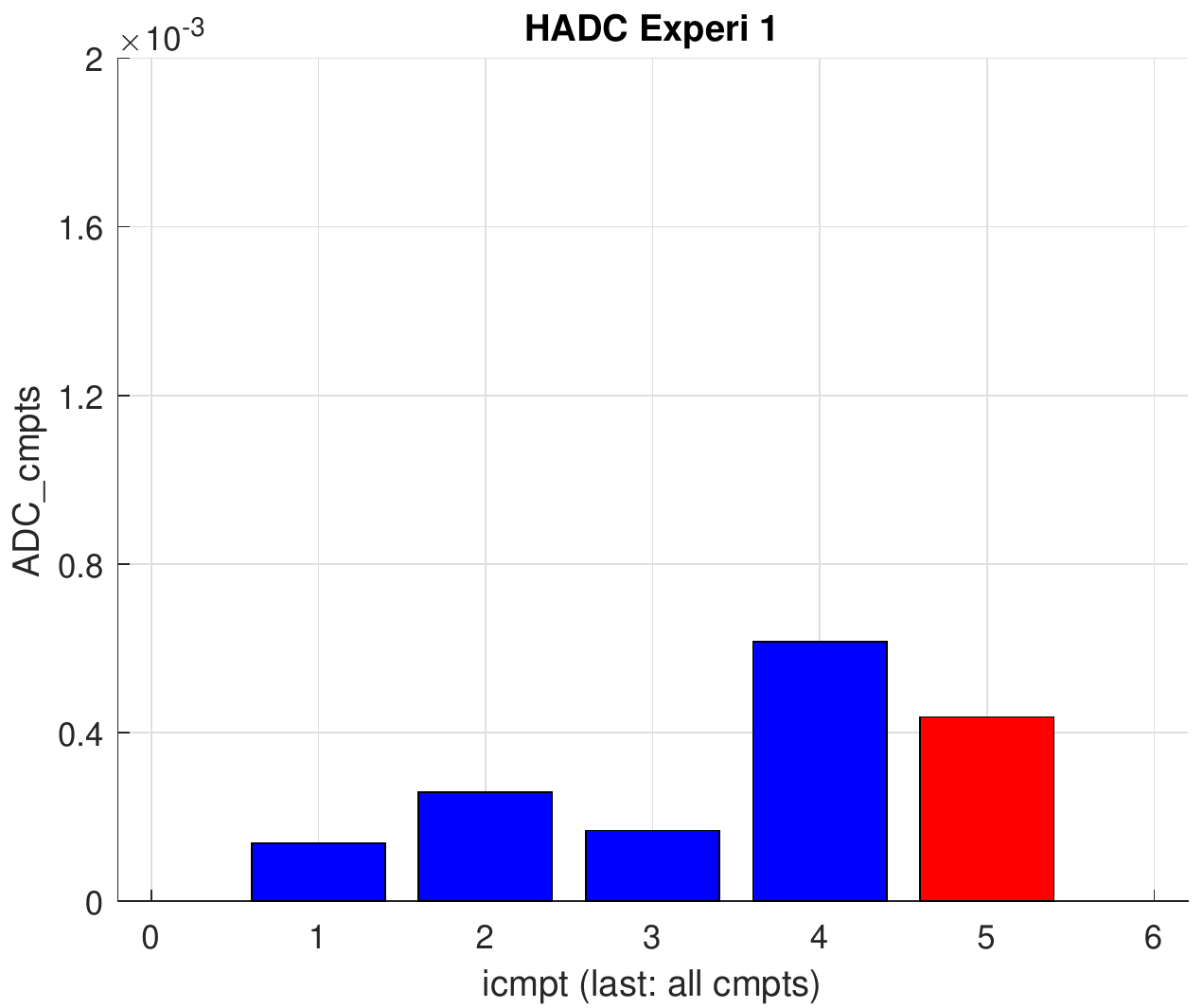}
  \label{2layer_sphere}
  \caption{Geometry: 3 spheres, tight wrap ECS, ECS gap = 0.3, $\bug = [1,1,0]$, $\sigma^{in}=\sigma^{ecs}=2\times10^{-3}\dunit$, $\kappa = 1\e{-3} \kunit$ (left), $\kappa = 0\kunit$ (right).  PGSE ($\delta = 5\tunit, \Delta = 5\tunit$).
The vertical bars indicate the ADC in each compartment.  
The ADC in the rightmost position is the ADC that takes into account the diffusion in all the compartments.
\label{fig:perm}}
\end{figure}

In Fig. \ref{fig:perm} we show the effect of permeability: the BTPDE model 
includes permeable membranes ($\kappa = 1\e{-3}\kunit$)
whereas the 
\soutnew{H-$ADC$}{HADC} 
has impermeable membranes.  We see in the permeable case, the ADC in the spheres are 
higher than in the impermeable case, whereas the ECS show reduced ADC because the faster diffusing 
spins in the ECS are allowed to moved into the slowly diffusing spherical cells.
We note that in the permeable case, the ADC in each compartment is 
obtained by using the fitting formula involving the logarithm 
of the dMRI signal, and we defined the "signal" in a compartment as the total magnetization in that compartment
at TE, which is just the integral of the solution of the BTPDE in that compartment. 

\subsection{Myelin layer}
\begin{figure}[!htb]
  \centering  
\includegraphics[width=0.45\textwidth]{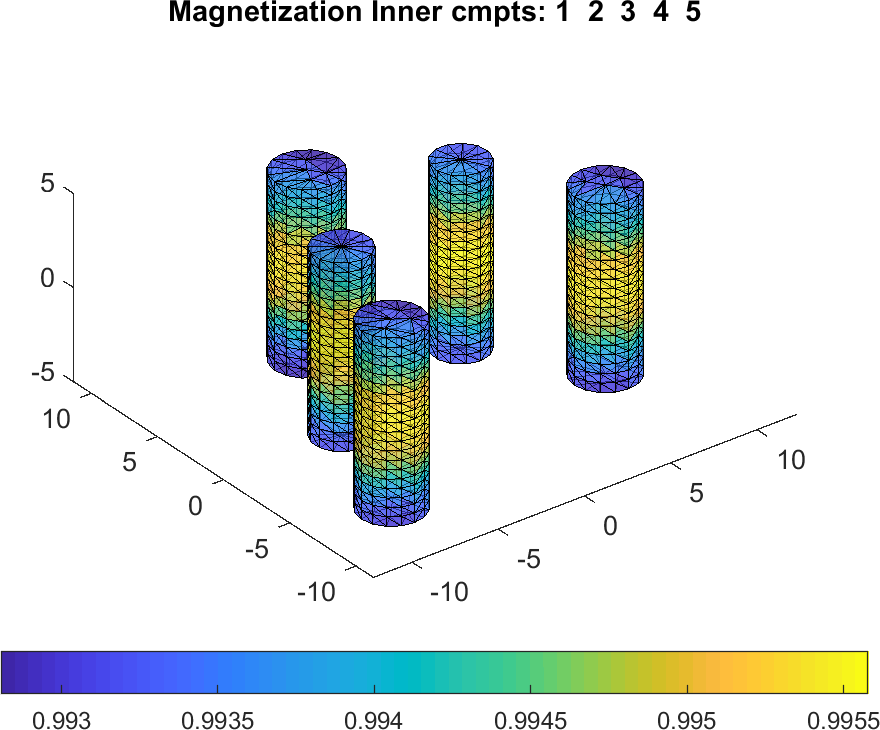}\quad
\includegraphics[width=0.45\textwidth]{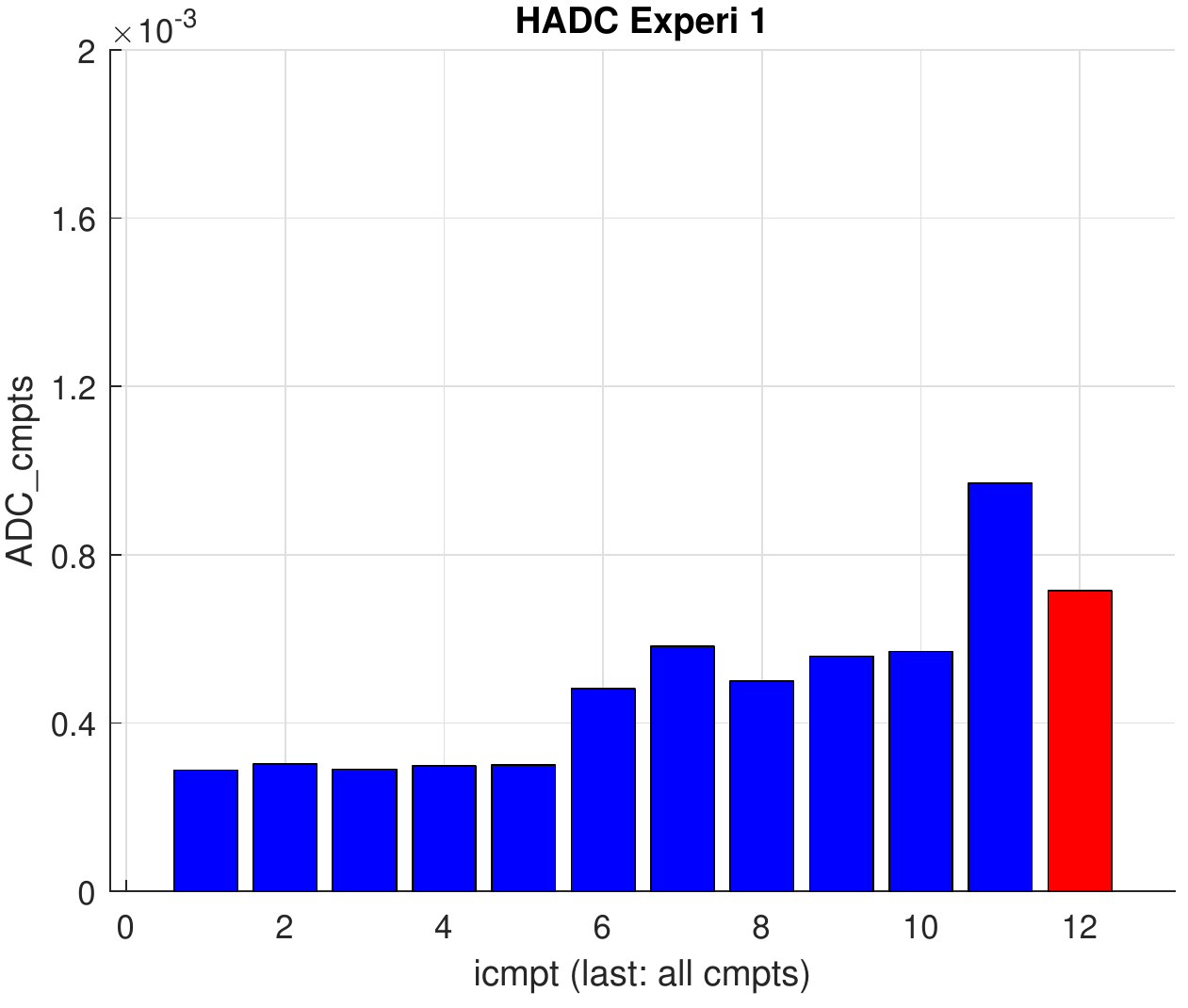}\\
\vspace{0.2cm}
\includegraphics[width=0.45\textwidth]{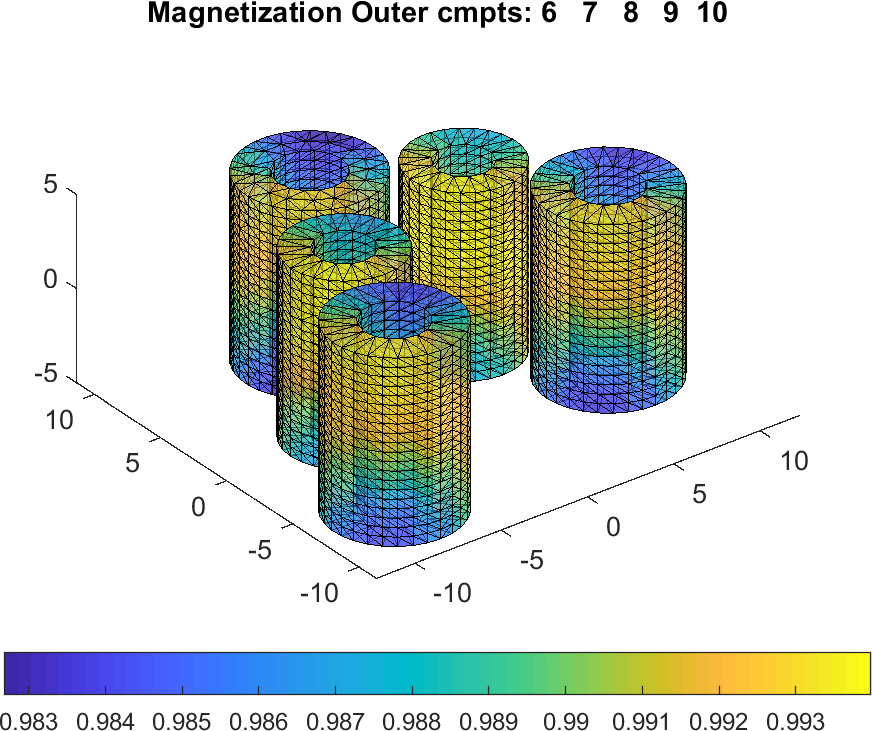} \quad
\includegraphics[width=0.45\textwidth]{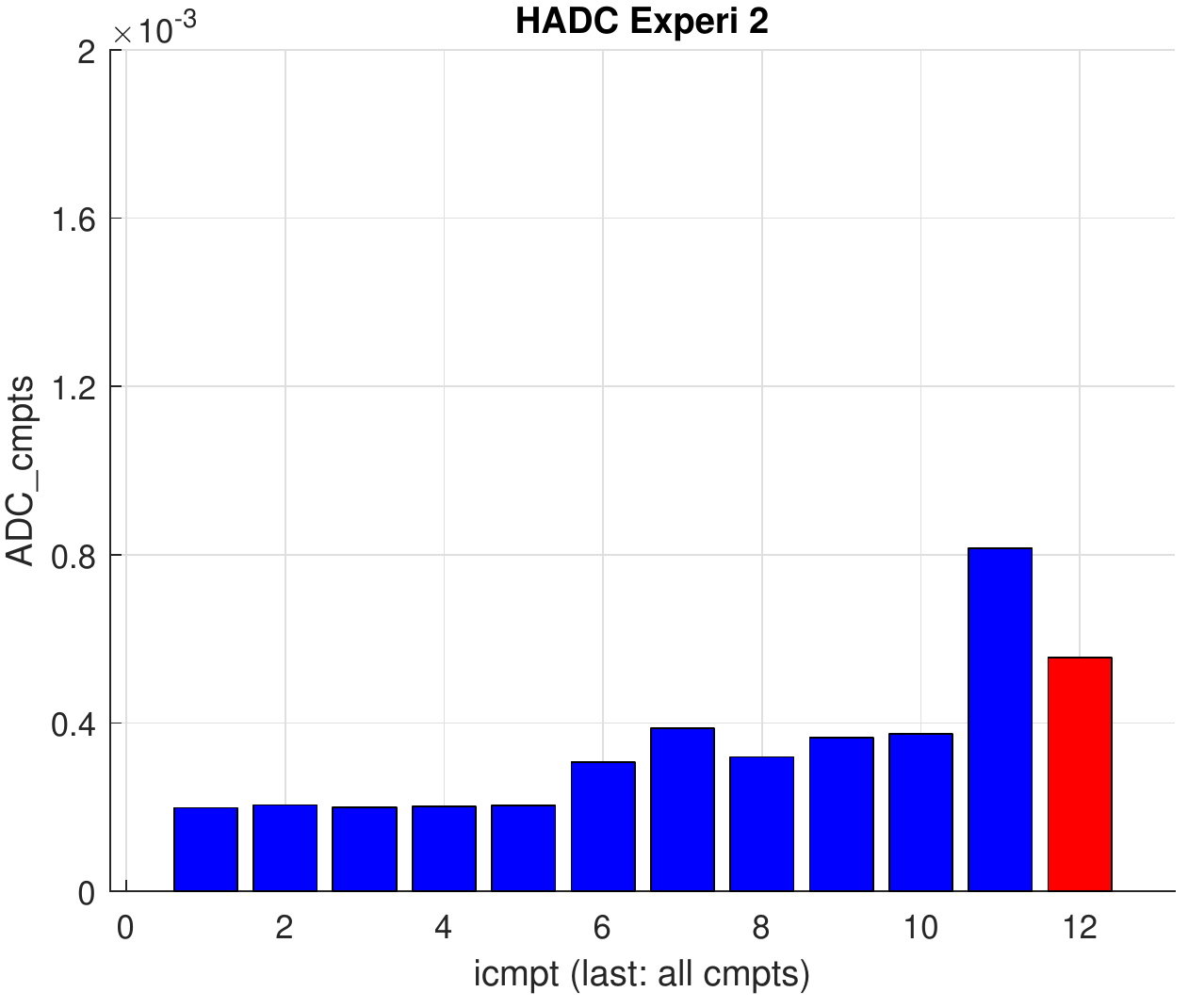}\\
\vspace{0.2cm}
\includegraphics[width=0.45\textwidth]{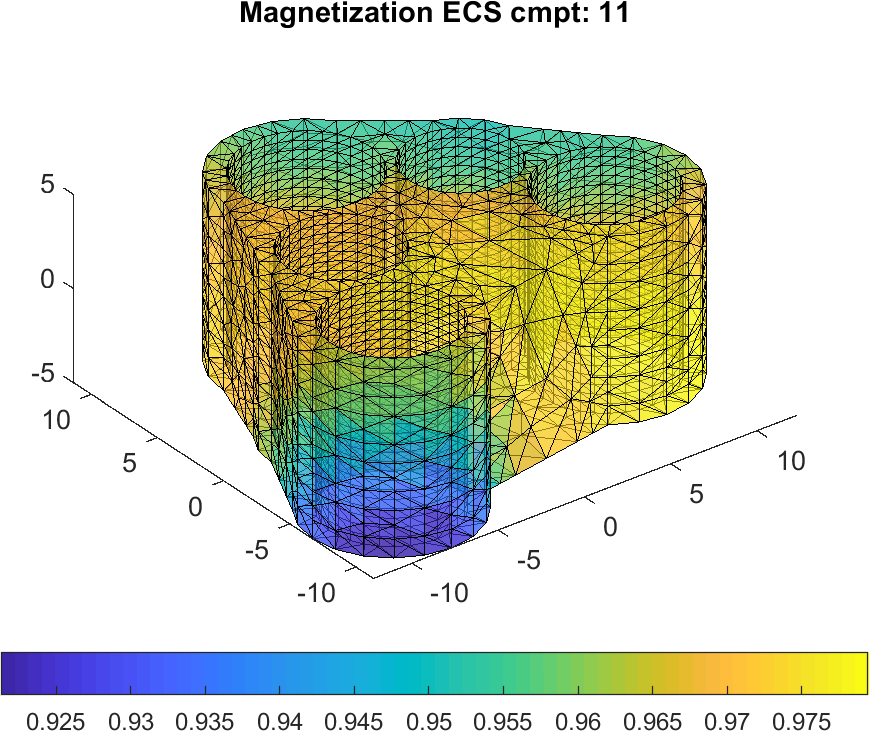} \quad
\includegraphics[width=0.45\textwidth]{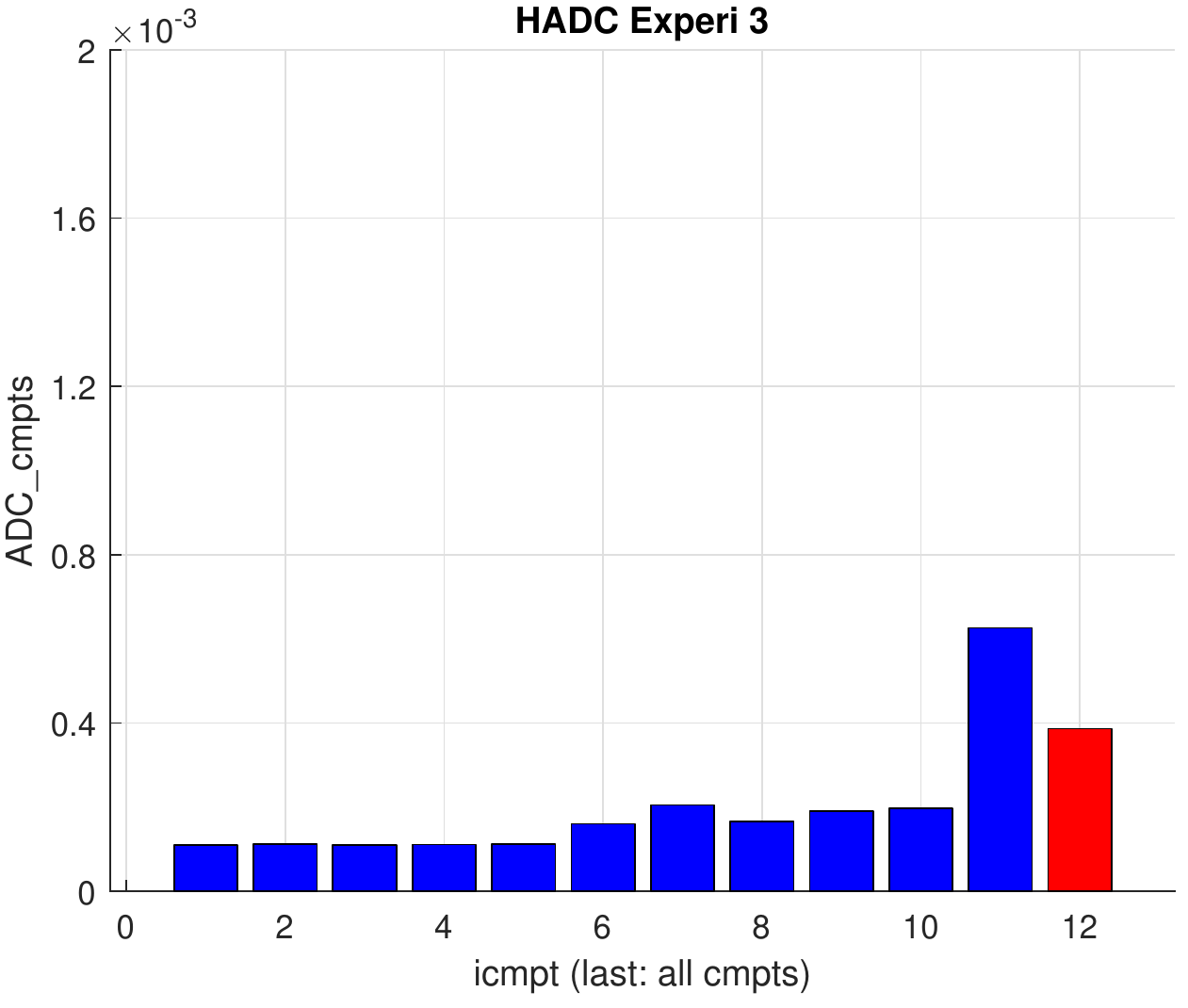}
  \caption{Geometry: 5 cylinders, myelin layer, $R_{in}/R_{out} = 0.5$, tight wrap ECS, ECS gap = 0.3, $\kappa = 0\kunit$,  $\bug = [1,1,1]$, $\sigma^{in}=\sigma^{out}=\sigma^{ecs}=2\times10^{-3}\dunit$, 3 experiments: PGSE ($\delta = 5\tunit, \Delta = 5, 10, 20\tunit$).  \soutnew{}{Left: the magnetization at $\Delta = 5\tunit$.  Right: the ADC values.}
The vertical bars indicate the ADC in each compartment.  
The ADC in the rightmost position is the ADC that takes into account the diffusion in all the compartments.
\label{fig:myelin}}
\end{figure}

In Fig. \ref{fig:myelin} we show the diffusion in cylindrical cells, the myelin layer, and the ECS. 
The ADC is higher in the myelin layer than in the cells, because for spins in the myelin layer diffusion occurs in the tangential 
direction (around the circle).  At longer diffusion times, the ADC of both the myelin layer and the cells becomes very low. 
The ADC is the highest in the ECS, because the diffusion distance
can be longer than the diameter of a cell, since the diffusing spins can move around multiple cells.  

\subsection{Twisting and bending}
\begin{figure}[!htb]
  \centering
\includegraphics[width=0.2448\textwidth]{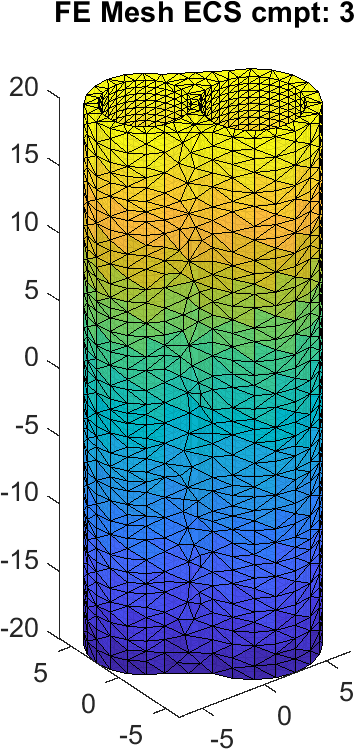}\quad\quad
\includegraphics[width=0.37\textwidth]{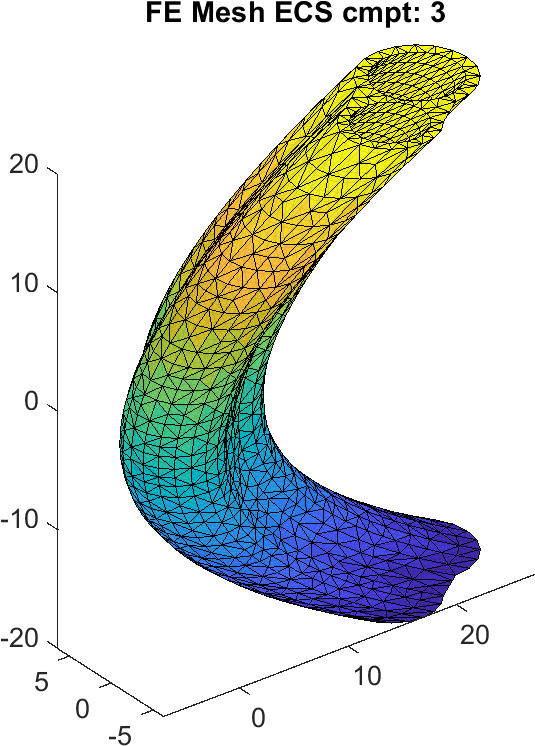}\quad\quad
\includegraphics[width=0.2725\textwidth]{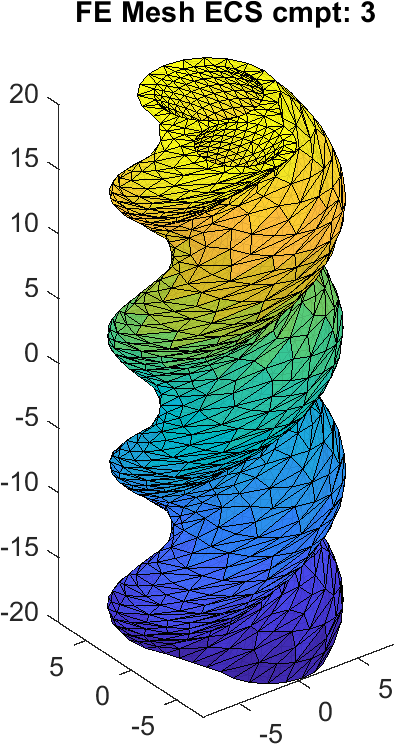}\\
\vspace{0.5cm}
\includegraphics[width=0.2874\textwidth]{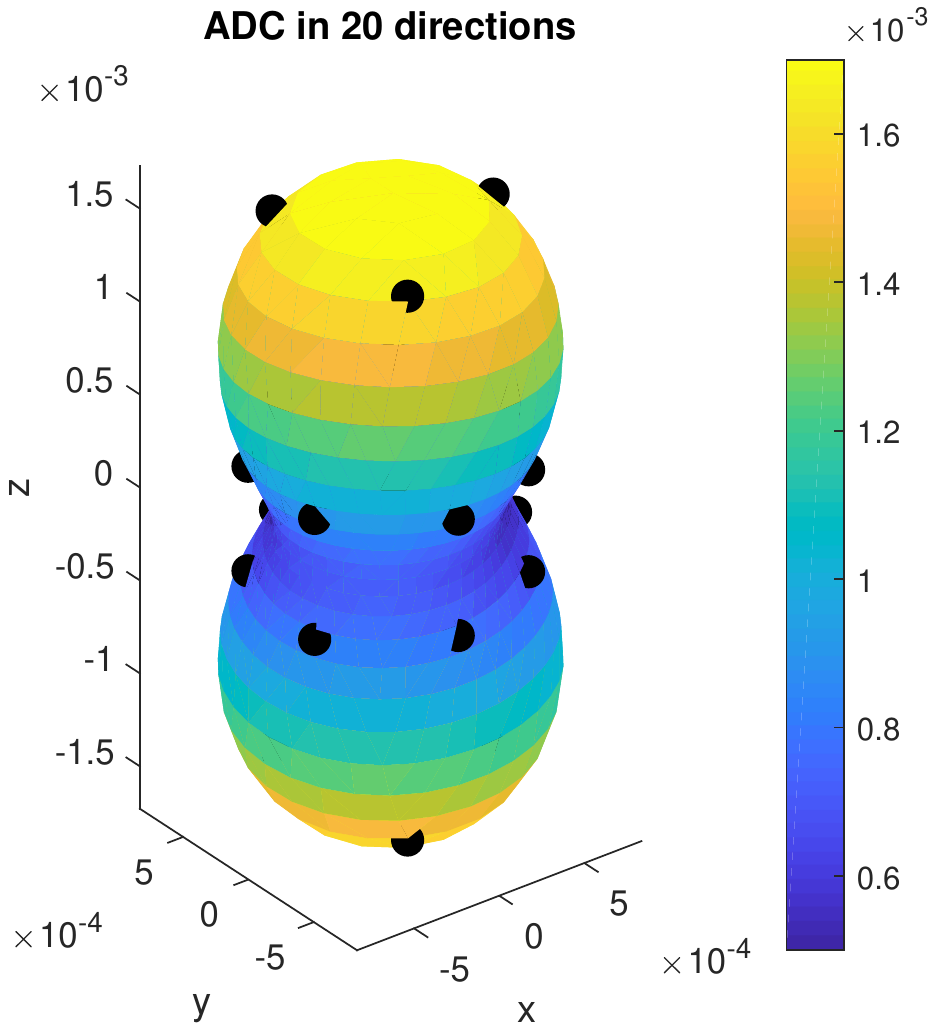} \ 
\includegraphics[width=0.3358\textwidth]{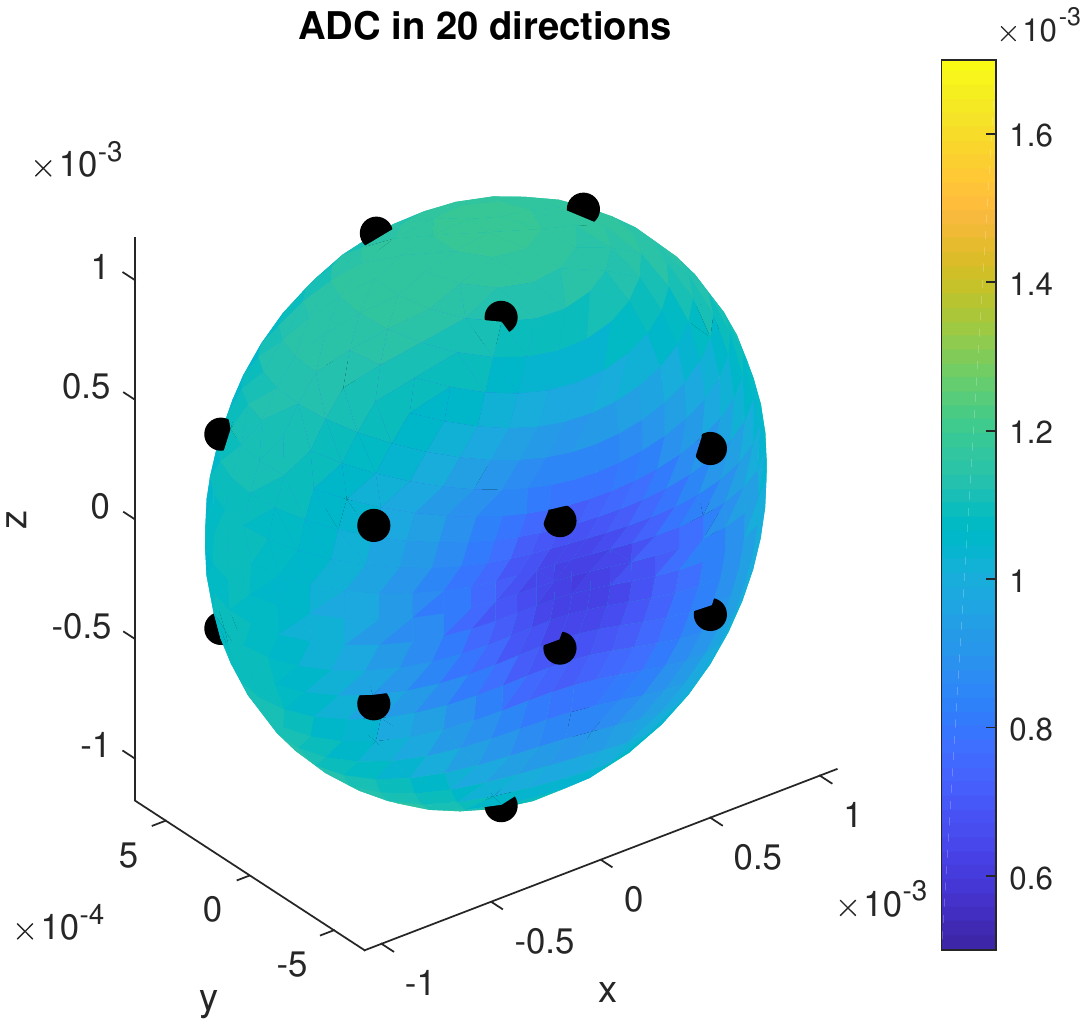} \ 
\includegraphics[width=0.3368\textwidth]{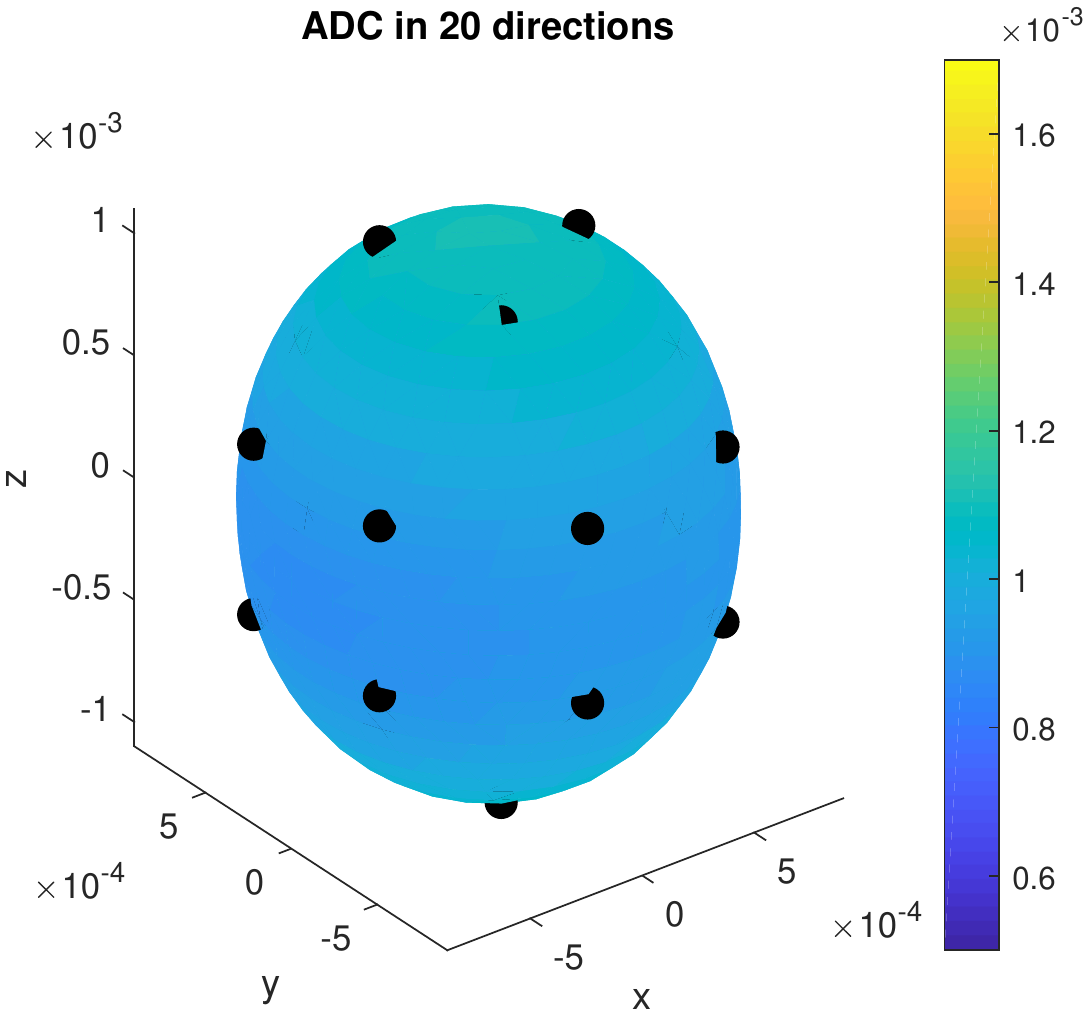}
\caption{Geometry: 2 cylinders, no myelin layer, tight wrap ECS, ECS gap = 0.3, $\kappa = 0\kunit$, $\sigma^{out}=\sigma^{ecs}=2\times10^{-3}\dunit$, PGSE ($\delta = 2.5\tunit, \Delta = 5\tunit$). \\ 
Left: canonical configuration. Middle: bend parameter = 0.05.
Right: twist parameter = 0.30. \label{fig:bend_twist}
Top: FE mesh of the ECS (the FE mesh of the axon compartments numbered 1 and 2 not shown).  
Bottom: interpolated values of the HADC on the unit sphere, and then the sphere was 
distorted to reflect the value of the HADC.  The color axis also gives the value of the HADC in the 
various gradient directions.  The black dots indicate the 20 original gradient-directions in which the HADC was 
simulated.  The spherical harmonics interpolation takes the 20 original directions into 900 directions 
uniformly distributed on the sphere.}
\end{figure}

In Fig. \ref{fig:bend_twist} we show the effect of bending and twising in cylindrical cells in multiple gradient directions.
The 
\soutnew{H-$ADC$}{HADC} 
is obtained in 20 directions uniformly distributed in the sphere.  We used spherical harmonics 
interpolation to interpolate the HADC in the entire sphere.  Then we deformed the radius of the unit sphere 
to be proportional to the interpolated HADC and plotted the 3D shape.  The color axis also indicates the value 
of the interpolated HADC.

\subsection{Timing}
In Table \ref{table:computational_times} we give the average computational times for solving the BTPDE and the HADC. All simulations were performed on \soutnew{}{a laptop computer with the processor} Intel(R) Core(TM) i5-4210U CPU @ 1.70 GHz 2.40 GHz, running Windows 10 (1809).
The geometrical configuration includes 2 axons and a tight wrap ECS, the simulated sequence is PGSE ($\delta = 2.5\tunit, \Delta = 5\tunit$).  In the impermeable case, the compartments are uncoupled, and the computational times are given separately
for each compartment.  In the permeable membrane case, the compartments are coupled, and the computational time 
is for the coupled system (relevant to the BTPDE only).
\begin{table}[!htb]
\begin{center}
\begin{tabular}{|l|c|c|c|c|}
\hline 
& \multirow{2}{*}{FE mesh size} & BTPDE &  BTPDE & \multirow{2}{*}{HADC} \\ 
&  &  $b = 50\bunit$ &   $b = 1000\bunit$ & \\ 
\hline 
 Uncoupled: Axons& 5865 nodes, 19087 ele& 7.89 sec & 9.07 sec& 8.80 sec \\ 
\hline 
Uncoupled: ECS & 6339 nodes, 19618 ele & 10.14 sec & 13.95 sec & 11.87 sec\\ 
\hline 
 Coupled: Axons+ECS & 7344 nodes, 38705 ele  & 39.14 sec &43.24 sec &N\slash A\\ 
\hline 
\end{tabular} 
\caption{Computational times for solving the BTPDE and the HADC. All simulations were performed on Intel(R) Core(TM) i5-4210U CPU @ 1.70 GHz 2.40 GHz, running Windows 10 (1809).
The geometrical configuration includes 2 axons and a tight wrap ECS, the simulated sequence is PGSE ($\delta = 2.5\tunit, \Delta = 5\tunit$).\label{table:computational_times}}
\end{center}
\end{table}

\section{Numerical validation of SpinDoctor}
\marginparnew{New section added in revised version}
\label{sec:mf}

In this section, we validate SpinDoctor by comparing SpinDoctor with the Matrix Formalism 
method \cite{Callaghan1997,Barzykin1999} in a simple geometry.
The Matrix Formalism method is a closed form representation of the dMRI signal 
based on the eigenfunctions of the Laplace operator subject to homogeneous Neumann boundary conditions.
These eigenfunctions are available in explicit form for elementary geometries such as the line segment, 
the disk, and the sphere \cite{Grebenkov2007,Ozarslan2009,Drobnjak2011a,Grebenkov2010a}.
The dMRI signal obtained using the Matrix Formalism method will be considered the reference solution in this section.

The accuracy of the SpinDoctor simulations can be tuned using three simulation parameters:
\begin{enumerate}
\item $Htetgen$ controls the finite element mesh size;
\begin{enumerate}
\item $Htetgen = -1$ means the FE mesh size is determined automatically by the internal algorithm of Tetgen to ensure a good quality mesh (subject to the constraint that the radius to edge ratio of tetrahedra is no larger than 2.0).
\item $Htetgen = h$  requests a desired FE mesh tetrahedra height of $h$ $\lunit$ (in later versions of Tetgen, this parameter has been changed to the desired volume of the tetrahedra).  
\end{enumerate}
\item $rtol$ controls the accuracy of the ODE solve.  It is the relative residual tolerance at all points of the FE mesh at each time step of the 
ODE solve;
\item $atol$ controls the accuracy of the ODE solve.  It is the absolute residual tolerance at all points of the FE mesh at each time step of the 
ODE solve;
\end{enumerate}
We varied the finite element mesh size and the ODE solve accuracy of SpinDoctor and ran
6 simulations with the following simulation parameters: 
\begin{enumerate}[label=SpinD Simul \ref{sec:mf}-\arabic*:, wide=0pt, font=\textbf]
\item $rtol = 10^{-3}$, $atol = 10^{-6}$, $Htetgen = -1$;
\item $rtol = 10^{-6}$, $atol = 10^{-9}$, $Htetgen = -1$;
\item $rtol = 10^{-9}$, $atol = 10^{-12}$, $Htetgen = -1$;
\item $rtol = 10^{-3}$, $atol = 10^{-6}$, $Htetgen = 1$;
\item $rtol = 10^{-6}$, $atol = 10^{-9}$, $Htetgen = 1$;
\item $rtol = 10^{-9}$, $atol = 10^{-12}$, $Htetgen = 1$;
\end{enumerate}

The geometry simulated is the following:
\begin{itemize}
\item {\it 3LayerCylinder} is a 3-layer cylindrical geometry of height $1\lunit$ and 
the layer radii, $R_1 = 2.5\lunit$, $R_2 = 5\lunit$ and $R_3 = 10\lunit$.  The middle layer is subject to permeable interface conditions on both the interior and the exterior interfaces,
with permeability coefficient $\kappa$.  The exterior boundary $R = R_3$ is subject to impermeable boundary conditions.  The top and bottom boundaries are also subject to impermeable boundary conditions.
\item For this geometry, $Htetgen = -1$ gives finite elements mesh size ($n_{nodes} = 440, n_{elem} = 1397$).
$Htetgen = 1$ gives finite elements mesh size ($n_{nodes} = 718, n_{elem} = 2088$).
\end{itemize}

The dMRI experimental parameters are the following:
\begin{itemize}
\item the diffusion coefficient in all compartments is $2 \times 10^{-3}\dunit$;
\item the diffusion-encoding sequence is PGSE ($\delta = 10\tunit$, $\Delta = 13\tunit$);
\item 8 b-values: $b =\{ 0, 100, 500, 1000, 2000, 3000, 6000, 10000\}\bunit$;
\item 1 gradient direction: $[1,1,0]$;
\end{itemize}

\begin{figure}[!htb]
 \centering
\includegraphics[width=0.49\textwidth]{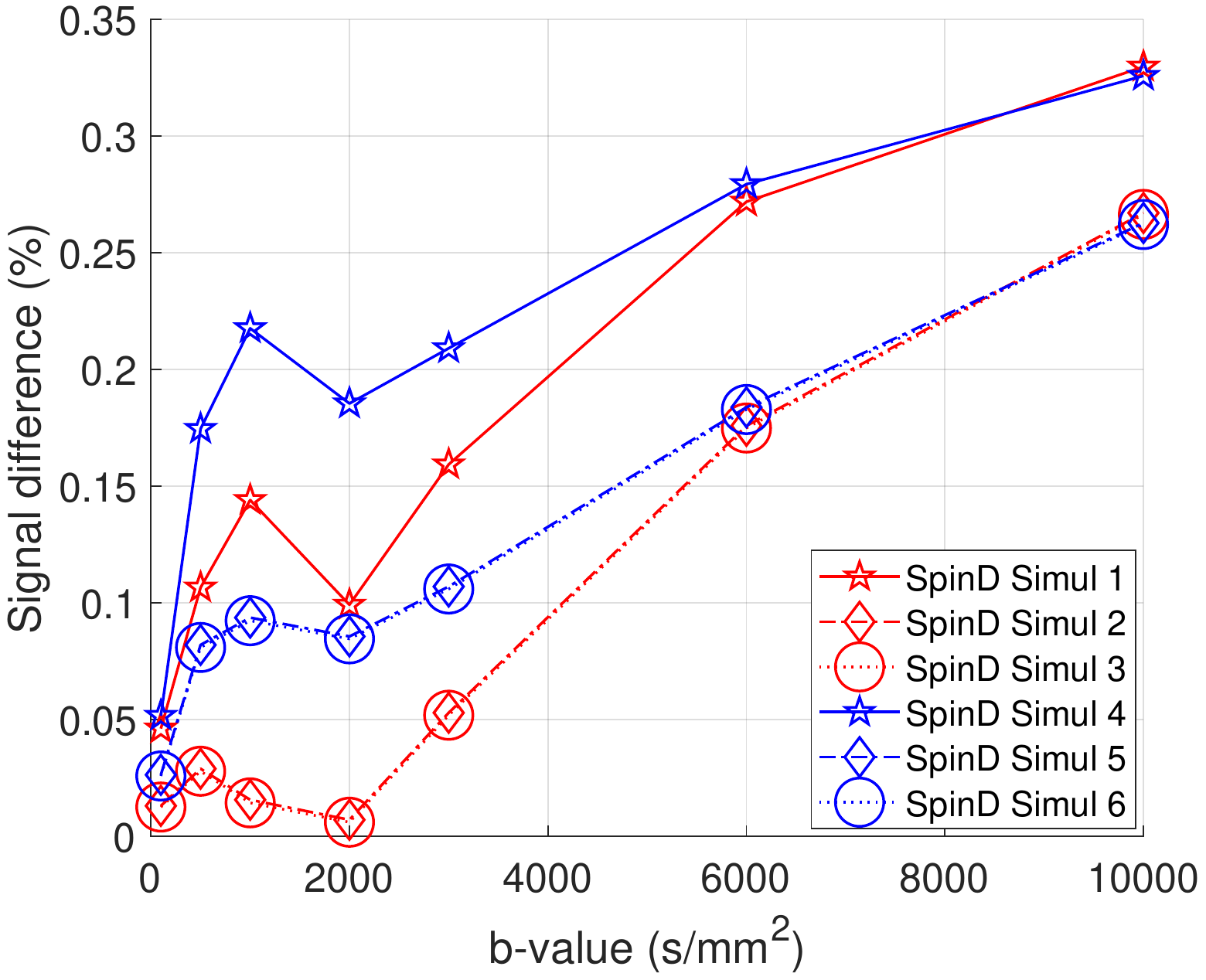}
\includegraphics[width=0.49\textwidth]{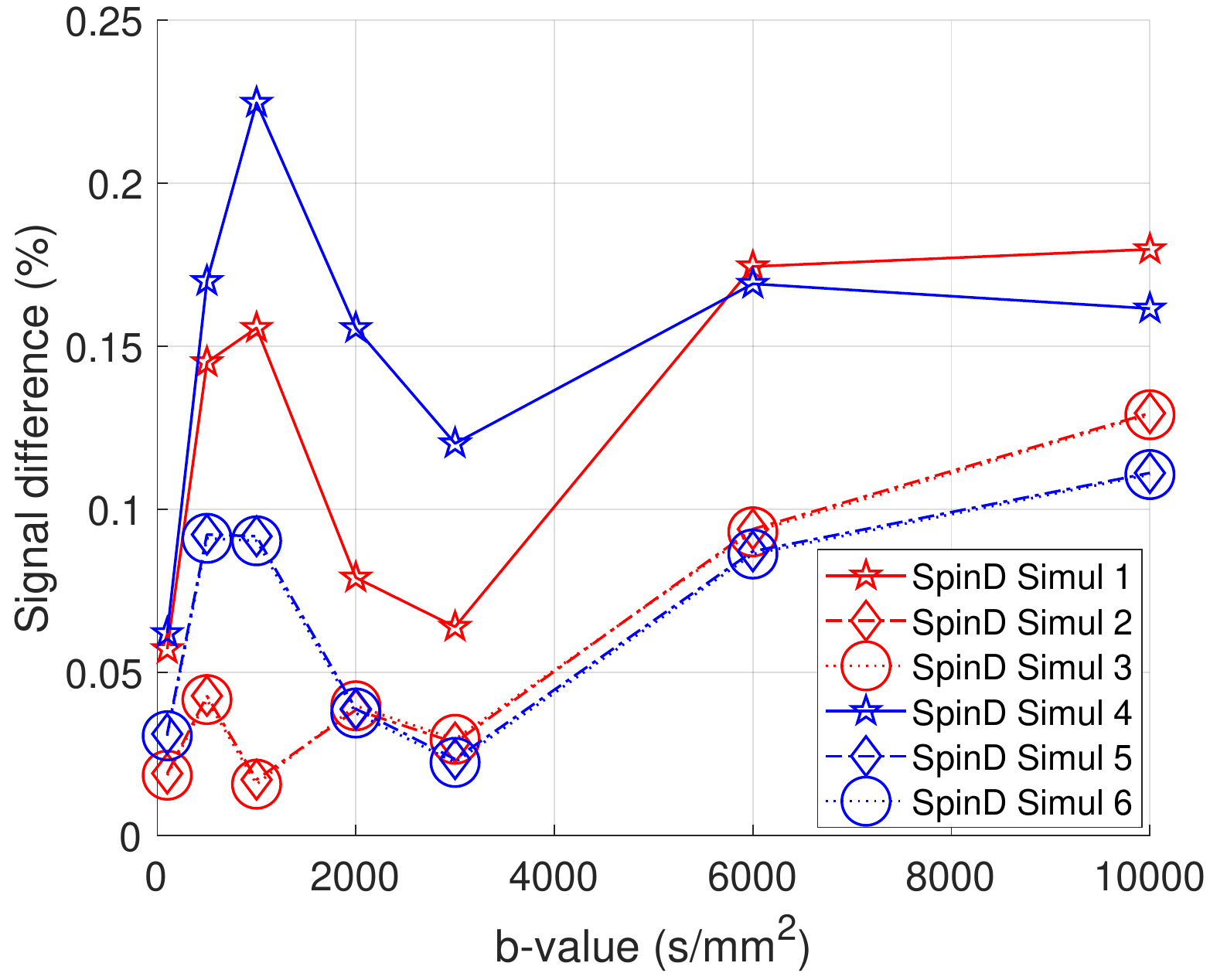}
 \caption{Signal difference between the Matrix Formalism signal (reference) and the SpinDoctor signal.  
Left: $\kappa = 10^{-5}\kunit$.  Right: $\kappa = 10^{-4}\kunit$.   
The geometry is {\it 3LayerCylinder}.  The diffusion coefficient in all compartments is $2 \times 10^{-3}\dunit$;
the diffusion-encoding sequence is PGSE ($\delta = 10\tunit$, $\Delta = 13\tunit$);
Simul 1: $rtol = 10^{-3}$, $atol = 10^{-6}$, $Htetgen = -1$;
Simul 2: $rtol = 10^{-6}$, $atol = 10^{-9}$, $Htetgen = -1$;
Simul 3: $rtol = 10^{-9}$, $atol = 10^{-12}$, $Htetgen = -1$;
Simul 4: $rtol = 10^{-3}$, $atol = 10^{-6}$, $Htetgen = 1$;
Simul 5: $rtol = 10^{-6}$, $atol = 10^{-9}$, $Htetgen = 1$;
Simul 6: $rtol = 10^{-9}$, $atol = 10^{-12}$, $Htetgen = 1$;
\label{fig:validation_mf}}
\end{figure}

In Figure \ref{fig:validation_mf} we show the signal differences (in percent) of the reference Matrix 
Formalism method and the SpinDoctor simulations, normalized by the reference signal at $b=0$:
\be{}
E(b) = \frac{\left |S^{MF}(b)-S^{SpinD}(b)\right |}{S^{MF}(b=0)}\times 100.
\ee
We see that the signal difference is less than $0.35\%$ for  $\kappa = 10^{-5}\kunit$  and 
it is less than $0.25\%$ for  $\kappa = 10^{-4}\kunit$ for all 6 SpinDoctor simulations.  
The signal difference becomes smaller when the ODE solve tolerances are changed from  ($rtol = 10^{-3}$, $atol = 10^{-6}$) to ($rtol = 10^{-6}$, $atol = 10^{-9}$), but there is no change when the tolerances are further reduced to  
($rtol = 10^{-9}$, $atol = 10^{-12}$).
If we refine the FE mesh, but keep the ODE solve tolerances the same, the signal difference is in fact larger using the refined mesh  than using the coarse mesh at the smaller b-values, though this effect disappears at higher b-values
and larger permeability.  This is probably due to parasitic oscillatory modes on the finer mesh that need smaller time steps to be sufficiently damped.

\section{\dcom{Computational time and comparison with Monte-Carlo simulation}}
\marginparnew{New section added in revised version}

In this section, we compare SpinDoctor with Monte-Carlo simulation 
using the publicly available software package Camino Diffusion MRI Toolkit \cite{Hall2009},
downloaded from \url{http://cmic.cs.ucl.ac.uk/camino}.
All the simulations were performed on a server computer with 12 processors (Intel (R) Xeon (R) E5-2667 @2.90 GHz), 192 GB of RAM, running CentOS 7.  SpinDoctor was run using MATLAB R2019a on the same computer.  

We give SpinDoctor computational times for three relatively complicated geometries.  
We also give Camino computational times for the first two geometries.   
We did not use Camino for the third geometry due to the excessive time required by Camino.  

The number of the degrees of freedom in the SpinDoctor simulations is the finite element mesh size (the number of nodes and the number of elements).  For Camino it is the number of spins.  
The time stepping choice of the SpinDoctor simulations is given by the ODE solve tolerances.
For Camino it is given by the number of time steps.
Camino has an initialization step where it places the spins and we give the time of 
this initialization step separately from the Camino random walk simulation time.

\begin{figure}[!htb]
  \centering
\includegraphics[width=0.9\textwidth]{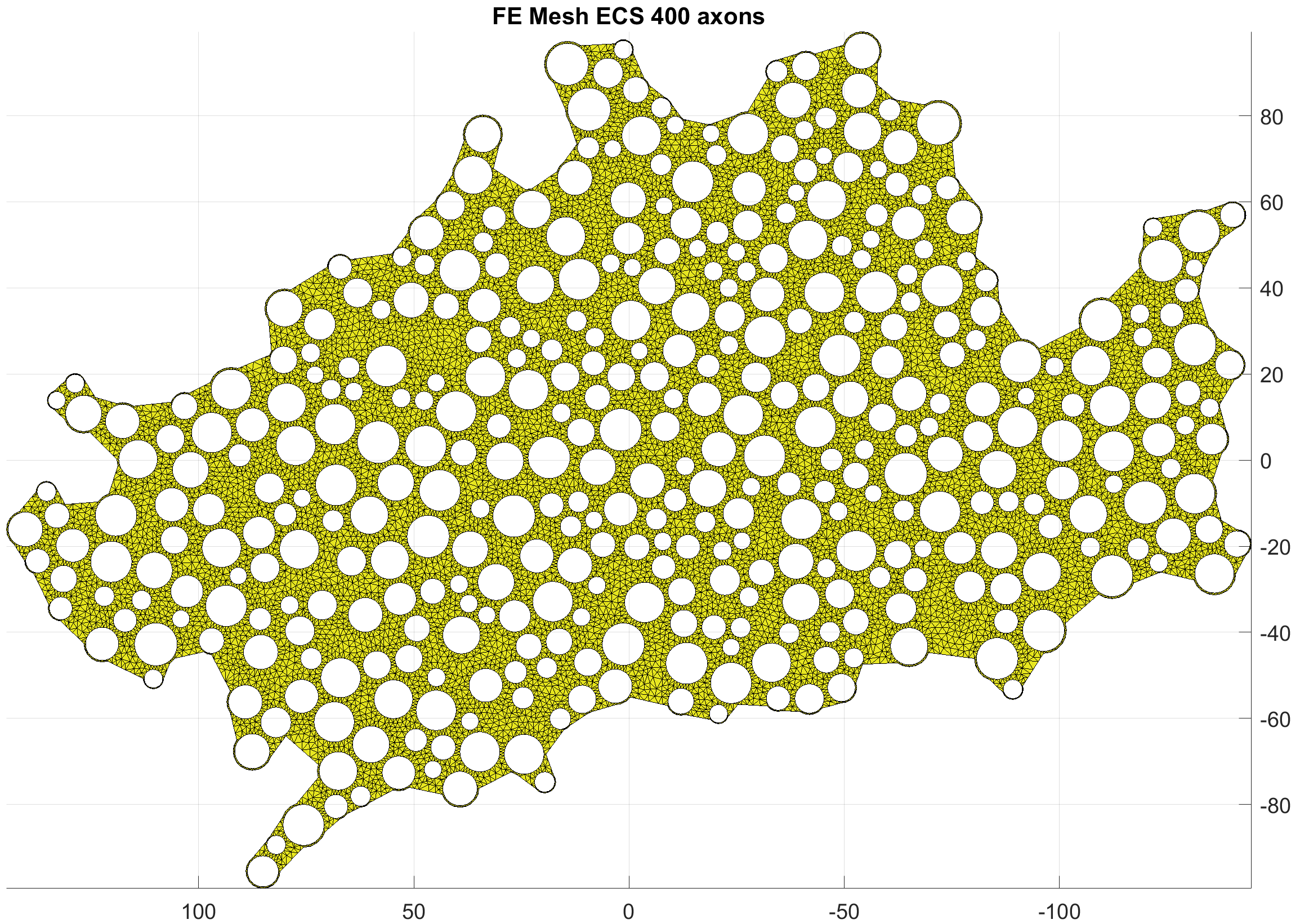} 
  \caption{\label{fig:FE_ecs400axons}  The geometry is {\it ECS400axons}.  This finite elements mesh size is ($n_{nodes} = 53280, n_{elem} = 125798$).
  }
\end{figure}

\begin{figure}[!htb]
  \centering
\includegraphics[width=0.9\textwidth]{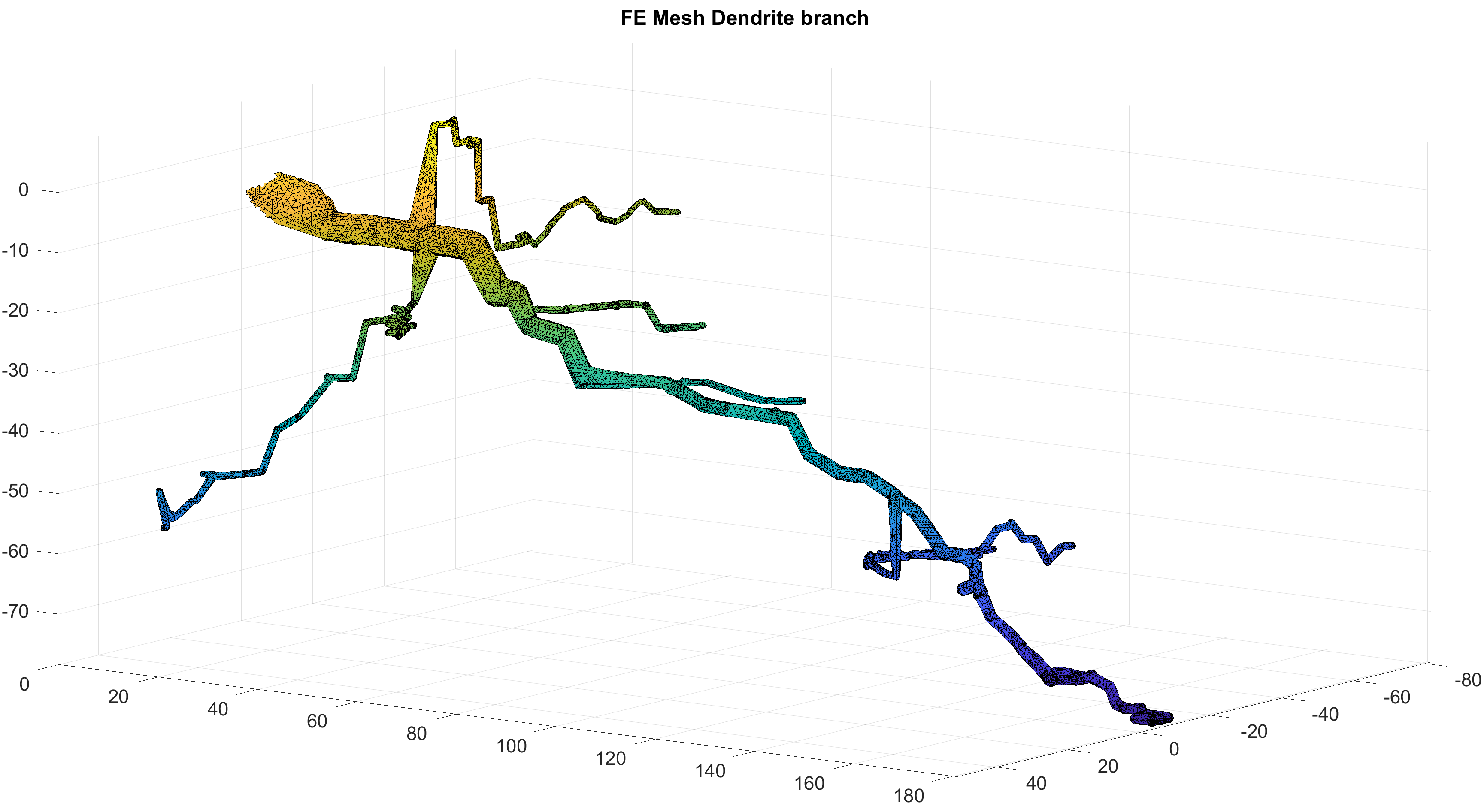} 
  \caption{ \label{fig:FE_dendritebranch}  The geometry is {\it DendriteBranch}.  This finite elements mesh size is ($n_{nodes} = 24651, n_{elem} = 91689$)
}
\end{figure}

\begin{figure}[!htb]
  \centering
\includegraphics[width=0.9\textwidth]{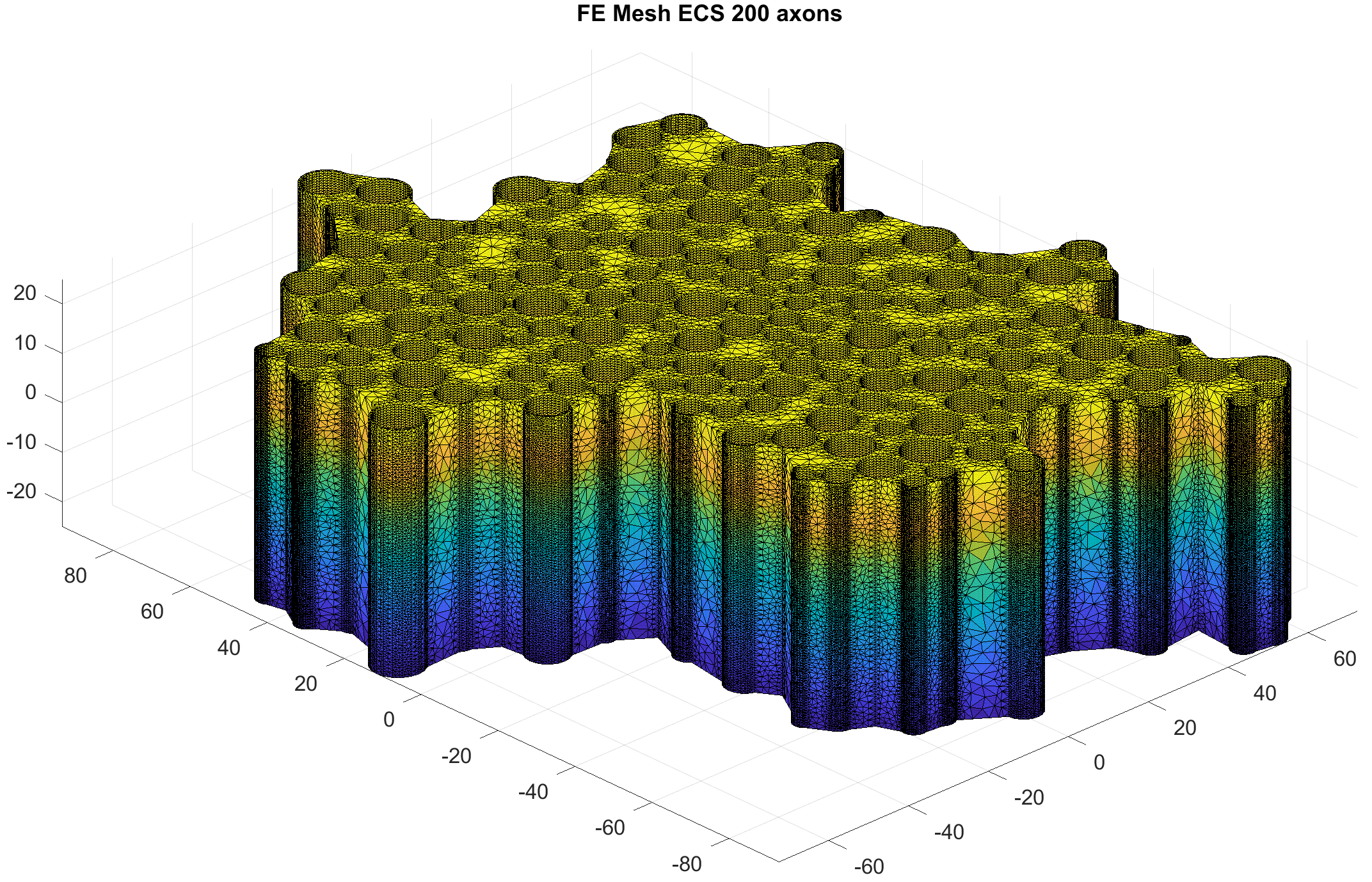}
\caption{\label{fig:FE_ecs200axons}  The geometry is {\it ECC200axons}.  This finite elements mesh size is ($n_{nodes} = 846298, n_{elem} = 2997386$)}
\end{figure}

Given the interest of the dMRI community in the extra-cellular space \cite{Burcaw2015} and 
neuron simulations, we chose the following three geometries:
\begin{enumerate}
\item {\it ECS400axons}.  See Figure \ref{fig:FE_ecs400axons}.  This models the extra-cellular space outside of 400 axons.  We generated 400 cylinders with height $1\lunit$ and radii ranging from $2-5\lunit$, randomly placed according 
to Algorithm \ref{algo:create_cells}.  The small height of the cylinders 
means that this geometry should be used only for studying transverse diffusion.
We used a tight-wrap ECS: this choice means we do not need to have a complicated algorithm to avoid large empty spaces as would be the case when the ECS is box-shaped. 
\item {\it DendriteBranch}.  See Figure \ref{fig:FE_dendritebranch}.  This is a dendrite branch whose original morphological reconstruction SWC file published in NeuroMorpho.Org \cite{Ascoli2007}.  By wrapping the geometry described in the SWC file in a new watertight surface and using the external FE meshing package GMSH \cite{Geuzaine2009}, we created a FE mesh for this dendrite branch. 
The FE mesh was in imported and used in SpinDoctor.   We note this is an externally generated FE mesh and this illustrate the capacity of SpinDoctor to simulate the dMRI on general geometries provided by the user.
\item {\it ECS200axons}.  See Figure \ref{fig:FE_ecs200axons}.  This models the extra-cellular space outside of 200 axons.  
To study 3-dimensional diffusion, the height of the cylinders was increased to $50\lunit$. 
To keep the finite element mesh size reasonable, we decreased the number of axons to 200, keeping the 
range of radii between $2-5$ microns, placed randomly as above, with a tight-wrap ECS.  
 
\end{enumerate}

The dMRI experimental parameters are the following:
\begin{itemize}
\item the diffusion coefficient is $2 \times 10^{-3}\dunit$;
\item the diffusion-encoding sequence is PGSE ($\delta = 10\tunit$, $\Delta = 13\tunit$);
\item 8 b-values: $b =\{ 0, 100, 500, 1000, 2000, 3000, 6000, 10000\}\bunit$;
\item 1 gradient direction: $[1,1,0]$.
\end{itemize}

The SpinDoctor simulations were done using one compartment.  The boundary of 
compartment is subject to impermeable boundary conditions.  We took the surface triangulations associated with the finite element mesh for the SpinDoctor simulations and used them as the input PLY files for Camino.  Camino is called with the command \texttt{datasynth}.  The options of Camino that are relevant to the simulations in the above three geometries are the following: 
\begin{itemize}
\item \texttt{-walkers \$\{N\}}: $N$ is the number of walkers ;
\item \texttt{-tmax \$\{T\}}: $T$ is the number of time steps;
\item \texttt{-p \$\{P\}}: $P$ is the probability that a spin will step through a barrier.  We set $P$ to zero;
\item \texttt{-voxels 1}: using 1 voxel for the experiment;
\item \texttt{-initial intra}: random walkers are placed uniformly inside the geometry and none outside of it;
\dcom{In the case of the extra-cellular space, \texttt{intra} means inside the geometry, with the geometry representing the extracellular space};
\item  \texttt{-voxelsizefrac 1}: the signal is computed using all the spins inside the geometry described by the PLY file, and not just in a center region;
\item \texttt{-diffusivity 2E-9}: the diffusion coefficient ($m^2/s$); 
\item \texttt{-meshsep \$\{xsep\} \$\{ysep\} \$\{zsep\}}: specifies the seperation between bounding box for mesh substrates.  We used a box that fully contains the geometry described by the PLY file;
\item \texttt{-substrate ply}: mesh substrates are constructed using a PLY file;
\item \texttt{-plyfile \$\{plyfile\}}:  the name of the PLY file.  We wrote a MATLAB function that outputs the list of triangles that make up the boundary of the finite element mesh and formatted it as a PLY file.  We note these triangles form a surface triangulation;
\end{itemize}

\subsection{ECS of 400 axons}
\label{sec:ecs400}
SpinDoctor was run with the following 3 sets of simulation parameters: 
\begin{enumerate}[label=SpinD Simul \ref{sec:ecs400}-\arabic*:, wide=0pt, font=\textbf]
\item  $rtol = 10^{-3}$, $atol = 10^{-6}$, $Htetgen = 0.5$; 
\item  $rtol = 10^{-3}$, $atol = 10^{-6}$, $Htetgen = 1$;
\item  $rtol = 10^{-3}$, $atol = 10^{-6}$, $Htetgen = -1$;
\end{enumerate}
 For this geometry, $Htetgen = -1$ gives finite elements mesh size ($n_{nodes} = 53280, n_{elem} = 125798$).   $Htetgen = 1$ gives finite elements mesh size ($n_{nodes} = 58018, n_{elem} = 139582$). $Htetgen = 0.5$ gives finite elements mesh size ($n_{nodes} = 70047, n_{elem} = 177259$).

Camino was run with the following 2 sets of simulation parameters: 
\begin{enumerate}[label=Camino Simul \ref{sec:ecs400}-\arabic*:, wide=0pt, font=\textbf]
\item  $N = 1000$, $T = 200$;
\item  $N = 4000$, $T = 800$;
\end{enumerate} 

The reference signals are \textbf{SpinD Simul \ref{sec:ecs400}-1}, the SpinDoctor signals computed on the finest FE mesh ($Htetgen = 0.5$).

\begin{figure}[!htb]
  \centering
\includegraphics[width=0.49\textwidth]{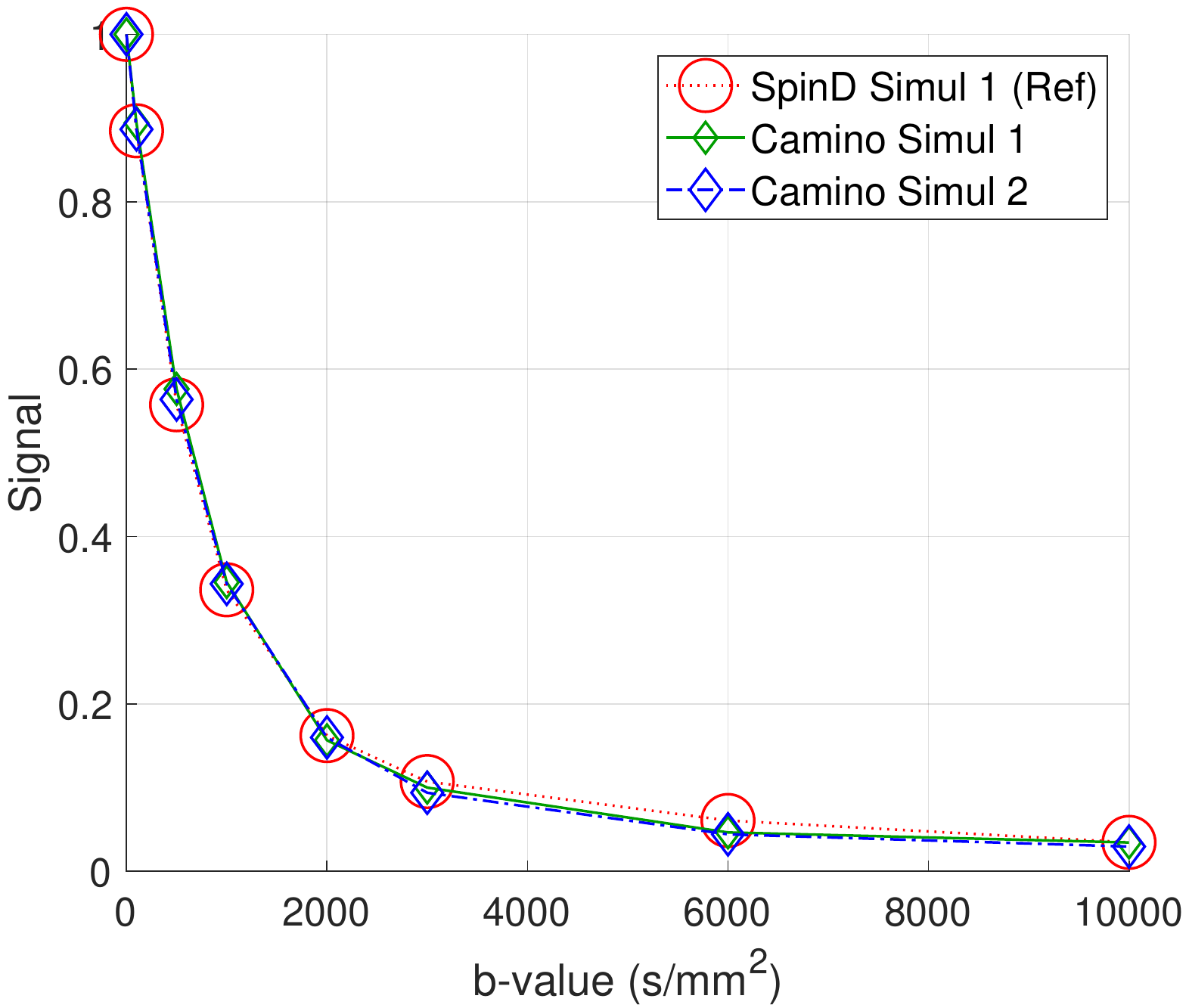}\\
 \includegraphics[width=0.49\textwidth]{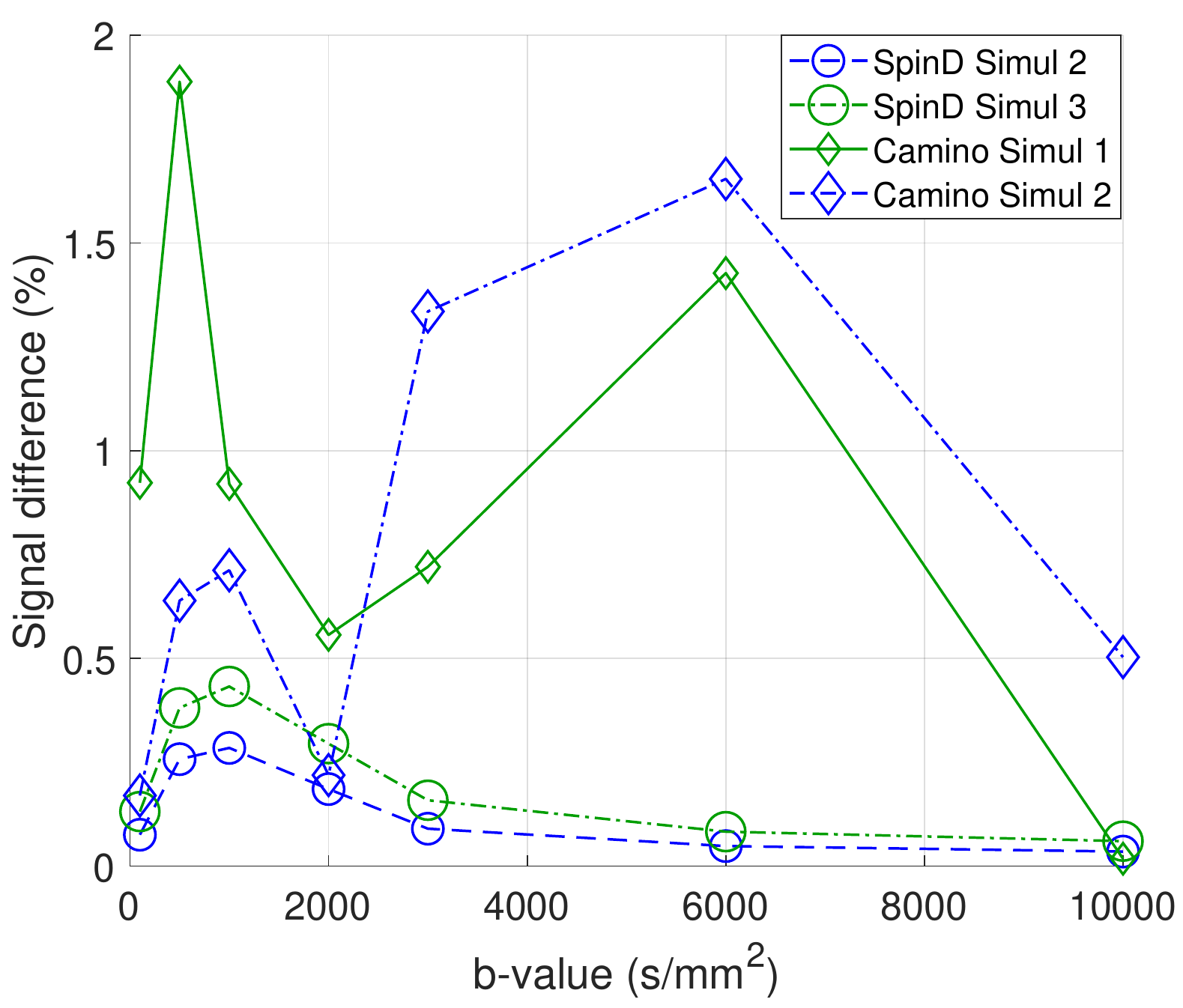}
\includegraphics[width=0.49\textwidth]{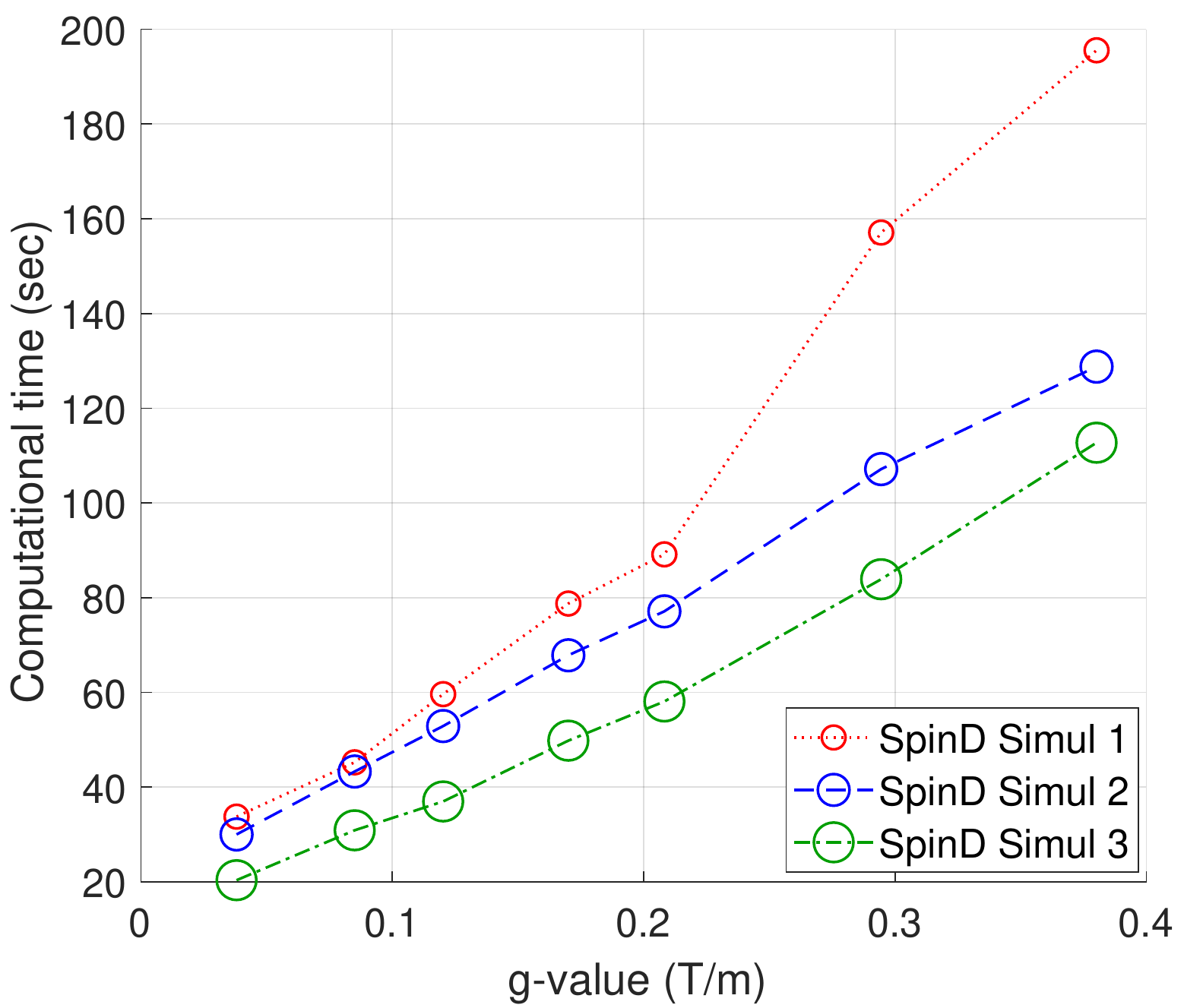}
  \caption{\label{fig:ecs400}The geometry is {\it ECS400axons}.  Top: SpinD Simul 1 is the reference signal, compared to two Camino simulations.  Bottom left: the signal difference between the reference simulation and two SpinDoctor simulations and two Camino simulations.  Bottom right: the computational times of SpinDoctor simulations as a function of the gradient amplitude.  The diffusion coefficient is $2 \times 10^{-3}\dunit$;
The diffusion-encoding sequence is PGSE ($\delta = 10\tunit$, $\Delta = 13\tunit$);
The gradient direction is $[1,1,0]$.
SpinD Simul 1:  $rtol = 10^{-3}$, $atol = 10^{-6}$, $Htetgen = 0.5$; 
SpinD Simul 2:  $rtol = 10^{-3}$, $atol = 10^{-6}$, $Htetgen = 1$;
SpinD Simul 3:  $rtol = 10^{-3}$, $atol = 10^{-6}$, $Htetgen = -1$;
Camino Simul 1:  $N = 1000$, $T = 200$;
Camino Simul 2:  $N = 4000$, $T = 800$;
}
\end{figure}

We computed the signal differences between the reference simulations and the 2 remaining SpinDoctor simulations as well as the two Camino signals: 
\be{}
E(b) = \left|\frac{S(b)}{S(0)} - \frac{S^{ref}(b)}{S^{ref}(0)}\right|\times 100.
\ee
In Figure \ref{fig:ecs400} we see $E(b)$ for the SpinDoctor simulation on the coarsest mesh  ($Htetgen = -1$) 
is less than $0.4\%$ for all b-values and for the SpinDoctor simulation on the mesh ($Htetgen = 1$) it is less than $0.2\%$.  
The Camino simulation with ($N=1000$, $T=200$) has a signal difference of 1.9\% for b-value up to $2000\bunit$, and 
the Camino simulation with ($N=4000$, $T=800$) has a signal difference of 0.7\% for b-value up to $2000\bunit$.  However,
for b-value $b=3000\bunit$ and greater, it seems the first Camino simulation is closer to the reference signal than the second Camino simulation.
It likely means that 4000 spins and 800 time steps are not enough to achieve signal convergence at higher b-values.
\dcom{In fact, they are below the recommended values for Monte-Carlo simulations \cite{Hall2009}, but we chose them
to keep the Camino simulations running within a reasonable amount of time.} 
On the other hand, the refinement of the FE mesh for the SpinDoctor achieves convergence 
for all b-values up to 10000$\bunit$.  There is a significant increase of the computational time of 
SpinDoctor as the diffusion-encoding amplitude is increased from 0.03 T/m to 0.37 T/m.  
At the finest mesh, the computational
time increased from 35 seconds to 200 seconds.  At the coarsest mesh, the computational
time increased from 20 seconds to 115 seconds.  This is due to the fact that at higher gradient
amplitudes, the magnetization is more oscillatory, so to achieve a fixed ODE solver tolerance,
smaller time steps are needed.

In Table \ref{table:ecs400} we show the total computational time to compute the 
dMRI signals at the 8 b-values for 2 SpinDoctor and 2 Camino simulations.  We also include 
the time for Camino to place the initial spins in the geometry described by the PLY file. 
We include in the Table the maximum signal differences for b-values up to $2000\bunit$ instead of all the b-values because 
Camino is not convergent for b-values greater than $3000\bunit$.
We see that at a similar level of signal difference ($0.4\%$ for SpinDoctor versus $0.7\%$ for Camino), the total computational time of SpinDoctor (438 seconds)
is more than 100 times faster than Camino (59147 seconds).
\begin{table}[!htb]
\begin{center}
\begin{tabular}{|c|c|c|c|c|}
\hline
\multirow{2}{*}{ {\it ECS400axons}} & \multicolumn{2}{c|}{SpinDoctor} & \multicolumn{2}{c|}{Camino}  \\ 
\cline{2-5}
 & Htet = -1  & Htet = 0.5 & $T=200$  & $T=800$\\
\hline
Degrees  & 53280 nodes &  70047 nodes & \multirow{2}{*}{1000 spins}  & \multirow{2}{*}{4000 spins}  \\
of freedom& 125798 elements &  177259 elements &   & \\
\hline Max signal difference & \multirow{2}{*}{0.4\%}   &  \multirow{2}{*}{Ref signal}  & \multirow{2}{*}{1.9\%} & \multirow{2}{*}{0.7\%} \\
($b\leq 2000\bunit$) &    &  & & \\
  \hline  Initialization time (sec) & \multicolumn{2}{c|}{}  & 69& 305\\
\hline
Solve time (sec), 8 bvalues & 438 & 667  & 3949  & 58842\\
\hline
Total time (sec)  & 438   &  667  & 4018 & 59147\\
\hline 
\end{tabular} 
\caption{The geometry is {\it ECS400axons}.  The total computational times  (in seconds) to simulate the dMRI signal at 8 b-values using SpinDoctor and Camino.  The initialization time is the time for Camino to place initial
spins inside the geometry described by the PLY file.  The b-values simulated are $b=\{0, 100, 500, 1000, 2000, 3000, 6000, 10000\}\bunit$.  
The maximum signal differences are given for b-values up to $2000\bunit$ because 
Camino is not convergent for b-values greater than $3000\bunit$.  The diffusion coefficient is $2 \times 10^{-3}\dunit$;
The diffusion-encoding sequence is PGSE ($\delta = 10\tunit$, $\Delta = 13\tunit$);
The gradient direction is $[1,1,0]$.
\label{table:ecs400}}
\end{center}
\end{table}

\subsection{Dendrite branch}

\label{sec:dendritebranch}
SpinDoctor was run with the following 2 sets of simulation parameters: 
\begin{enumerate}[label=SpinD Simul \ref{sec:dendritebranch}-\arabic*:, wide=0pt, font=\textbf]
\item $rtol = 10^{-3}$, $atol = 10^{-6}$;
\item $rtol = 10^{-2}$, $atol = 10^{-4}$;
\end{enumerate}
The finite elements mesh was generated by an external package and imported into SpinDoctor.  
The finite elements mesh size is ($n_{nodes} = 24651, n_{elem} = 91689$).  We do not refine
the FE mesh, rather, we vary the ODE solve tolerances in the SpinDoctor simulations.

Camino was run with the following 3 sets of simulation parameters: 
\begin{enumerate}[label=Camino Simul \ref{sec:dendritebranch}-\arabic*:, wide=0pt, font=\textbf]
\item $N = 1000$, $T = 200$;
\item $N = 2000$, $T = 400$;
\item $N = 4000$, $T = 800$;
\end{enumerate} 

The reference signal is \textbf{SpinD Simul \ref{sec:dendritebranch}-1}, 
the SpinDoctor signal with the higher ODE solve tolerances ($rtol = 10^{-3}$, $atol = 10^{-6}$).

\begin{figure}[!htb]
  \centering
\includegraphics[width=0.49\textwidth]{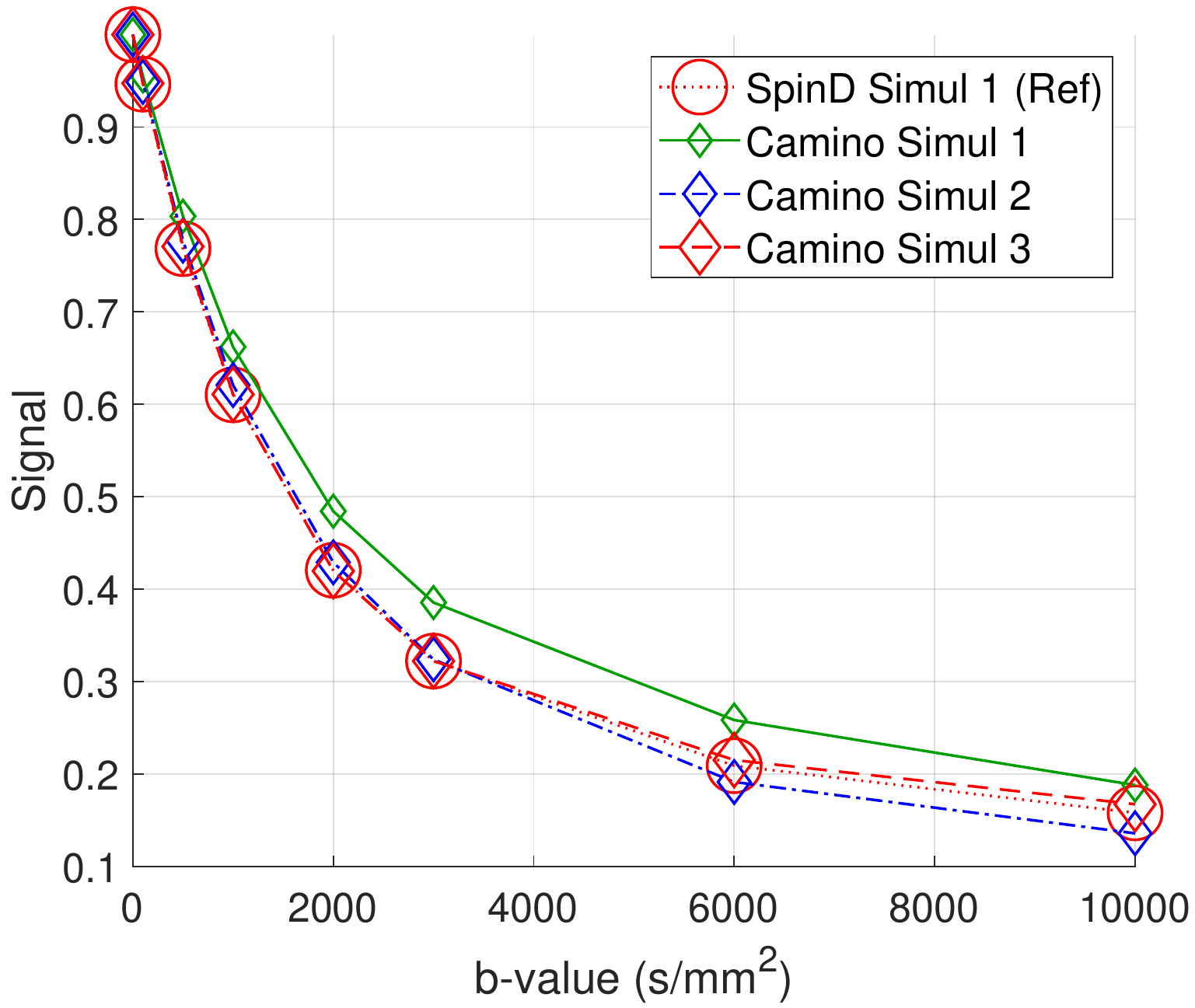} \\
\includegraphics[width=0.49\textwidth]{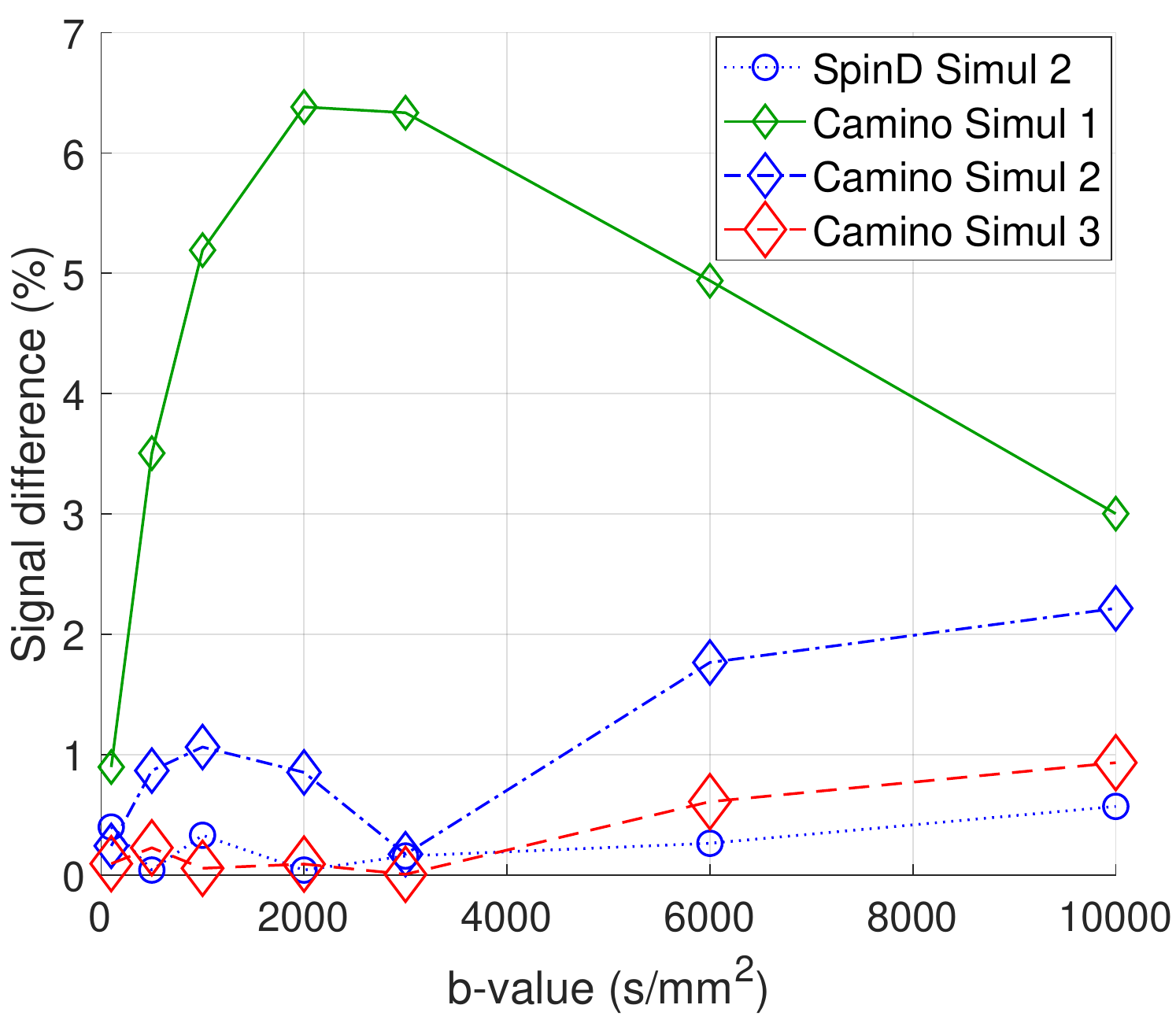} 
\includegraphics[width=0.49\textwidth]{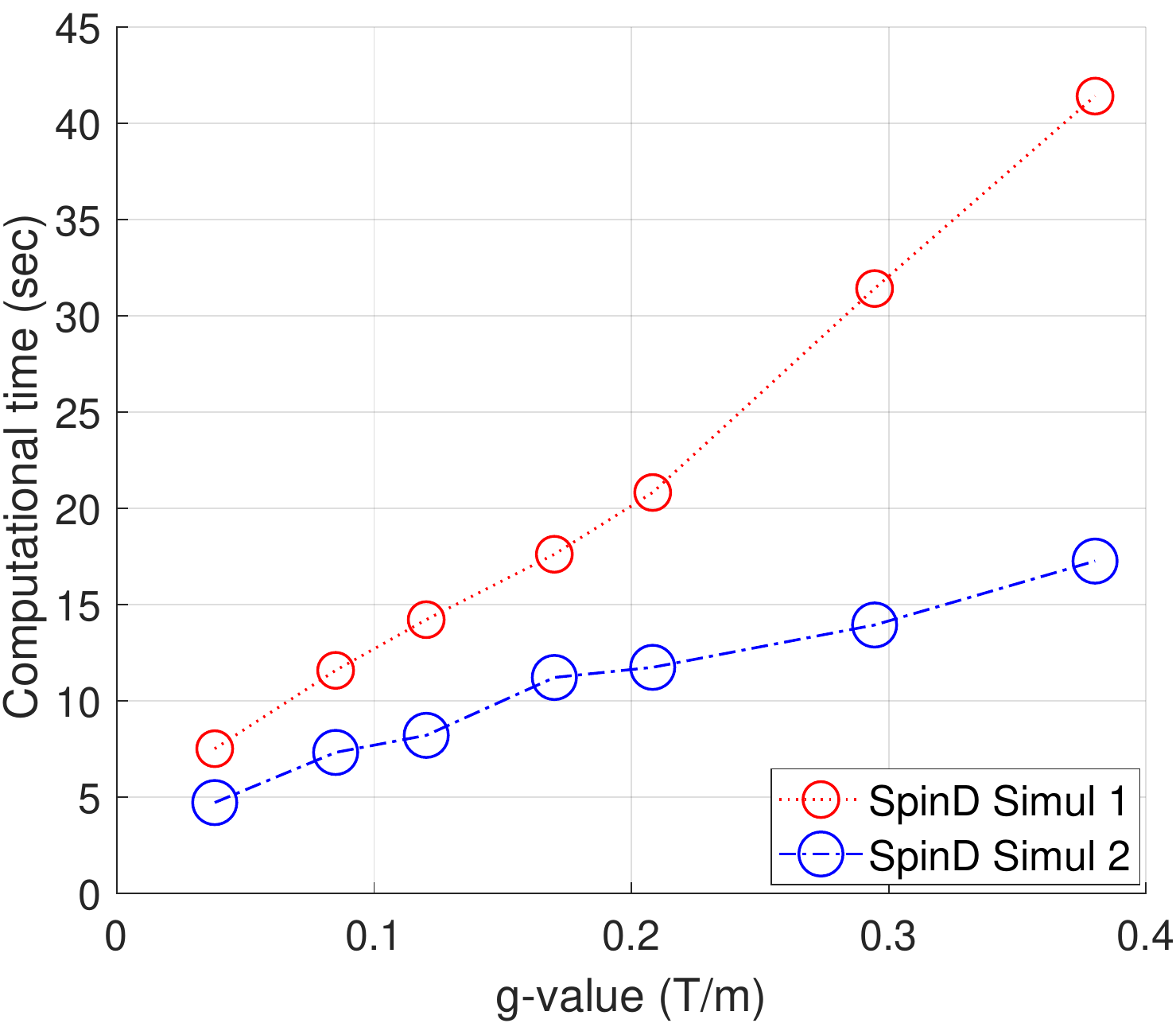} 
  \caption{ \label{fig:dendritebranch}The geometry is {\it DendriteBranch}.  Top: SpinD Simul 1 is the reference signal, compared to three Camino simulations.  Bottom left: the signal difference between the reference simulation and a SpinDoctor simulation and three Camino simulations.  Bottom right: the computational times of SpinDoctor simulations as a function of the gradient amplitudes.  The diffusion coefficient is $2 \times 10^{-3}\dunit$;
The diffusion-encoding sequence is PGSE ($\delta = 10\tunit$, $\Delta = 13\tunit$);
The gradient direction is $[1,1,0]$.
SpinD Simul 1: $rtol = 10^{-3}$, $atol = 10^{-6}$;
SpinD Simul 2: $rtol = 10^{-2}$, $atol = 10^{-4}$; 
Camino Simul 1: $N = 1000$, $T = 200$;
Camino Simul 2: $N = 2000$, $T = 400$;
Camino Simul 3: $N = 4000$, $T = 800$;
}
\end{figure}

In Figure \ref{fig:dendritebranch} we see the signal difference $E(b)$ 
for the SpinDoctor simulation with the bigger ODE solve tolerances ($rtol = 10^{-2}$, $atol = 10^{-4}$) 
is less than $0.6\%$ for all b-values.  
The Camino simulation with ($N=1000$, $T=200$) has a maximum signal difference of 6.4\%,  
the Camino simulation with ($N=4000$, $T=800$) has a maximum signal difference of 1.0\%.
As the gradient amplitude is increased from 0.03 T/m to 0.37 T/m, 
at the larger ODE solve tolerances, the computational
time increased from 5 seconds to 17 seconds.  At smaller ODE solve tolerances, the computational
time increased from 7 seconds to 42 seconds.  Again, this increase is due to the fact that at higher gradient
amplitudes, the magnetization is more oscillatory, so to achieve a fixed ODE solver tolerance,
smaller time steps are needed.
In Table \ref{table:dendritebranch} 
we see for the same level of accuracy ($0.6\%$ for SpinDoctor and and $1 \%$ for Camino), 
SpinDoctor (109 seconds) is 400 times faster than Camino (43918 seconds).

\begin{table}[!htb]
\begin{center}
\begin{tabular}{|c|c|c|c|c|c|}
\hline
$Dendrite$ & \multicolumn{2}{c|}{SpinDoctor} & \multicolumn{3}{c|}{Camino}  \\ 
\cline{2-6}
$Branch$ & $rtol = 10^{-2}$ & $rtol = 10^{-3}$ & $ T=200$ & $ T=400$ & $T=800$ \\
\hline
Degrees &\multicolumn{2}{c|}{24651 nodes}& \multirow{2}{*}{1000 spins} & \multirow{2}{*}{2000 spins} & \multirow{2}{*}{4000 spins}   \\
of freedom  &  \multicolumn{2}{c|}{91689 elements}  &  &  & \\
\hline
Max signal difference  & 0.6\% & Ref signal & 6.4\% & 2.2\% & 1.0\%\\\hline
Initialization time (sec) & \multicolumn{2}{c|}{} &5897 &11739 & 23702 \\
\hline
Solve time (sec), 8 bvalues & 109 & 207  &1336 &5138 &20216 \\
\hline
Total time (sec) & 109 & 207 & 7233 & 16877& 43918\\
\hline
\end{tabular} 
\caption{The geometry is  {\it DendriteBranch}. 
The total computational times  in seconds to simulate the dMRI signal at 8 b-values using SpinDoctor and Camino.  The initialization time is the time for Camino to place initial
spins inside the geometry described by the PLY file.  The b-values simulated are $b=\{0, 100, 500, 1000, 2000, 3000, 6000, 10000\}\bunit$.  The diffusion coefficient is $2 \times 10^{-3}\dunit$;
The diffusion-encoding sequence is PGSE ($\delta = 10\tunit$, $\Delta = 13\tunit$);
The gradient direction is $[1,1,0]$.
\label{table:dendritebranch}}
\end{center}
\end{table}

\subsection{Three dimensional ECS of 200 axons}
\label{sec:ecs200}
Due to computational time limitations, we only computed 4 b-values, $b = \{0, 100, 500, 1000\}\bunit$,
for the geometry {\it ECS200axons} (see Figure \ref{fig:FE_ecs200axons} for the finite element mesh).   

SpinDoctor was run with the following 2 sets of simulation parameters: 
\begin{enumerate}[label=SpinD Simul \ref{sec:ecs200}-\arabic*:, wide=0pt, font=\textbf]
\item  $rtol = 10^{-3}$, $atol = 10^{-6}$, $Htetgen = -1$ 
\item  $rtol = 10^{-3}$, $atol = 10^{-6}$, $Htetgen = 0.3$.
\end{enumerate}
 For this geometry, $Htetgen = -1$ gives finite elements mesh size ($n_{nodes} = 846298, n_{elem} = 2997386$).
$Htetgen = 0.3$ gives finite elements mesh size ($n_{nodes} = 1017263, n_{elem} = 3950572$).

The difference inthe  signals between the two simulations is less than $0.35\%$ (not plotted),
meaning the FE meshes are fine enough to produce accurate signals.  
In Table \ref{table:ecs200}, we see that using about 846K
nodes required 1.8 hours at $b=100\bunit$, 2.7 hours at $b = 500\bunit$, 3.3 hours at $b = 1000\bunit$.  
We did not use Camino for  {\it ECS200axons}
due to the excessive time required by Camino.

\begin{table}[!htb]
\begin{center}
\begin{tabular}{|c|c|c|}
\hline
\multirow{2}{*}{{\it ECS200 axons}} & \multicolumn{2}{c|}{SpinDoctor}  \\ 
\cline{2-3}
 & Htet = -1 & Htet = 0.3 \\
\hline
\multirow{2}{*}{Mesh} & 846298 nodes & 1017263 nodes \\
 & 2997386 elements&  3950572 elements \\
\hline 
Max signal difference & 0.35\%& Ref signal \\
\hline
Solve time (sec), $b=100, 500, 1000\bunit$ & (6611, 9620, 12107) & (16978, 23988, 32044) \\
\hline 
\end{tabular} 
\caption{The geometry is {\it ECS200axons}.  The computational times  in seconds
to simulate the dMRI signal at 3 b-values $b=\{100, 500, 1000\}\bunit$ 
using SpinDoctor.  The times are listed separately for each b-value.
The diffusion coefficient is $2 \times 10^{-3}\dunit$;
The diffusion-encoding sequence is PGSE ($\delta = 10\tunit$, $\Delta = 13\tunit$);
The gradient direction is $[1,1,0]$.\label{table:ecs200}}
\end{center}
\end{table}

\subsection{SpinDoctor computational time}

We collected the computational times of the SpinDoctor simulations for  
{\it ECS400axons}, {\it DendriteBranch}, and {\it ECS200axons},
that had the ODE solve tolerances ($rtol = 10^{-3}$, $atol = 10^{-6}$). 
In addition, for {\it ECS400axons} and {\it DendriteBranch}, we performed simulations for
another PGSE sequence ($\delta = 10\tunit$, $\Delta = 23\tunit$).

\begin{figure}[!htb]
  \centering
 \includegraphics[width=0.65\textwidth]{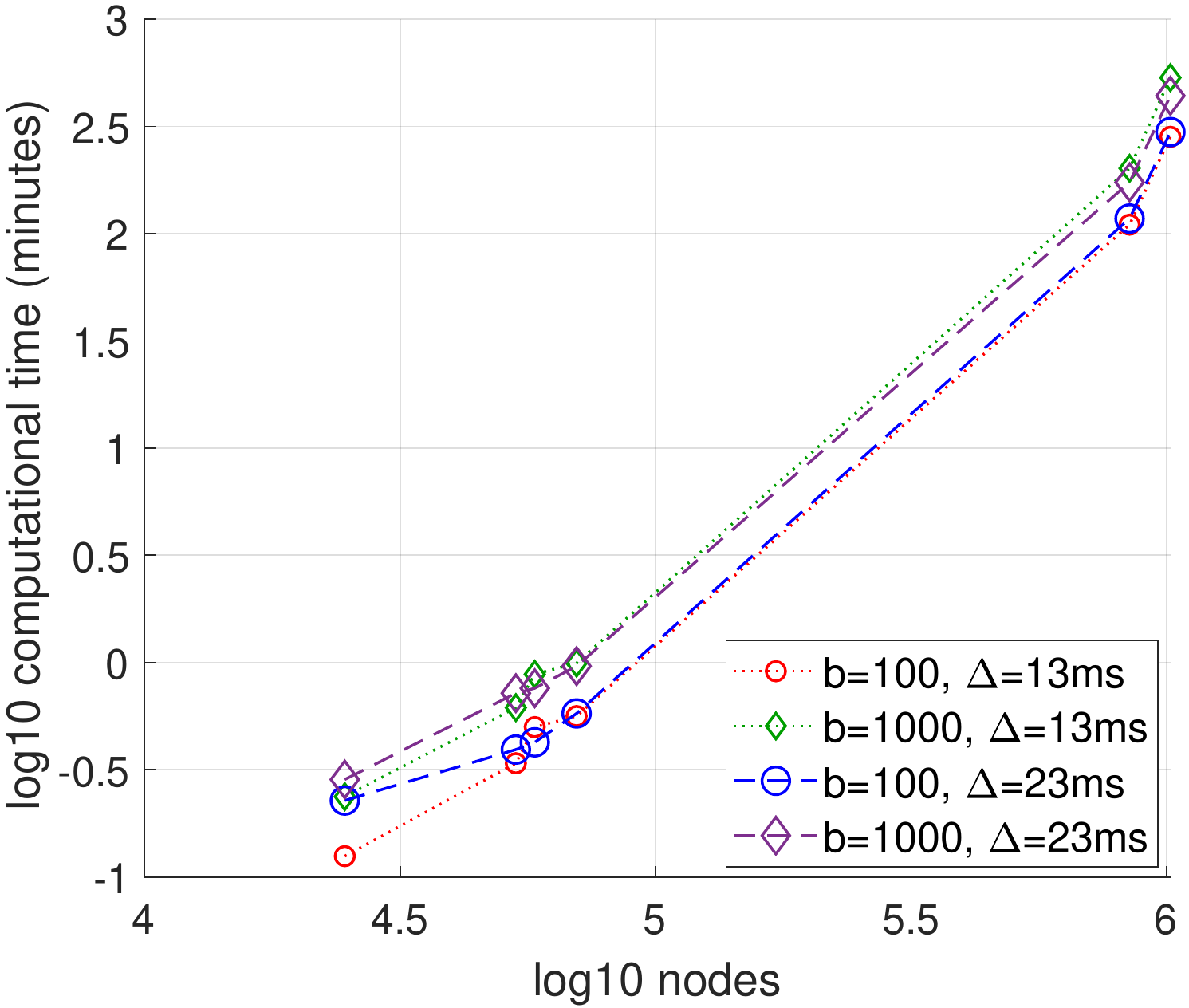} 
  \caption{Computational times of SpinDoctor to simulate one b-value (either $b=100\bunit$ or $b=1000\bunit$).  The x-axis gives log 10 of the number of finite elements nodes.  The data include 3 FE meshes of {\it ECS400axons}, 1 FE mesh of {\it DendriteBranch}, and 2 FE meshes of 
{\it ECS200axons}.  The y-axis gives the log 10 of the comptational time in minutes.   Below $y=0$ are computational times that are less than one minute.  The two sequences simulated are PGSE sequence ($\delta = 10\tunit$, $\Delta = 13\tunit$) and PGSE sequence ($\delta = 10\tunit$, $\Delta = 23\tunit$).  The diffusion coefficient is $2 \times 10^{-3}\dunit$;
The gradient direction is $[1,1,0]$.  
\label{fig:spindoctor_all_ctime}}
\end{figure}

Now we examine the computational time as a function of the finite element mesh size
for those simulations with ODE solve tolerances ($rtol = 10^{-3}$, $atol = 10^{-6}$). 
There are 3 FE meshes of {\it ECS400axons}, 1 FE mesh of {\it DendriteBranch}, and 2 FE meshes of 
{\it ECS200axons}.  In Figure \ref{fig:spindoctor_all_ctime} we plot the computational times to simulate the dMRI signal 
at two b-values ($b=100\bunit$ and $b=1000\bunit$) as a function of the 
number of FE nodes.  We see at fewer than 100K finite element nodes, 
the SpinDoctor simulation time is less than 1 minute per b-value.  
At 1 million FE nodes, the SpinDoctor simulation time is about 4.7 hours for $b=100\bunit$ and 8.9 hours for $b=1000\bunit$.

\section{SpinDoctor permeability and Monte-Carlo transmission probability}

\marginparnew{New section added in revised version}
\label{sec:perm}
Here we illustrate the link between the membrane permeability of Spindoctor and the transmission 
probability of crossing a membrane in the Camino simulation.
The geometry is the following:
\begin{itemize}
\item {\it Permeable Sphere} involves uniformly placed initial spins inside a sphere of radius $5\lunit$, subject to permeable interface condition on the surface of the sphere, 
with permeability coefficient $\kappa$.  No spins are initially placed outside of this sphere.
In the SpinDoctor simulation, this sphere is enclosed inside a sphere of diameter $30\lunit$, subject to impermeable boundary condition on the outermost interface.  
In the Camino simulation, this sphere is enclosed in a box of side length $30\lunit$, subject to periodic boundary conditions.  The inner sphere is far enough from the outer sphere in SpinDoctor and from the outer box in Camino so 
that there is no influence of the outer surface during the simulated diffusion times.
\end{itemize}
The dMRI experimental parameters are the following:
\begin{itemize}
\item the diffusion coefficient in all compartments is $2 \times 10^{-3}\dunit$;
\item the diffusion-encoding sequence is PGSE ($\delta = 10\tunit$, $\Delta = 13\tunit$);
\item 8 b-values: $b =\{ 0, 100, 500, 1000, 2000, 3000, 6000, 10000\}\bunit$;
\item 1 gradient direction: $[1,1,0]$.
\end{itemize}
SpinDoctor was run with the following 3 sets of simulation parameters: 
\begin{enumerate}[label=SpinD Simul \ref{sec:perm}-\arabic*: , wide=0pt, font=\textbf]
\item $rtol = 10^{-3}$, $atol = 10^{-6}$, $Htetgen = 0.5$;
\item $rtol = 10^{-3}$, $atol = 10^{-6}$, $Htetgen = 1$;
\item $rtol = 10^{-3}$, $atol = 10^{-6}$, $Htetgen = -1$;
\end{enumerate}
 For this geometry, $Htetgen = -1$ gives finite elements mesh size ($n_{nodes} = 46384, n_{elem} = 196920$). $Htetgen = 1$ gives finite elements mesh size ($n_{nodes} = 49618, n_{elem} = 218007$). $Htetgen = 0.5$ gives finite elements mesh size ($n_{nodes} = 52803, n_{elem} = 237613$).
 
Camino was run with the following 2 sets of simulation parameters: 
\begin{enumerate}[label=Camino Simul \ref{sec:perm}-\arabic*:, wide=0pt, font=\textbf]
\item $N = 4000$, $T = 800$;
\item $N = 8000$, $T = 3200$;
\end{enumerate}

The reference signal is \textbf{SpinD Simul \ref{sec:perm}-1}, the SpinDoctor signal on the finest FE mesh.

In \cite{FIEREMANS201839}, there is a discussion about the transmission probability of random walkers 
as they encounter a permeable membrane with permeability $\kappa$.  
The formula found in that paper is (for three dimensions)
\be{eq:pex}
P_{EX} = C_{dim}\frac{\kappa}{\sigma} \sqrt{2 \,dim \,\sigma \,\delta t } , \quad C_{dim} = \frac{2}{3}, \quad dim = 3,
\ee
$\sigma$ being the intrinsic diffusion coefficient, $\delta t$ is the time step.

\begin{figure}[!htb]
  \centering
\includegraphics[width=0.49\textwidth]{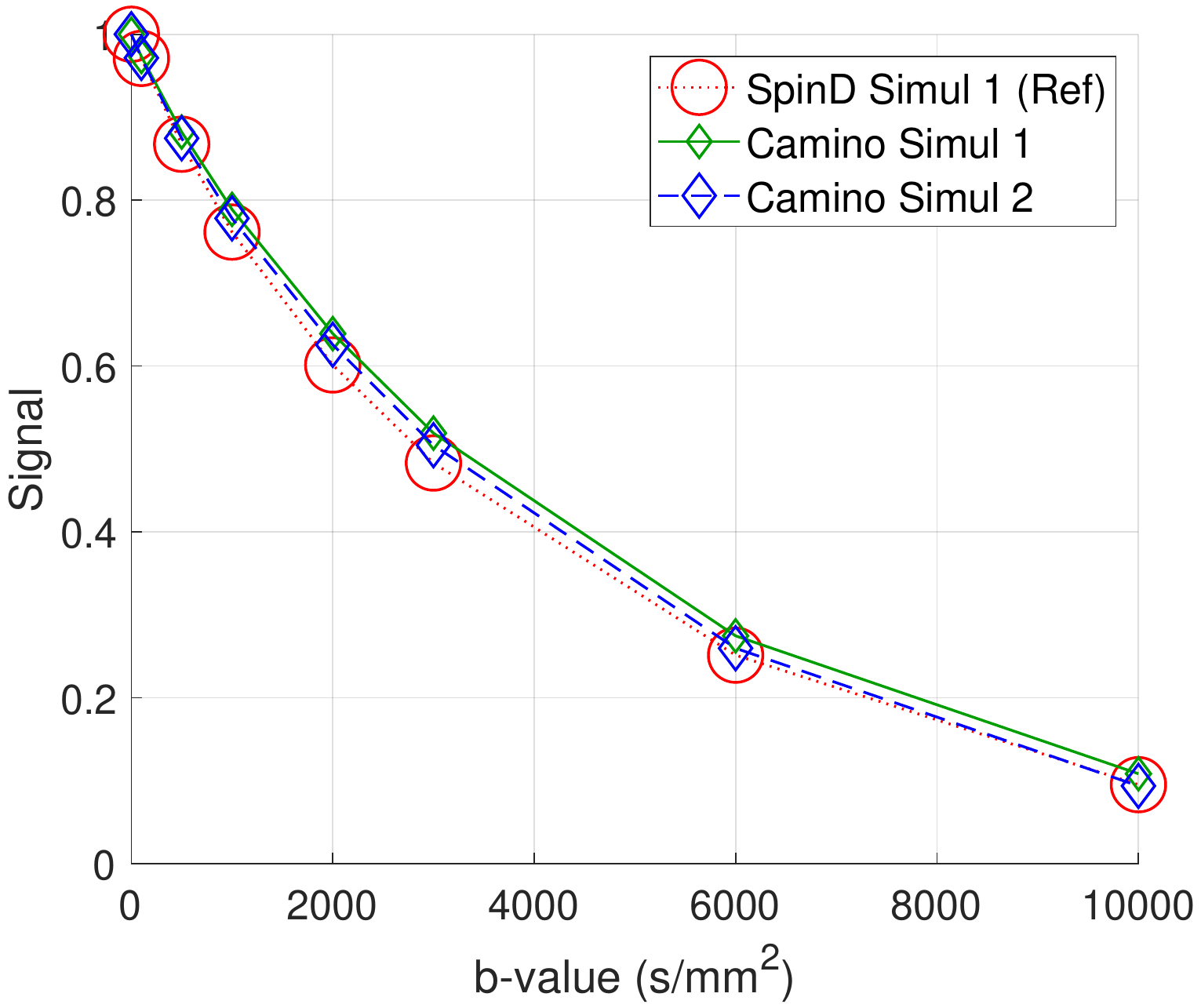} 
\includegraphics[width=0.49\textwidth]{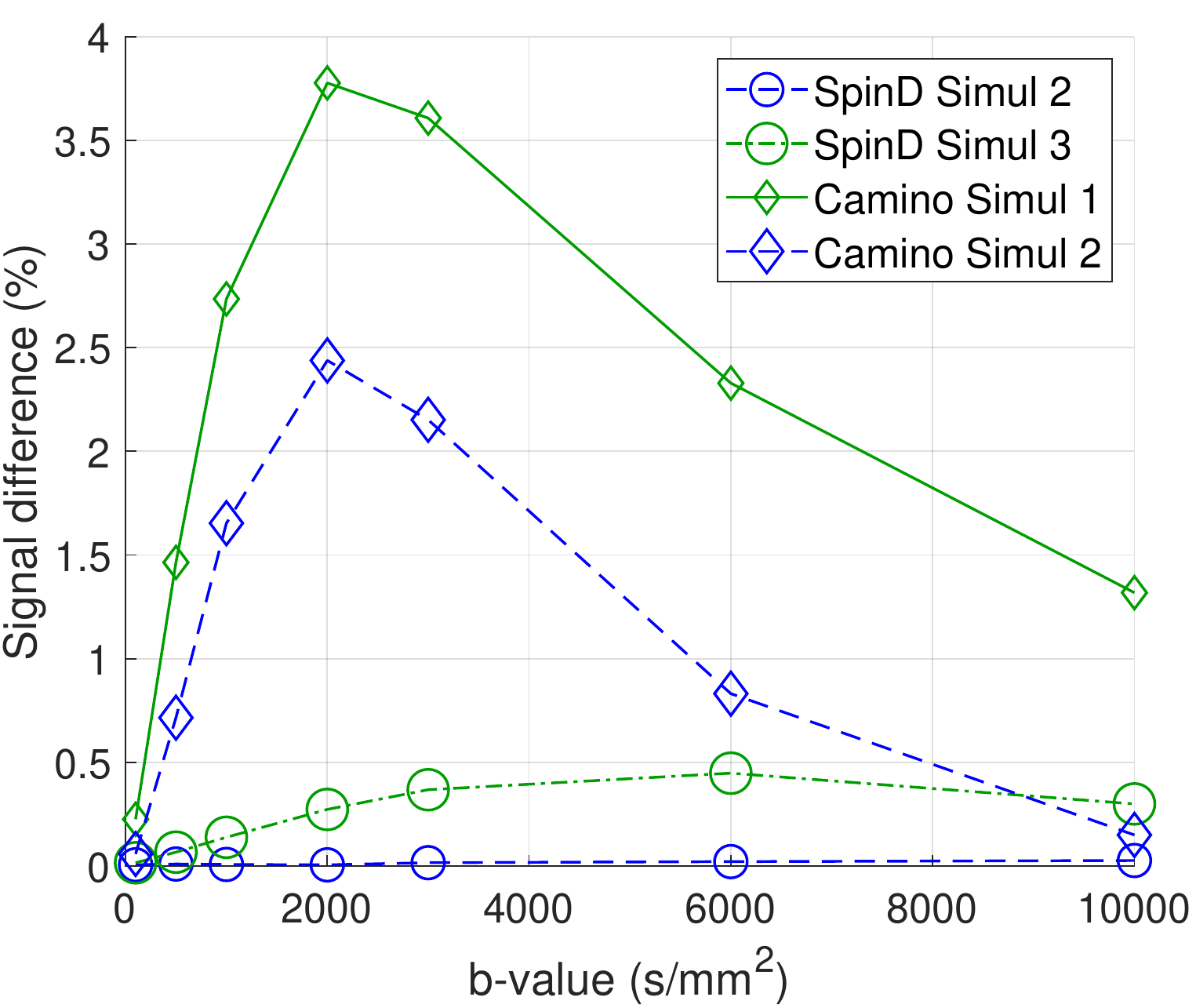} 
  \caption{The {\it Permeable Sphere} example involves uniformly placed initial spins inside a sphere of radius $5\lunit$, subject to permeable interface condition on the surface of the sphere, with permeability coefficient $\kappa= 10^{-5}\kunit$. 
Left: the SpinDoctor simulation on the finest mesh as the reference signal and two Camino signals.  Right: the signal difference between the reference signal and two SpinDoctor simulations and two Camino simulations.
SpinD Simul 1: $rtol = 10^{-3}$, $atol = 10^{-6}$, $Htetgen = 0.5$;
SpinD Simul 2:  $rtol = 10^{-3}$, $atol = 10^{-6}$, $Htetgen = 1$;
SpinD Simul 3:  $rtol = 10^{-3}$, $atol = 10^{-6}$, $Htetgen = -1$;
Camino Simul 1: $N = 4000$, $T = 800$;
Camino Simul 2: $N = 8000$, $T = 3200$; 
The diffusion coefficient in all compartments is $2 \times 10^{-3}\dunit$;
The diffusion-encoding sequence is PGSE ($\delta = 10\tunit$, $\Delta = 13\tunit$);
The gradient direction is $[1,1,0]$.
 \label{fig:camino_perm}}
\end{figure}

In Figure \ref{fig:camino_perm} we show the three SpinDoctor simulations at $\kappa = 10^{-5}\kunit$ and the 
two Camino simulations using the above formula for $P_{EX}$.  We considered the SpinDoctor signal computed 
on the finest FE mesh as the reference signal and we computed the signal differences between the reference signal and the other two SpinDoctor signals and the Camino signals.  We see that the Camino signals approach the reference signal as the number of spins and times steps in Camino are increased, the maximum difference decreasing from 3.8\% to 2.4\%.  The SpinDoctor signals have signal differences of 
less than 0.5\% and 0.1\%, respectively.

\section{Extensions of SpinDoctor} 
\marginparnew{New section added in revised version}
Here we mention two extensions in the functionalities of SpinDoctor that are planned for a future release.

\subsection{Non-standard diffusion-encoding sequences}

Given the interest in nonstandard diffusion sequences beyond PGSE and OGSE, such as double diffusion encoding (see \cite{Shemesh2016, Dhital2019, Novikov2019, Henriques2019}) and multidimensional diffusion encoding (see \cite{Topgaard2017}), it is natural that SpinDoctor should easily support arbitrary diffusion-encoding sequences.  Besides the PGSE and 
the sine and cosine OGSE sequences that are provided in the SpinDoctor package, 
new sequences can be straightforwardly implemented by changing three files in the SpinDoctor package
\begin{itemize}
\item \texttt{SRC/DMRI/seqprofile.m} defines $f(t)$
\item \texttt{SRC/DMRI/seqintprofile.m} defines the integral $F(t) = \int_0^t f(s) ds$
\item \texttt{SRC/DMRI/seqbvaluenoq.m} defines the associated $b-$value.
\end{itemize}

\begin{figure}[!htb]
  \centering  
\includegraphics[width=0.65\textwidth]{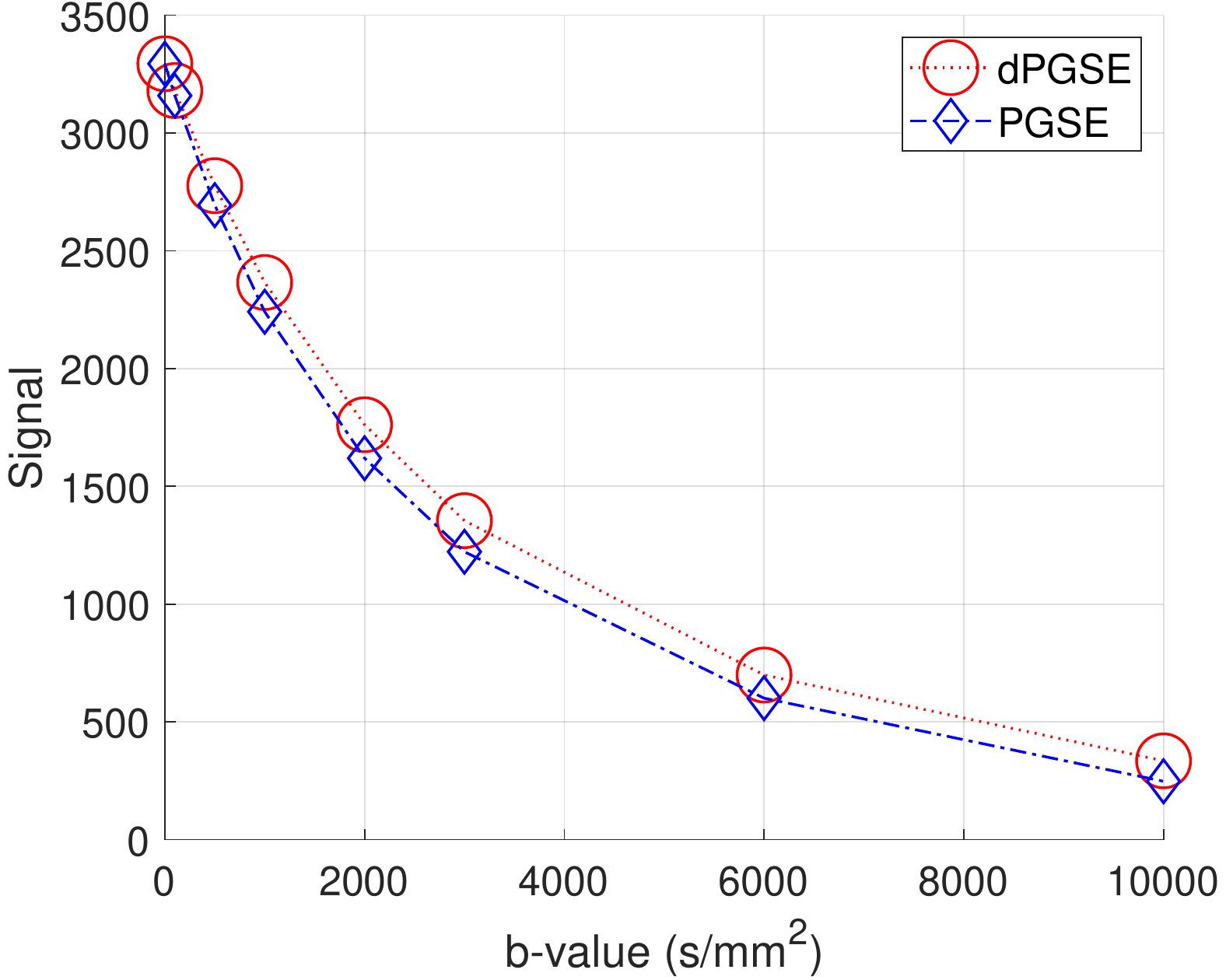} 
  \caption{DMRI signals of the PGSE and the double PGSE diffusion-encoding sequences.  
The geometry is made of cylindrical cells, the myelin layer, and the ECS (see Figure \ref{fig:myelin}).
The diffusion coefficient in all compartments is $2\times10^{-3}\dunit$ and the compartments do not experience 
spin exchange, with all permeability coefficients set to zero.  The diffusion-encoding sequeces are PGSE ($\delta = 10\tunit, \Delta = 13\tunit$) and dPGSE sequences ($\delta = 10\tunit, \Delta = 13\tunit, \tau=\delta+\Delta$), the diffusion-encoding direction is $\bug = [1,1,1]$.  
\label{fig:dpgse}}
\end{figure}

In the example below, we simulate the double-PGSE (Eq. \ref{eq:dpgse}) sequence:
\begin{equation}\label{eq:dpgse}
f(t) =
\begin{cases}
1, \quad &0 \leq t \leq \delta, \\
-1,
\quad & \Delta < t \leq \Delta+\delta,\\
1, \quad &\tau \leq t \leq \delta + \tau, \\
-1,
\quad & \Delta + \tau < t \leq \Delta+\delta + \tau,\\
0, \quad & \text{otherwise.}
\end{cases}
\end{equation}
Here $\delta$ is the duration of the diffusion-encoding gradient pulse, 
$\Delta$ is the time delay between the start of the two pulses, and $\tau$ is the distance between the 
two pairs of pulses ($\tau\geq \delta+\Delta$).
The geometry is made of cylindrical cells, the myelin layer, and the ECS (see Figure \ref{fig:myelin}).
In Figure \ref{fig:dpgse} we show the dMRI signals for the PGSE ($\delta = 10\tunit, \Delta = 13\tunit$) and dPGSE sequences ($\delta = 10\tunit, \Delta = 13\tunit, \tau=\delta+\Delta$), the diffusion-encoding direction is 
$\bug = [1,1,1]$.  

\subsection{$T_2$ relaxation}
When $T_2-$relaxation is considered, the Bloch-Torrey PDE (Eq. \ref{eq:btpde}) takes the following form
\begin{alignat}{4}
\label{eq:btpde_withT2}
&\frac{\partial}{\partial t}{M^{in}_i(\bx,t)} &&= -I\gamma f(t) \bg \cdot \bx \,M^{in}_i(\bx,t) -\frac{M^{in}}{T^{in}_2}
+ \nabla \cdot (\sigma^{in} \nabla M^{in}_i(\bx,t)), &&\bx \in \Omega^{in}_i,&\\
&\frac{\partial}{\partial t}{M^{out}_i(\bx,t)} &&= -I\gamma f(t) \bg \cdot \bx \,M^{out}_i(\bx,t) -\frac{M^{out}}{T^{out}_2}
+ \nabla \cdot (\sigma^{out} \nabla M^{out}_i(\bx,t)),  \;&&\bx \in \Omega^{out}_i,&\\
&\frac{\partial}{\partial t}{M^e(\bx,t)} &&= -I\gamma f(t) \bg \cdot \bx \,M^e(\bx,t)
-\frac{M^{e}}{T^{e}_2} 
+ \nabla \cdot (\sigma^{e} \nabla M^e(\bx,t)), & &\bx \in \Omega^{e},&
\end{alignat}

We plan to incorporate $T_2$ relaxation effects in the next official release of SpinDoctor.  In the meantime, this additional functionality 
can be found in a development branch of SpinDoctor available on GitHub.  
The source code in this development branch allows 
 the ability to add relaxation, with different relaxivities in
 the different compartments \cite{Veraart2018,Lampinen2019}.

\begin{figure}[!htb]
  \centering  
\includegraphics[width=0.65\textwidth]{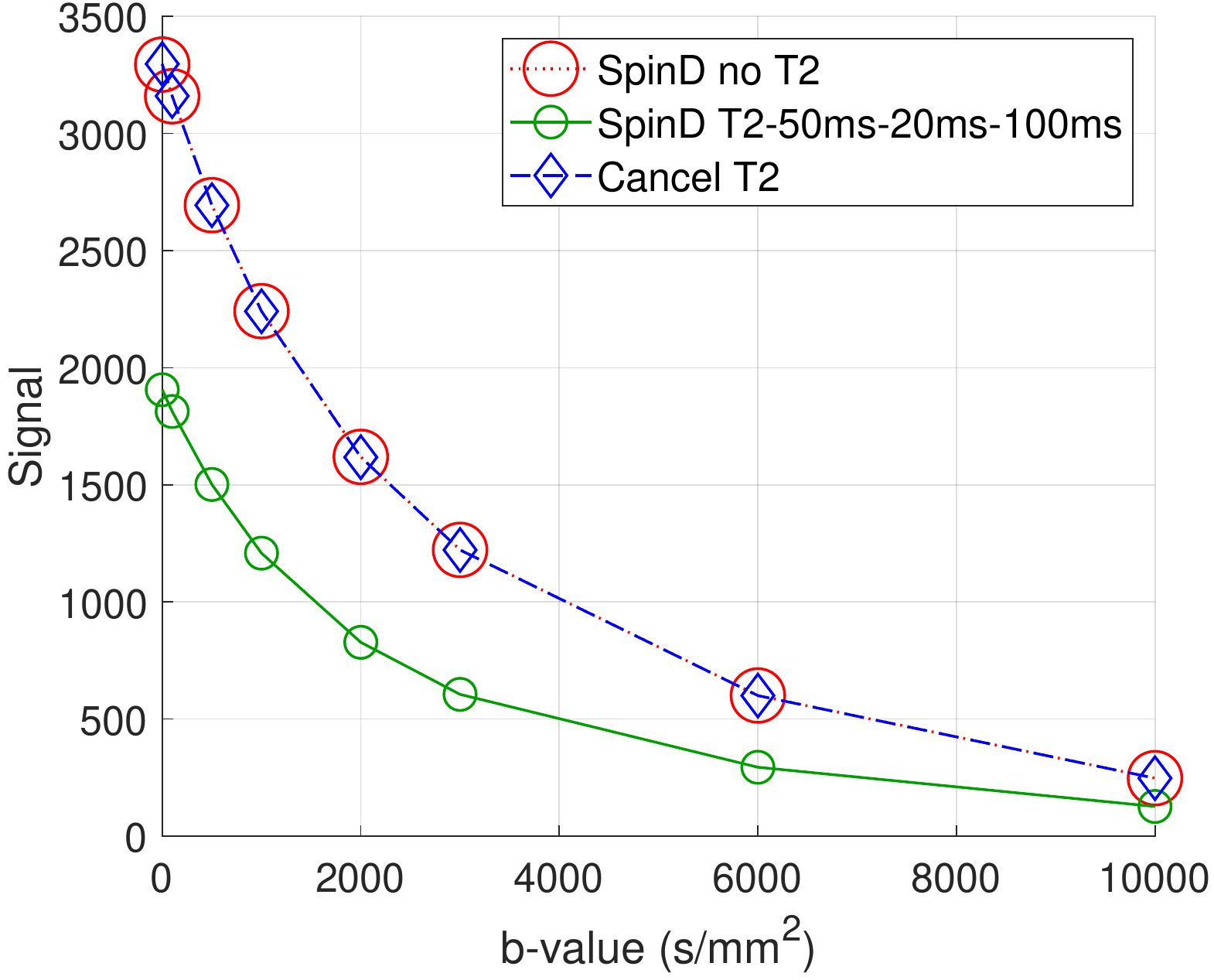} 
  \caption{DMRI signal including $T_2=[50\tunit, 20\tunit, 100\tunit]$ relaxation is lower than the signal without
relaxation effects ("no T2").  The $T_2$ effects are completely canceled out using Eq. \ref{eq:t2_cancellation}
so that the curve "cancel T2" coincides with the no relaxation signal.
The geometry is made of cylindrical cells, the myelin layer, and the ECS (see Figure \ref{fig:myelin}).
The diffusion coefficient in all compartments is $2\times10^{-3}\dunit$ and the compartments do not experience 
spin exchange, with all permeability coefficients set to zero.  The diffusion-encoding sequece is PGSE ($\delta = 10\tunit, \Delta = 13\tunit$), the diffusion-encoding direction is $\bug = [1,1,1]$.  
\label{fig:t2}}
\end{figure}

$T_2$ relaxation is incorporated using the format $T_2=[T_2^{in}, T_2^{out}, T_2^{e}]$ where $T_2^{in}, T_2^{out}, T_2^{e}$ are the $T_2$ values for the three compartments, respectively. 
To verify the correctness of our implementation, we check the following. 
Let $S_\text{No-T2}(b)$ be the signal without $T_2$ effects.  If there is no exchange between compartments, then the $T_2$ effects can be cancelled from the signals in the three compartments that include $T_2$ effects in the 
following way:
\begin{equation}\label{eq:t2_cancellation}
\begin{aligned}
{S_{cancel}}(b) &=\frac{S^{in}(b)}{ e^{-\frac{TE}{T_2^{in}}}} + \frac{S^{out}(b) }{ e^{-\frac{TE}{T_2^{out}}}}
+ \frac{S^{e}(b)}{e^{-\frac{TE}{T_2^{e}}}}.
\end{aligned}
\end{equation}
In Fig. \ref{fig:t2}, we compare $S_\text{No-T2}(b)$ with $S(b)$ where $T_2=[50\tunit, 20\tunit, 100\tunit]$, for the PGSE sequences ($\delta = 10\tunit, \Delta = 13\tunit$) and  $\bug = [1,1,1]$.  We also compute 
${S_{cancel}}(b)$, using Eq. \ref{eq:t2_cancellation}. 
The geometry (see Figure \ref{fig:myelin}) is made of cylindrical cells, the myelin layer, and the ECS.
The $T_2$ effects on the signal $S(b)$ are clearly seen.  The $T_2$ effects are completely canceled out using Eq. (\ref{eq:t2_cancellation}).

\section{Discussion}

Built upon MATLAB, SpinDoctor is a software package that seeks to reduce the work required to 
perform numerical simulations for dMRI for prototyping purposes.
There have been software packages for dMRI simulation that implements the random walkers approach.
A \dcom{detailed} comparison of the Monte-Carlo/random walkers approach with the FEM approach is beyond the scope of this paper.  SpinDoctor offers an alternative, solving the same physics problem using PDEs. 

After surveying other works on dMRI simulations, we saw a need to have a simulation toolbox that provides a way to easily define geometrical configurations.  
In SpinDoctor we have tried to offer useful configurations, without being overly general.  Allowing too much generality in the geometrical 
configurations would have made code robustness very difficult to achieve due to the difficulties
related to problems in computational geometry (high quality surface triangulation, robust FE mesh generation).  
The geometrical configuration routines provided by SpinDoctor are a helpful front end, 
to enable dMRI researchers to get started quickly to perform numerical simulations.
Those users who already have a high quality surface triangulation can use the other parts of SpinDoctor without
passing through this front end.

The bulk of SpinDoctor is the numerical solutions of two PDEs.  When one is only interested in the ADC, then computing the 
\soutnew{H-$ADC$}{HADC} model is the good option.   When one is interested in higher order behavior
in the dMRI signal, then the BTPDE model is a good option for accessing high b-value behavior.

Because time stepping methods for semi-discretized linear systems arising from finite element discretization is a well-studied subject in the mathematical literature, the ODE solvers implemented in MATLAB already optimize for such linear systems.  For example, 
the mass matrix is passed into the ODE solver as an optional parameter so as to avoid explicit matrix inversion.  In addition, 
the ODE solution is guaranteed to stay within a user-requested residual tolerance.  We believe this type of optimization and 
error control is clearly advantageous over simulation codes that do not have it.

To mimic the phenomenon where the water molecules can enter and exit the computational domain, the pseudo-periodic boundary conditions were implemented in \cite{Xu2007, Li2014, Nguyen2014}.  At this stage, we have chosen not to implement this in SpinDoctor, instead,
spins are not allowed to leave the computational domain.  Implementing pseudo-periodic boundary 
conditions would make the code more complicated, and it remains to be seen if it is a desired feature
among potential users.  If it is, then it could be part of a future development.

The twising and bending of the canonical configuration is something unique to SpinDoctor.  It removes
many computational geometry difficulties by meshing first the canonical configuration before deforming
the FE mesh via an analytical coordinate transformation. This is a way to simulate fibers that are not 
parallel, that bend, for example.  For fibers that disperse, perhaps more complicated analytical 
coordinate transformations can be performed on the canonical configuration to mimic that situation.  
This is a possible future direction to explore.  

SpinDoctor depends on MATLAB for the ODE solve routines as well as for the computational
geometry routines to produce the tight wrap ECS.  To implement SpinDoctor outside of MATLAB would require replacing 
these two sets of MATLAB routines.  Other routines of SpinDoctor can be easily implemented in another programming language.

SpinDoctor can be downloaded at \url{https://github.com/jingrebeccali/SpinDoctor}.

\soutnew{}{In summary, 
we have validated SpinDoctor simulations using reference signals from the Matrix Formalism method,
in particular in the case of permeable membranes.  We then compared SpinDoctor with the Monte-Carlo 
simulations produced by the publicly available software package Camino Diffusion MRI Toolkit \cite{Hall2009}.
We showed that the membrane permeability of SpinDoctor
can be straightforwardly linked to the transmission probability in Monte-Carlo simulations.
For numerous examples, it was seen that the SpinDoctor and the Camino simulations can be made 
close to each other if one increases
the degrees of freedom (the finite element mesh size for SpinDoctor and the number of spins for Camino) and increase
the accuracy of the time stepping (by tightening the ODE solve tolerances in SpinDoctor and by
increasing the number of time steps in Camino).  }

\soutnew{}{
At high gradient amplitudes, the ocsillatory nature of the magnetization requires 
the use of smaller time steps to maintain accuracy.  For this reason, the computational time to simulate the dMRI 
signal at high gradient amplitudes must be longer than at low gradient amplitudes. 
This adaptivity in the time stepping as a function of gradient amplitude is done automatically in SpinDoctor. }

\soutnew{}{
We have computed the dMRI signals on several complicated geometries on a stand-alone computer.
For these examples, we have shown that SpinDoctor can be more than 100 times faster than 
Camino.}  \dcom{Of course, in simple configurations such as straight, parallel cylinders, it is much more efficient to use an analytical representation of the diffusion environment rather than a triangulated mesh in Camino.  In addition, some recent 
implementations of random walk simulations \cite{Ginsburger2019, Rensonnet2019} should be faster than Camino.}

\soutnew{}{With a finite element mesh of 100K nodes, SpinDoctor takes less than one minute per b-value.  
At 1 million finite element nodes, limited computer memory resulted in a computational time 4.7 hours for $b=100\bunit$ and 8.9 hours for $b=1000\bunit$.
This issue will be taken into account in the future with high
performance computing techniques in MATLAB and on other platforms. 
One of our recent works \cite{Nguyen2018} is promising for this purpose.}

\soutnew{}{We also illustrated several extensions 
of SpinDoctor functionalities, including the incorporation of $T_2$ relaxation, the simulation of
non-standard diffusion-encoding sequences.  We note the dendrite branch example illustrates SpinDoctor's ability to import 
and use externally generated meshes provided by the user.}  
This capability will be 
very useful given the most recent developments in simulating ultra-realistic virtual tissues \cite{Palombo2019,Ginsburger2019}.

\section{Conclusion}\label{Conclusion}
This paper describes a publicly available MATLAB toolbox called SpinDoctor that can be used 
to solve the BTPDE to obtain
the dMRI signal and to solve the \dcom{diffusion equation} of the 
\soutnew{H-$ADC$}{HADC} model to obtain the ADC.
SpinDoctor is a software package that seeks to reduce the work required to 
perform numerical simulations for dMRI for prototyping purposes.

SpinDoctor provides built-in options of including spherical cells with a nucleus,
cylindrical cells with a myelin layer, an extra-cellular space enclosed either 
in a box or in a tight wrapping around the cells. 
The deformation of canonical cells by bending and twisting is implemented via an analytical coordinate 
transformation of the FE mesh.
Permeable membranes for the BTPDE is implemented using double nodes on the compartment interfaces.
Built-in diffusion-encoding pulse sequences include the Pulsed Gradient Spin Echo and the Ocsillating Gradient Spin Echo.  Error control in the time stepping is done using built-in MATLAB ODE solver routines.

User feedback to improve SpinDoctor is welcomed.

\section*{Acknowledgment}
The authors gratefully acknowledge the {\it French-Vietnam Master in Applied Mathematics} program 
whose students (co-authors on this paper, {Van-Dang Nguyen}, {Try Nguyen Tran},
{Bang Cong Trang},
{Khieu Van Nguyen},
{Vu Duc Thach Son},
{Hoang An Tran},
{Hoang Trong An Tran},
{Thi Minh Phuong Nguyen}) 
have contributed to the SpinDoctor project during their internships in France in the past 
several years, as well
as the {\it Vice-Presidency for Marketing and International Relations} at Ecole Polytechnique for financially supporting
a part of the students' stay.  Jan Valdman was supported by the Czech Science Foundation (GACR), through the grant GA17-04301S.  Van-Dang Nguyen was supported by the Swedish Energy Agency, Sweden with the project ID P40435-1 and MSO4SC with the grant number 731063.






\bibliographystyle{elsarticle-num}
\bibliography{myref}

\begin{thebibliography}{10}
\expandafter\ifx\csname url\endcsname\relax
  \def\url#1{\texttt{#1}}\fi
\expandafter\ifx\csname urlprefix\endcsname\relax\def\urlprefix{URL }\fi
\expandafter\ifx\csname href\endcsname\relax
  \def\href#1#2{#2} \def\path#1{#1}\fi

\bibitem{Hahn1950}
E.~L. Hahn, Spin echoes, Phys. Rev. 80 (1950) 580--594.

\bibitem{Stejskal1965}
E.~O. Stejskal, J.~E. Tanner, Spin diffusion measurements: Spin echoes in the
  presence of a time-dependent field gradient, The Journal of Chemical Physics
  42~(1) (1965) 288--292.
\newblock \href {https://doi.org/10.1063/1.1695690}
  {\path{doi:10.1063/1.1695690}}.

\bibitem{Bihan1986}
D.~L. Bihan, E.~Breton, D.~Lallemand, P.~Grenier, E.~Cabanis, M.~Laval-Jeantet,
  {MR} imaging of intravoxel incoherent motions: application to diffusion and
  perfusion in neurologic disorders., Radiology 161~(2) (1986) 401--407, pMID:
  3763909.

\bibitem{Does2003}
M.~D. Does, E.~C. Parsons, J.~C. Gore, Oscillating gradient measurements of
  water diffusion in normal and globally ischemic rat brain, Magn. Reson. Med.
  49~(2) (2003) 206--215.
\newblock \href {https://doi.org/10.1002/mrm.10385}
  {\path{doi:10.1002/mrm.10385}}.

\bibitem{Jensen2005}
J.~H. Jensen, J.~A. Helpern, A.~Ramani, H.~Lu, K.~Kaczynski, Diffusional
  kurtosis imaging: The quantification of non-{G}aussian water diffusion by
  means of magnetic resonance imaging, Magnetic Resonance in Medicine 53~(6)
  (2005) 1432--1440.
\newblock \href {https://doi.org/10.1002/mrm.20508}
  {\path{doi:10.1002/mrm.20508}}.

\bibitem{Tuch2002}
D.~S. Tuch, T.~G. Reese, M.~R. Wiegell, N.~Makris, J.~W. Belliveau, V.~J.
  Wedeen, \href{https://onlinelibrary.wiley.com/doi/abs/10.1002/mrm.10268}{High
  angular resolution diffusion imaging reveals intravoxel white matter fiber
  heterogeneity}, Magnetic Resonance in Medicine 48~(4) (2002) 577--582.
\newblock \href
  {http://arxiv.org/abs/https://onlinelibrary.wiley.com/doi/pdf/10.1002/mrm.10268}
  {\path{arXiv:https://onlinelibrary.wiley.com/doi/pdf/10.1002/mrm.10268}},
  \href {https://doi.org/10.1002/mrm.10268} {\path{doi:10.1002/mrm.10268}}.
\newline\urlprefix\url{https://onlinelibrary.wiley.com/doi/abs/10.1002/mrm.10268}

\bibitem{Assaf2008}
Y.~Assaf, T.~Blumenfeld-Katzir, Y.~Yovel, P.~J. Basser,
  \href{https://onlinelibrary.wiley.com/doi/abs/10.1002/mrm.21577}{Axcaliber:
  {A} method for measuring axon diameter distribution from diffusion {MRI}},
  Magnetic Resonance in Medicine 59~(6) (2008) 1347--1354.
\newblock \href
  {http://arxiv.org/abs/https://onlinelibrary.wiley.com/doi/pdf/10.1002/mrm.21577}
  {\path{arXiv:https://onlinelibrary.wiley.com/doi/pdf/10.1002/mrm.21577}},
  \href {https://doi.org/10.1002/mrm.21577} {\path{doi:10.1002/mrm.21577}}.
\newline\urlprefix\url{https://onlinelibrary.wiley.com/doi/abs/10.1002/mrm.21577}

\bibitem{Alexander2010}
D.~C. Alexander, P.~L. Hubbard, M.~G. Hall, E.~A. Moore, M.~Ptito, G.~J.
  Parker, T.~B. Dyrby,
  \href{http://www.sciencedirect.com/science/article/pii/S1053811910007755}{Orientationally
  invariant indices of axon diameter and density from diffusion {MRI}},
  NeuroImage 52~(4) (2010) 1374--1389.
\newline\urlprefix\url{http://www.sciencedirect.com/science/article/pii/S1053811910007755}

\bibitem{Zhang2011}
H.~Zhang, P.~L. Hubbard, G.~J. Parker, D.~C. Alexander,
  \href{http://www.sciencedirect.com/science/article/pii/S1053811911001376}{Axon
  diameter mapping in the presence of orientation dispersion with diffusion
  {MRI}}, NeuroImage 56~(3) (2011) 1301--1315.
\newline\urlprefix\url{http://www.sciencedirect.com/science/article/pii/S1053811911001376}

\bibitem{Zhang2012}
H.~Zhang, T.~Schneider, C.~A. Wheeler-Kingshott, D.~C. Alexander,
  \href{http://www.sciencedirect.com/science/article/pii/S1053811912003539}{{NODDI}:
  {Practical} in vivo neurite orientation dispersion and density imaging of the
  human brain}, NeuroImage 61~(4) (2012) 1000--1016.
\newline\urlprefix\url{http://www.sciencedirect.com/science/article/pii/S1053811912003539}

\bibitem{Burcaw2015}
L.~M. Burcaw, E.~Fieremans, D.~S. Novikov, Mesoscopic structure of neuronal
  tracts from time-dependent diffusion, NeuroImage 114 (2015) 18 -- 37.
\newblock \href {https://doi.org/10.1016/j.neuroimage.2015.03.061}
  {\path{doi:10.1016/j.neuroimage.2015.03.061}}.

\bibitem{Palombo2017a}
M.~Palombo, C.~Ligneul, J.~Valette, {Modeling diffusion of intracellular
  metabolites in the mouse brain up to very high diffusion-weighting: Diffusion
  in long fibers (almost) accounts for non-monoexponential attenuation},
  Magnetic Resonance in Medicine 77~(1) (2017) 343--350.
\newblock \href {https://doi.org/10.1002/mrm.26548}
  {\path{doi:10.1002/mrm.26548}}.

\bibitem{Palombo2016}
M.~Palombo, C.~Ligneul, C.~Najac, J.~Le~Douce, J.~Flament, C.~Escartin,
  P.~Hantraye, E.~Brouillet, G.~Bonvento, J.~Valette, {New paradigm to assess
  brain cell morphology by diffusion-weighted MR spectroscopy in vivo},
  Proceedings of the National Academy of Sciences 113~(24) (2016) 6671--6676.
\newblock \href
  {http://arxiv.org/abs/http://www.pnas.org/content/113/24/6671.full.pdf}
  {\path{arXiv:http://www.pnas.org/content/113/24/6671.full.pdf}}, \href
  {https://doi.org/10.1073/pnas.1504327113}
  {\path{doi:10.1073/pnas.1504327113}}.

\bibitem{Ning2017}
L.~Ning, E.~Özarslan, C.-F. Westin, Y.~Rathi, Precise inference and
  characterization of structural organization (picaso) of tissue from molecular
  diffusion, NeuroImage 146 (2017) 452 -- 473.
\newblock \href {https://doi.org/10.1016/j.neuroimage.2016.09.057}
  {\path{doi:10.1016/j.neuroimage.2016.09.057}}.

\bibitem{McHugh2015}
D.~J. McHugh, F.~Zhou, P.~L. Hubbard~Cristinacce, J.~H. Naish, G.~J.~M. Parker,
  Ground truth for diffusion {MRI} in cancer: A model-based investigation of a
  novel tissue-mimetic material, in: S.~Ourselin, D.~C. Alexander, C.-F.
  Westin, M.~J. Cardoso (Eds.), Information Processing in Medical Imaging,
  Springer International Publishing, Cham, 2015, pp. 179--190.

\bibitem{Reynaud2017}
O.~Reynaud,
  \href{https://www.frontiersin.org/article/10.3389/fphy.2017.00058}{Time-dependent
  diffusion {MRI} in cancer: Tissue modeling and applications}, Frontiers in
  Physics 5 (2017) 58.
\newblock \href {https://doi.org/10.3389/fphy.2017.00058}
  {\path{doi:10.3389/fphy.2017.00058}}.
\newline\urlprefix\url{https://www.frontiersin.org/article/10.3389/fphy.2017.00058}

\bibitem{Fieremans2011}
E.~Fieremans, J.~H. Jensen, J.~A. Helpern, White matter characterization with
  diffusional kurtosis imaging, NeuroImage 58~(1) (2011) 177 -- 188.

\bibitem{Panagiotaki2012}
E.~Panagiotaki, T.~Schneider, B.~Siow, M.~G. Hall, M.~F. Lythgoe, D.~C.
  Alexander,
  \href{http://www.sciencedirect.com/science/article/pii/S1053811911011566}{Compartment
  models of the diffusion {MR} signal in brain white matter: A taxonomy and
  comparison}, NeuroImage 59~(3) (2012) 2241--2254.
\newline\urlprefix\url{http://www.sciencedirect.com/science/article/pii/S1053811911011566}

\bibitem{Jespersen2007}
S.~N. Jespersen, C.~D. Kroenke, L.~Astergaard, J.~J. Ackerman, D.~A.
  Yablonskiy,
  \href{http://www.sciencedirect.com/science/article/pii/S1053811906010950}{Modeling
  dendrite density from magnetic resonance diffusion measurements}, NeuroImage
  34~(4) (2007) 1473--1486.
\newline\urlprefix\url{http://www.sciencedirect.com/science/article/pii/S1053811906010950}

\bibitem{Ianus2016}
A.~Ianu{\c{s}}, D.~C. Alexander, I.~Drobnjak, Microstructure imaging sequence
  simulation toolbox, in: S.~A. Tsaftaris, A.~Gooya, A.~F. Frangi, J.~L. Prince
  (Eds.), Simulation and Synthesis in Medical Imaging, Springer International
  Publishing, Cham, 2016, pp. 34--44.

\bibitem{Drobnjak2011}
I.~Drobnjak, H.~Zhang, M.~G. Hall, D.~C. Alexander,
  \href{http://www.sciencedirect.com/science/article/pii/S1090780711000838}{The
  matrix formalism for generalised gradients with time-varying orientation in
  diffusion {NMR}}, Journal of Magnetic Resonance 210~(1) (2011) 151 -- 157.
\newblock \href {https://doi.org/10.1016/j.jmr.2011.02.022}
  {\path{doi:10.1016/j.jmr.2011.02.022}}.
\newline\urlprefix\url{http://www.sciencedirect.com/science/article/pii/S1090780711000838}

\bibitem{Mercredi2018}
M.~Mercredi, M.~Martin, \href{https://doi.org/10.1007/s10334-018-0680-1}{Toward
  faster inference of micron-scale axon diameters using {Monte} {Carlo}
  simulations}, Magnetic Resonance Materials in Physics, Biology and Medicine
  31~(4) (2018) 511--530.
\newblock \href {https://doi.org/10.1007/s10334-018-0680-1}
  {\path{doi:10.1007/s10334-018-0680-1}}.
\newline\urlprefix\url{https://doi.org/10.1007/s10334-018-0680-1}

\bibitem{Rensonnet2018}
G.~Rensonnet, B.~Scherrer, S.~K. Warfield, B.~Macq, M.~Taquet,
  \href{https://onlinelibrary.wiley.com/doi/abs/10.1002/mrm.26832}{Assessing
  the validity of the approximation of diffusion-weighted-{MRI} signals from
  crossing fascicles by sums of signals from single fascicles}, Magnetic
  Resonance in Medicine 79~(4) (2018) 2332--2345.
\newblock \href
  {http://arxiv.org/abs/https://onlinelibrary.wiley.com/doi/pdf/10.1002/mrm.26832}
  {\path{arXiv:https://onlinelibrary.wiley.com/doi/pdf/10.1002/mrm.26832}},
  \href {https://doi.org/10.1002/mrm.26832} {\path{doi:10.1002/mrm.26832}}.
\newline\urlprefix\url{https://onlinelibrary.wiley.com/doi/abs/10.1002/mrm.26832}

\bibitem{Hughes1995}
B.~D. Hughes, Random walks and random environments, Clarendon Press Oxford ;
  New York, 1995.

\bibitem{Yeh2013}
C.-H. Yeh, B.~Schmitt, D.~Le~Bihan, J.-R. Li-Schlittgen, C.-P. Lin, C.~Poupon,
  Diffusion microscopist simulator: A general monte carlo simulation system for
  diffusion magnetic resonance imaging, PLoS ONE 8~(10) (2013) e76626.
\newblock \href {https://doi.org/10.1371/journal.pone.0076626}
  {\path{doi:10.1371/journal.pone.0076626}}.

\bibitem{Hall2009}
M.~Hall, D.~Alexander, Convergence and parameter choice for {Monte-Carlo}
  simulations of diffusion {MRI}, IEEE Transactions on Medical Imaging 28~(9)
  (2009) 1354 --1364.
\newblock \href {https://doi.org/10.1109/TMI.2009.2015756}
  {\path{doi:10.1109/TMI.2009.2015756}}.

\bibitem{Balls2009}
G.~T. Balls, L.~R. Frank, \href{http://dx.doi.org/10.1002/mrm.22033}{A
  simulation environment for diffusion weighted {MR} experiments in complex
  media}, Magn. Reson. Med. 62~(3) (2009) 771--778.
\newline\urlprefix\url{http://dx.doi.org/10.1002/mrm.22033}

\bibitem{Nguyen2018a}
K.~V. Nguyen, E.~H. Garzon, J.~Valette,
  \href{http://www.sciencedirect.com/science/article/pii/S1090780718302386}{Efficient
  {GPU}-based {M}onte-{C}arlo simulation of diffusion in real astrocytes
  reconstructed from confocal microscopy}, Journal of Magnetic Resonance
  (2018).
\newblock \href {https://doi.org/10.1016/j.jmr.2018.09.013}
  {\path{doi:10.1016/j.jmr.2018.09.013}}.
\newline\urlprefix\url{http://www.sciencedirect.com/science/article/pii/S1090780718302386}

\bibitem{Waudby2011}
C.~A. Waudby, J.~Christodoulou,
  \href{http://www.sciencedirect.com/science/article/pii/S1090780711001376}{{GPU}
  accelerated {M}onte {C}arlo simulation of pulsed-field gradient {NMR}
  experiments}, Journal of Magnetic Resonance 211~(1) (2011) 67 -- 73.
\newblock \href {https://doi.org/https://doi.org/10.1016/j.jmr.2011.04.004}
  {\path{doi:https://doi.org/10.1016/j.jmr.2011.04.004}}.
\newline\urlprefix\url{http://www.sciencedirect.com/science/article/pii/S1090780711001376}

\bibitem{Hagslatt2003}
H.~Hagslatt, B.~Jonsson, M.~Nyden, O.~Soderman,
  \href{http://www.sciencedirect.com/science/article/pii/S1090780702000393}{Predictions
  of pulsed field gradient {NMR} echo-decays for molecules diffusing in various
  restrictive geometries. simulations of diffusion propagators based on a
  finite element method}, Journal of Magnetic Resonance 161~(2) (2003)
  138--147.
\newline\urlprefix\url{http://www.sciencedirect.com/science/article/pii/S1090780702000393}

\bibitem{Loren2005}
N.~Loren, H.~Hagslatt, M.~Nyden, A.-M. Hermansson, Water mobility in
  heterogeneous emulsions determined by a new combination of confocal laser
  scanning microscopy, image analysis, nuclear magnetic resonance diffusometry,
  and finite element method simulation, The Journal of Chemical Physics 122~(2)
  (2005) --.
\newblock \href {https://doi.org/10.1063/1.1830432}
  {\path{doi:10.1063/1.1830432}}.

\bibitem{Moroney2013}
B.~F. Moroney, T.~Stait-Gardner, B.~Ghadirian, N.~N. Yadav, W.~S. Price,
  \href{http://www.sciencedirect.com/science/article/pii/S1090780713001572}{Numerical
  analysis of {NMR} diffusion measurements in the short gradient pulse limit},
  Journal of Magnetic Resonance 234~(0) (2013) 165--175.
\newline\urlprefix\url{http://www.sciencedirect.com/science/article/pii/S1090780713001572}

\bibitem{Xu2007}
J.~Xu, M.~Does, J.~Gore,
  \href{http://view.ncbi.nlm.nih.gov/pubmed/17374905}{Numerical study of water
  diffusion in biological tissues using an improved finite difference method},
  Physics in medicine and biology 52~(7) (Apr. 2007).
\newline\urlprefix\url{http://view.ncbi.nlm.nih.gov/pubmed/17374905}

\bibitem{Li2014}
J.-R. Li, D.~Calhoun, C.~Poupon, D.~L. Bihan,
  \href{http://stacks.iop.org/0031-9155/59/i=2/a=441}{Numerical simulation of
  diffusion {MRI} signals using an adaptive time-stepping method}, Physics in
  Medicine and Biology 59~(2) (2014) 441.
\newline\urlprefix\url{http://stacks.iop.org/0031-9155/59/i=2/a=441}

\bibitem{Nguyen2014}
D.~V. Nguyen, J.-R. Li, D.~Grebenkov, D.~Le~Bihan,
  \href{http://www.sciencedirect.com/science/article/pii/S0021999114000308}{A
  finite elements method to solve the {B}loch-{T}orrey equation applied to
  diffusion magnetic resonance imaging}, Journal of Computational Physics
  263~(0) (2014) 283--302.
\newline\urlprefix\url{http://www.sciencedirect.com/science/article/pii/S0021999114000308}

\bibitem{Beltrachini2015}
L.~Beltrachini, Z.~A. Taylor, A.~F. Frangi,
  \href{http://www.sciencedirect.com/science/article/pii/S1090780715001743}{A
  parametric finite element solution of the generalised {B}loch–{T}orrey
  equation for arbitrary domains}, Journal of Magnetic Resonance 259 (2015) 126
  -- 134.
\newblock \href {https://doi.org/https://doi.org/10.1016/j.jmr.2015.08.008}
  {\path{doi:https://doi.org/10.1016/j.jmr.2015.08.008}}.
\newline\urlprefix\url{http://www.sciencedirect.com/science/article/pii/S1090780715001743}

\bibitem{Russell2012}
G.~Russell, K.~D. Harkins, T.~W. Secomb, J.-P. Galons, T.~P. Trouard,
  \href{http://stacks.iop.org/0031-9155/57/i=4/a=N35}{A finite difference
  method with periodic boundary conditions for simulations of
  diffusion-weighted magnetic resonance experiments in tissue}, Physics in
  Medicine and Biology 57~(4) (2012) N35.
\newline\urlprefix\url{http://stacks.iop.org/0031-9155/57/i=4/a=N35}

\bibitem{Verwer1990}
J.~G. Verwer, W.~H. Hundsdorfer, B.~P. Sommeijer, {Convergence properties of
  the Runge-Kutta-Chebyshev method}, Numerische Mathematik 57 (1990) 157--178.
\newblock \href {https://doi.org/10.1007/BF01386405}
  {\path{doi:10.1007/BF01386405}}.

\bibitem{Nguyen2016a}
V.~D. Nguyen, \href{https://www.eccomas2016.org/}{{A FEniCS-HPC framework for
  multi-compartment Bloch-Torrey models}}, Vol.~1, 2016, pp. 105--119, {QC
  20170509}.
\newline\urlprefix\url{https://www.eccomas2016.org/}

\bibitem{Nguyen2018}
V.-D. Nguyen, J.~Jansson, J.~Hoffman, J.-R. Li,
  \href{http://www.sciencedirect.com/science/article/pii/S0021999118305709}{A
  partition of unity finite element method for computational diffusion {MRI}},
  Journal of Computational Physics (2018).
\newblock \href {https://doi.org/10.1016/j.jcp.2018.08.039}
  {\path{doi:10.1016/j.jcp.2018.08.039}}.
\newline\urlprefix\url{http://www.sciencedirect.com/science/article/pii/S0021999118305709}

\bibitem{Nguyen2019}
V.-D. Nguyen, J.~Jansson, H.~T.~A. Tran, J.~Hoffman, J.-R. Li,
  \href{http://www.sciencedirect.com/science/article/pii/S1090780719300023}{Diffusion
  {MRI} simulation in thin-layer and thin-tube media using a discretization on
  manifolds}, Journal of Magnetic Resonance 299 (2019) 176 -- 187.
\newblock \href {https://doi.org/https://doi.org/10.1016/j.jmr.2019.01.002}
  {\path{doi:https://doi.org/10.1016/j.jmr.2019.01.002}}.
\newline\urlprefix\url{http://www.sciencedirect.com/science/article/pii/S1090780719300023}

\bibitem{Si2015}
H.~Si, \href{http://doi.acm.org/10.1145/2629697}{{TetGen}, a {D}elaunay-based
  quality tetrahedral mesh generator}, ACM Trans. Math. Softw. 41~(2) (2015)
  11:1--11:36.
\newblock \href {https://doi.org/10.1145/2629697} {\path{doi:10.1145/2629697}}.
\newline\urlprefix\url{http://doi.acm.org/10.1145/2629697}

\bibitem{RahmanValdman2013}
T.~Rahman, J.~Valdman,
  \href{http://www.sciencedirect.com/science/article/pii/S0096300311010836}{Fast
  {MATLAB} assembly of {FEM} matrices in {2D} and {3D}: Nodal elements},
  Applied Mathematics and Computation 219~(13) (2013) 7151 -- 7158, eSCO 2010
  Conference in Pilsen, June 21- 25, 2010.
\newblock \href {https://doi.org/https://doi.org/10.1016/j.amc.2011.08.043}
  {\path{doi:https://doi.org/10.1016/j.amc.2011.08.043}}.
\newline\urlprefix\url{http://www.sciencedirect.com/science/article/pii/S0096300311010836}

\bibitem{Nguyen2014d}
D.~V. Nguyen, J.-R. Li, D.~Grebenkov, D.~L. Bihan,
  \href{http://www.sciencedirect.com/science/article/pii/S0021999114000308}{A
  finite elements method to solve the {B}loch–{T}orrey equation applied to
  diffusion magnetic resonance imaging}, Journal of Computational Physics
  263~(0) (2014) 283 -- 302.
\newblock \href {https://doi.org/10.1016/j.jcp.2014.01.009}
  {\path{doi:10.1016/j.jcp.2014.01.009}}.
\newline\urlprefix\url{http://www.sciencedirect.com/science/article/pii/S0021999114000308}

\bibitem{Callaghan1995}
P.~T. Callaghan, J.~Stepianik,
  \href{http://www.sciencedirect.com/science/article/pii/S1064185885799597}{Frequency-domain
  analysis of spin motion using modulated-gradient {NMR}}, Journal of Magnetic
  Resonance, Series A 117~(1) (1995) 118--122.
\newline\urlprefix\url{http://www.sciencedirect.com/science/article/pii/S1064185885799597}

\bibitem{schiavi2016}
H.~Haddar, J.-R. Li, S.~Schiavi, \href{https://doi.org/10.1137/15M1019398}{A
  macroscopic model for the diffusion {MRI} signal accounting for
  time-dependent diffusivity}, SIAM Journal on Applied Mathematics 76~(3)
  (2016) 930--949.
\newblock \href {https://doi.org/10.1137/15M1019398}
  {\path{doi:10.1137/15M1019398}}.
\newline\urlprefix\url{https://doi.org/10.1137/15M1019398}

\bibitem{Mitra1992}
P.~P. Mitra, P.~N. Sen, L.~M. Schwartz, P.~Le~Doussal, Diffusion propagator as
  a probe of the structure of porous media, Physical review letters 68~(24)
  (1992) 3555--3558.

\bibitem{Mitra1993}
P.~P. Mitra, P.~N. Sen, L.~M. Schwartz, Short-time behavior of the diffusion
  coefficient as a geometrical probe of porous media, Phys. Rev. B 47 (1993)
  8565--8574.

\bibitem{Callaghan1997}
P.~Callaghan, \href{http://dx.doi.org/10.1006/jmre.1997.1233}{A simple matrix
  formalism for spin echo analysis of restricted diffusion under generalized
  gradient waveforms}, Journal of Magnetic Resonance 129~(1) (1997) 74--84.
\newline\urlprefix\url{http://dx.doi.org/10.1006/jmre.1997.1233}

\bibitem{Barzykin1999}
A.~V. Barzykin,
  \href{http://www.sciencedirect.com/science/article/pii/S1090780799917780}{Theory
  of spin echo in restricted geometries under a step-wise gradient pulse
  sequence}, Journal of Magnetic Resonance 139~(2) (1999) 342--353.
\newline\urlprefix\url{http://www.sciencedirect.com/science/article/pii/S1090780799917780}

\bibitem{Grebenkov2007}
D.~Grebenkov, {NMR} survey of reflected brownian motion, Reviews of Modern
  Physics 79~(3) (2007) 1077--1137.
\newblock \href {https://doi.org/10.1103/RevModPhys.79.1077}
  {\path{doi:10.1103/RevModPhys.79.1077}}.

\bibitem{Ozarslan2009}
E.~Ozarslan, N.~Shemesh, P.~J. Basser,
  \href{https://www.ncbi.nlm.nih.gov/pmc/PMC2736571/}{A general framework to
  quantify the effect of restricted diffusion on the {NMR} signal with
  applications to double pulsed field gradient {NMR} experiments}, The Journal
  of chemical physics 130~(19292544) (2009) 104702--104702.
\newline\urlprefix\url{https://www.ncbi.nlm.nih.gov/pmc/PMC2736571/}

\bibitem{Drobnjak2011a}
I.~Drobnjak, H.~Zhang, M.~G. Hall, D.~C. Alexander,
  \href{http://www.sciencedirect.com/science/article/pii/S1090780711000838}{The
  matrix formalism for generalised gradients with time-varying orientation in
  diffusion {NMR}}, Journal of Magnetic Resonance 210~(1) (2011) 151--157.
\newline\urlprefix\url{http://www.sciencedirect.com/science/article/pii/S1090780711000838}

\bibitem{Grebenkov2010a}
D.~S. Grebenkov,
  \href{http://www.sciencedirect.com/science/article/pii/S1090780710001199}{Pulsed-gradient
  spin-echo monitoring of restricted diffusion in multilayered structures},
  Journal of Magnetic Resonance 205~(2) (2010) 181--195.
\newline\urlprefix\url{http://www.sciencedirect.com/science/article/pii/S1090780710001199}

\bibitem{Ascoli2007}
G.~A. Ascoli, D.~E. Donohue, M.~Halavi,
  \href{http://www.jneurosci.org/content/27/35/9247}{Neuromorpho.org: A central
  resource for neuronal morphologies}, Journal of Neuroscience 27~(35) (2007)
  9247--9251.
\newblock \href
  {http://arxiv.org/abs/http://www.jneurosci.org/content/27/35/9247.full.pdf}
  {\path{arXiv:http://www.jneurosci.org/content/27/35/9247.full.pdf}}, \href
  {https://doi.org/10.1523/JNEUROSCI.2055-07.2007}
  {\path{doi:10.1523/JNEUROSCI.2055-07.2007}}.
\newline\urlprefix\url{http://www.jneurosci.org/content/27/35/9247}

\bibitem{Geuzaine2009}
C.~Geuzaine, J.~F. Remacle, {Gmsh: a three-dimensional finite element mesh
  generator with built-in pre- and post-processing facilities}, International
  Journal for Numerical Methods in Engineering 79~(11) (2009) 1309--1331.

\bibitem{FIEREMANS201839}
E.~Fieremans, H.-H. Lee,
  \href{http://www.sciencedirect.com/science/article/pii/S1053811918305536}{Physical
  and numerical phantoms for the validation of brain microstructural {MRI}: A
  cookbook}, NeuroImage 182 (2018) 39 -- 61, microstructural Imaging.
\newblock \href
  {https://doi.org/https://doi.org/10.1016/j.neuroimage.2018.06.046}
  {\path{doi:https://doi.org/10.1016/j.neuroimage.2018.06.046}}.
\newline\urlprefix\url{http://www.sciencedirect.com/science/article/pii/S1053811918305536}

\bibitem{Shemesh2016}
N.~Shemesh, S.~N. Jespersen, D.~C. Alexander, Y.~Cohen, I.~Drobnjak, T.~B.
  Dyrby, J.~Finsterbusch, M.~A. Koch, T.~Kuder, F.~Laun, M.~Lawrenz,
  H.~Lundell, P.~P. Mitra, M.~Nilsson, E.~Özarslan, D.~Topgaard, C.-F. Westin,
  \href{https://onlinelibrary.wiley.com/doi/abs/10.1002/mrm.25901}{Conventions
  and nomenclature for double diffusion encoding {NMR} and {MRI}}, Magnetic
  Resonance in Medicine 75~(1) (2016) 82--87.
\newblock \href {https://doi.org/10.1002/mrm.25901}
  {\path{doi:10.1002/mrm.25901}}.
\newline\urlprefix\url{https://onlinelibrary.wiley.com/doi/abs/10.1002/mrm.25901}

\bibitem{Dhital2019}
B.~Dhital, M.~Reisert, E.~Kellner, V.~G. Kiselev,
  \href{http://www.sciencedirect.com/science/article/pii/S1053811919300151}{Intra-axonal
  diffusivity in brain white matter}, NeuroImage 189 (2019) 543 -- 550.
\newblock \href
  {https://doi.org/https://doi.org/10.1016/j.neuroimage.2019.01.015}
  {\path{doi:https://doi.org/10.1016/j.neuroimage.2019.01.015}}.
\newline\urlprefix\url{http://www.sciencedirect.com/science/article/pii/S1053811919300151}

\bibitem{Novikov2019}
D.~S. Novikov, E.~Fieremans, S.~N. Jespersen, V.~G. Kiselev,
  \href{https://onlinelibrary.wiley.com/doi/abs/10.1002/nbm.3998}{Quantifying
  brain microstructure with diffusion {MRI}: Theory and parameter estimation},
  NMR in Biomedicine 32~(4) (2019) e3998.
\newblock \href {https://doi.org/10.1002/nbm.3998}
  {\path{doi:10.1002/nbm.3998}}.
\newline\urlprefix\url{https://onlinelibrary.wiley.com/doi/abs/10.1002/nbm.3998}

\bibitem{Henriques2019}
R.~N. Henriques, S.~N. Jespersen, N.~Shemesh,
  \href{https://onlinelibrary.wiley.com/doi/abs/10.1002/mrm.27606}{Microscopic
  anisotropy misestimation in spherical-mean single diffusion encoding {MRI}},
  Magnetic Resonance in Medicine 81~(5) (2019) 3245--3261.
\newblock \href
  {http://arxiv.org/abs/https://onlinelibrary.wiley.com/doi/pdf/10.1002/mrm.27606}
  {\path{arXiv:https://onlinelibrary.wiley.com/doi/pdf/10.1002/mrm.27606}},
  \href {https://doi.org/10.1002/mrm.27606} {\path{doi:10.1002/mrm.27606}}.
\newline\urlprefix\url{https://onlinelibrary.wiley.com/doi/abs/10.1002/mrm.27606}

\bibitem{Topgaard2017}
D.~Topgaard,
  \href{http://www.sciencedirect.com/science/article/pii/S1090780716302701}{Multidimensional
  diffusion {MRI}}, Journal of Magnetic Resonance 275 (2017) 98 -- 113.
\newblock \href {https://doi.org/https://doi.org/10.1016/j.jmr.2016.12.007}
  {\path{doi:https://doi.org/10.1016/j.jmr.2016.12.007}}.
\newline\urlprefix\url{http://www.sciencedirect.com/science/article/pii/S1090780716302701}

\bibitem{Veraart2018}
J.~Veraart, D.~S. Novikov, E.~Fieremans,
  \href{http://www.sciencedirect.com/science/article/pii/S1053811917307784}{{TE}
  dependent diffusion imaging ({TEdDI}) distinguishes between compartmental
  {T2} relaxation times}, NeuroImage 182 (2018) 360 -- 369, microstructural
  Imaging.
\newblock \href
  {https://doi.org/https://doi.org/10.1016/j.neuroimage.2017.09.030}
  {\path{doi:https://doi.org/10.1016/j.neuroimage.2017.09.030}}.
\newline\urlprefix\url{http://www.sciencedirect.com/science/article/pii/S1053811917307784}

\bibitem{Lampinen2019}
B.~Lampinen, F.~Szczepankiewicz, M.~Novén, D.~van Westen, O.~Hansson,
  E.~Englund, J.~Mårtensson, C.-F. Westin, M.~Nilsson,
  \href{https://onlinelibrary.wiley.com/doi/abs/10.1002/hbm.24542}{Searching
  for the neurite density with diffusion {MRI}: Challenges for biophysical
  modeling}, Human Brain Mapping 40~(8) (2019) 2529--2545.
\newblock \href {https://doi.org/10.1002/hbm.24542}
  {\path{doi:10.1002/hbm.24542}}.
\newline\urlprefix\url{https://onlinelibrary.wiley.com/doi/abs/10.1002/hbm.24542}

\bibitem{Ginsburger2019}
K.~Ginsburger, F.~Matuschke, F.~Poupon, J.-F. Mangin, M.~Axer, C.~Poupon,
  \href{http://www.sciencedirect.com/science/article/pii/S105381191930151X}{{MEDUSA}:
  A gpu-based tool to create realistic phantoms of the brain microstructure
  using tiny spheres}, NeuroImage 193 (2019) 10 -- 24.
\newblock \href
  {https://doi.org/https://doi.org/10.1016/j.neuroimage.2019.02.055}
  {\path{doi:https://doi.org/10.1016/j.neuroimage.2019.02.055}}.
\newline\urlprefix\url{http://www.sciencedirect.com/science/article/pii/S105381191930151X}

\bibitem{Rensonnet2019}
G.~Rensonnet, B.~Scherrer, G.~Girard, A.~Jankovski, S.~K. Warfield, B.~Macq,
  J.-P. Thiran, M.~Taquet,
  \href{http://www.sciencedirect.com/science/article/pii/S1053811918319487}{Towards
  microstructure fingerprinting: Estimation of tissue properties from a
  dictionary of monte carlo diffusion {MRI} simulations}, NeuroImage 184 (2019)
  964 -- 980.
\newblock \href
  {https://doi.org/https://doi.org/10.1016/j.neuroimage.2018.09.076}
  {\path{doi:https://doi.org/10.1016/j.neuroimage.2018.09.076}}.
\newline\urlprefix\url{http://www.sciencedirect.com/science/article/pii/S1053811918319487}

\bibitem{Palombo2019}
M.~Palombo, D.~C. Alexander, H.~Zhang,
  \href{http://www.sciencedirect.com/science/article/pii/S1053811918321694}{A
  generative model of realistic brain cells with application to numerical
  simulation of the diffusion-weighted {MR} signal}, NeuroImage 188 (2019) 391
  -- 402.
\newblock \href
  {https://doi.org/https://doi.org/10.1016/j.neuroimage.2018.12.025}
  {\path{doi:https://doi.org/10.1016/j.neuroimage.2018.12.025}}.
\newline\urlprefix\url{http://www.sciencedirect.com/science/article/pii/S1053811918321694}

\end{thebibliography}






\end{document}